\documentclass{amsart}
\usepackage{xr}
\newcommand*\sref[1]{%
    \ref{#1} of [1]}

\usepackage{pstricks}
\def\FigureRho{%
\psset{unit = 1cm}
\pspicture(4,4)
\pscircle[linestyle=dotted, linecolor=red](2,2.5){1.3}
\psline[linestyle=dotted, linecolor=red](2,0)(1.54,1.2)
\pscircle(1.6,1.27){0.04}
\put(1.9,1.75){\n}
\put(1.0,0.75){\k}
\put(2.1,-0.1){$0$}
\psline{->}(1.22,0.92)(1.5,1.17)
\psline{->}(1.9,1.7)(1.63,1.33)
\psline{->}(0.4,1.3)(1.41,1.3)
\put(1.9,0.6){$\Stem$}
\put(3.5,2.5){$\L$}
\put(-1,1.2){$\s \n = \s \k$}
\endpspicture}

\begin{document}

\title{The Church numbers in NF Set Theory}     
\author{Michael Beeson}   
\date{\today}

\newtheorem{theorem}{Theorem}
\newtheorem{lemma}[theorem]{Lemma}
\numberwithin{theorem}{section} 
\newtheorem{corollary}[theorem]{Corollary}
\newtheorem{definition}[theorem]{Definition}
\let\sect\S
\def\s{{\sf S}}
\def\prog{{\sf prog}}
\def\bbS{{\mathbb S\,}}
\def\k{{\bf k}}
\def\n{{\bf n}}
\def\m{{\bf m}}
\def\K{{\bf K}}
\def\Stem{{\sf STEM}}
\def\L{{\mathcal L}}
\def\R{{\mathcal R}}
\def\Z{{\mathbb Z}}
\def\d{{\bf d}}
\def\C{{\mathcal C}}
\def\r{{\bf \ r\, }}
\def\q{{\bf q}}
\def\App{{\bf App}}
\def\lgn{{}^{\lceil}}   
\def\rgn{{}^{\rceil}}   
\def\defined{{\downarrow}}
\def\imp{\ {\rightarrow}\ }
\def\iff{\ {\leftrightarrow}\ } 
\def\FUNC{{\sf FUNC}}
\def\V{{\mathbb V}}
\def\N{{\mathbb N}}
\def\S{{\mathbb S}}
\def\T{{\mathbb T}}
\def\J{{\mathbb J}} 
\def\JC{{\mathcal J}}
\def\JLift{{\mathbb B}}
\def\F{{\mathbb F}}
\def\FregeN{{\mathbb F}}
\def\i{{\bf i\, }}
\def\j{{\bf j\,}}
\def\comp{{\mathit comp\,}} 
\def\INF{{\rm INF}}
\def\FINITE{{\sf FINITE}}
\def\PRIME{{\sf PRIME}}
\def\tail{{\sf tail}}
\def\cl{{\sf cl}}
\def\size{{\sf size}}
\def\order{{\sf order}}
\def\Order{{\sf Order}}
\def\empty{{\Lambda}}
\def\ox{{\,\otimes\,}}
\def\DECIDABLE{{\sf DECIDABLE}}
\def\zero{{\sf zero}}
\def\one{{\sf one}} 
\def\two{{\sf two}}
\def\three{{\sf three}}
\def\id{{\sf id}}
\def\max{{\sf max}}
\def\div{{\ | \ }}
\def\imp{\ {\rightarrow}\ }
\def\ChurchZero{{\bf 0}}
\def\twocases#1#2{\left\{
\begin{array}{ll} 
#1\\ 
{}\\
#2
\end{array} \right.}  
\def\threecases#1#2#3{\left\{
\begin{array}{ll} 
#1\\ 
#2\\
#3
\end{array} \right.}  

\begin{abstract} By NF we mean Quine's New Foundations set theory.
 We define the Church numerals (or better, Church numbers)
and elaborate their properties in INF. 
Here we investigate the question whether the set of Church numbers
is infinite.  Usually in NF, the natural numbers are represented by 
the finite Frege cardinals.  We answer the question by proving that 
if the set of finite Frege cardinals is infinite, then so is the
set of Church numbers.  

 Specker showed in 1953 that classical NF proves the
 set of finite Frege cardinals is infinite, so using classical
 logic the set of Church numbers is infinite.  It has long been
 an open problem whether any set can be proved infinite using 
NF with intuitionistic logic (INF).  Perhaps INF proves the set
of Church numbers is infinite; we tried to prove that, 
but we could only succeed with an additional assumption,
the ``Church counting axiom.''  That is
a fundamental counting principle:
it says
that iterating successor $n$ times, starting at zero, results in $n$.

  We also prove, without the 
aid of the counting axiom, that if the set of Church numbers is not finite,
then it is infinite, and Church successor is one-to-one.  Consequently,
Heyting's arithmetic is interpretable in INF plus the Church counting axiom.  
Finally, we show that the Church counting axiom is equivalent in INF to Rosser's
 counting axiom.   That equivalence is a new theorem even classically. 
 Since it is known that the Rosser counting axiom is not 
provable in NF (if NF is consistent),  it follows that the same is true of the 
Church counting axiom. 
The original question, whether INF proves the set of Church numbers is infinite,
remains open.  At least it has been shown equivalent to the question whether
INF proves there exists an infinite set.   
\end{abstract}

\maketitle

\begin{quotation}
\vskip-0.6cm
\hskip 2.7cm
\raggedleft
{\em In all the world there is nothing so interesting, so curious, and so beautiful as truth.}

--Hercule Poirot
\end{quotation}

\tableofcontents

\section{Introduction}
Quine's NF set theory is a first-order theory whose language contains only the binary predicate symbol $\in$,
and whose axioms are two in number:  extensionality and stratified comprehension.  
The definition of these axioms will be reviewed below; full details 
can be found in \cite{forster}.
Intuitionistic NF, or \INF,  is the theory
with the same language and axioms as NF, but with intuitionistic logic instead of classical.   

The ``axiom'' of infinity is a theorem of NF, proved by Rosser 
\cite{rosser1952} and Specker \cite{specker1953};
see also \cite{forster}, p.~49.  These proofs use classical logic
in an apparently essential way, and it is still an open question 
whether \INF\ proves the existence of an infinite set.  
In \cite{inf-basics},  the fundamental results of \INF\ have been 
developed,  and on that basis, we analyzed (in unpublished work)
the constructive content
of Specker's proof,  but that did not lead to a proof of infinity. 
It has also long been known that if $\V$ is not finite, then the logic
of stratified formulas must be classical, so Specker's proof could be 
done; hence the statement ``$\V$ is not finite'' can be proved in \INF.
But that does not lead to a proof of infinity, because there might be 
a maximum integer,  whose elements $U$ would be ``unenlargeable'' 
in the sense that we cannot find any set that is not in $U$.

We therefore have two reasons for being interested in the Church
numerals in \INF:  to use them to prove infinity,  and to study their 
structure for its own intrinsic interest.

   The Church numerals were introduced by Church (or technically, 
by his student Kleene) in the context of $\lambda$-calculus.  They are defined so that the $n$-th Church number is a function 
that takes inputs $f$ and $x$ and produces output $f^n(x)$, where $f^n$
is $f$ iterated $n$ times.  It turns out that the Church numerals,
or Church numbers as we shall call them%
\footnote{
The word ``numeral'' is usually used for a syntactic object, a name 
for a number.  The ``Church numerals'' that we define here are 
sets, not syntax.  We therefore refer to them as ``Church numbers''
instead.  Church himself in \cite{church1941} 
never mentions the word ``numeral'' but instead refers to ``the formula
representing the integer $n$''.

},  
and the set $\N$ of Church numbers, can be defined
by stratified comprehension straightforwardly in \INF. 
The Church numerals were introduced in 1935 \cite{kleene1935} 
and NF set theory in 1937 \cite{quine1937}, so the definition 
of Church numerals in NF could have been done at any time since 1937.%
\footnote{Church defines addition, multiplication, and exponentiation
on p.~10 of \cite{church1941}, but he attributes these definitions
to Rosser and refers to \cite{kleene1935} for detailed proofs. 
Rosser later worked on Quine's NF, so he had all the background
needed to write this paper 
right after he wrote \cite{rosser1952}. 
}  

The equation satisfied by Church successor $\s$ is 
$$ \s jfx = f(jfx). $$
 It was Church's student Kleene who defined
the predecessor function in $\lambda$-calculus, in his Ph.~D. thesis,
and in  \cite{kleene1935}.
The $\lambda$-calculus definition of that function does not lead
to a definition by stratified comprehension in NF, so to prove that
Church successor is injective in NF requires a new argument.

\begin{figure}[ht]
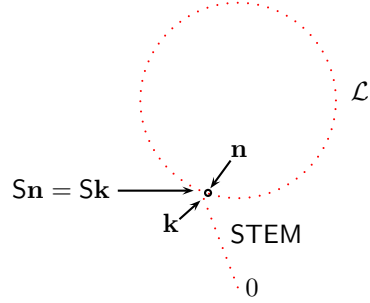

\caption{The stem $\Stem$, the loop $\L$, and unique double successor}
\label{figure:rho}
\FigureRho
\end{figure}

 In this paper we analyze the structure
of the Church numbers under successor.  First we 
introduce addition $x \oplus y$ on $\N$. Assuming $\N$ is finite,
there is exactly one ``double successor'', i.e., there are exactly
two numerals $\k$ and $\n$ with $\s \k = \s \n$.  The numbers from
0 to $\k$ behave ``normally'' (trichotomy holds there).  We 
call that the ``stem.''
The rest of the Church numbers form a ``loop'' $\L$, as a directed
graph in which the edges are pairs $\langle x, \s x\rangle$. 
See Fig.~\ref{figure:rho}.  

So trichotomy fails dramatically on $\L$: everything in $\L$ is 
less than everything else.  The ``length'' of the loop is a 
Church number $\m$ such that $x\oplus \m = x$ for every $x \in \L$.
We then prove the Annihilation Theorem:  any one-to-one map $f$ on 
any finite set,  when iterated $\m$ times, is the identity.
In symbols, $\m f = f$.  Here the $\m$-fold iteration of $f$
is expressed by the Church number $\m$. 

One could then reach a contradiction if one could simply exhibit
a finite set $X$ and a map whose $\m$-th iterate is not the identity.
The crucial question here is, what is the order of successor, considered
as a permutation of the loop?   Intuitively it seems that it should 
be $\m$, the length of the loop.  That assertion,  however,  is 
equivalent to the unstratified formula 
$$ j \in \N \imp  j \s \n = \n \oplus j,$$
which in turn is equivalent to the Church counting axiom, 
$$ j \in \N \imp j \s \ChurchZero = j.$$
That unstratified formula,  although expressing a fundamental truth 
about the Church numbers,  seems to be unprovable in NF.  
If we assume the Church counting axiom, then it is a routine exercise
to define a permutation of the loop that  leaves one element unmoved
while cycling the rest, and that permutation will have order $\m-1$,
and hence will not satisfy the Annihilation Theorem.  That will show 
that $\N$ is not finite; and we can go further and prove that $\N$ is 
actually infinite, and Heyting's arithmetic is thus interpretable in 
INF plus the counting axiom. 

Without the counting axiom,   one can exhibit 
specific permutations of specific finite sets, for example, a permutation
of order 3 or one of order 873  or one of order 24,569; so $\m$ 
must be larger than any integer with a name.  In general we must think
of the Church numbers as non-standard integers.  The counting axiom 
is a way to outlaw certain non-standard integers.   In general in NF set theory 
one cannot assert the existence of $\{ x \in \N :  P(x)\}$ unless $P$ is 
stratified,  though such sets do not lead to any known paradoxes.  Of course
this separation axiom would imply the counting axiom.  Without the counting 
axiom, the order $\q$ of the loop might be smaller than $\m$;  in fact the 
quotient $\m/\q$ would itself be a non-standard integer, and the map 
$x \mapsto x\s\n$ (which is not definable in NF because the formula is 
not stratified)  would  take a small initial part of the loop 
onto the whole loop,  and as $x$ increases past $q$,  $x \s \n$ would
wrap around the loop many times before $x$ reached $\m$.  It would 
be interesting to see models of NF in which the Church counting axiom fails;
but then again,  it would be interesting to see models of NF at all.

Previous studies of infinity in NF have used the ``Frege cardinals''  $\F$,
defined as the least set containing $\zero := \{\empty\}$ and closed under
inhabited successor, where the Frege successor of $x$ is $x^+ $, the set of all $x \cup\{c\}$ 
such that $c \not\in x$.%
\footnote{Actually, Rosser and Spector used Nn, defined as the least set 
containing $\zero$ and closed under successor;  so possibly Nn might contain 
$\empty$, while $\F$ certainly does not contain $\empty$.  Therefore we 
use a different letter.}
Of course there are natural functions defined from $\N$
to $\F$ and vice-versa,  but nothing about these functions is obvious.  If $\F$ is 
finite, then there is a maximum integer $\max$,  containing a set $U$ that is 
``unenlargeable'' in the sense that no $c$ can be produced such that $c \not\in U$,
in spite of the fact that the universe $\V$ is infinite.   In that case there are 
several possibilities:  $\max$ might come ``too soon'',  so that there are ``not enough''
Frege numbers to correspond to all the Church numbers.  Or, $\max$ might be ``too big'',
so that the correspondence is no longer one-to-one,  as the Frege numbers start to wrap 
the ``loop'' a second time.  Possibly $\max$ is ``just right'' and the Frege numbers
correspond to the Church numbers in one-to-one fashion.    If $\N$ and $\F$
are both finite,  there is no obvious relation between their cardinalities.
\smallskip

If $\N$ is infinite, then of course $\F$ is infinite too, since there are 
arbitrarily large finite sets, namely initial segments of $\N$. 
  On the other hand, if $\F$ is infinite, then we can prove that 
$\N$ is infinite too,  since if $\F$ is infinite, we can show that every Church
number is the order of some permutation,  but that contradicts the 
Annihilation Theorem.  Taken together with Specker's result that
 $\F$ is infinite in classical NF,  we see that classical NF does prove
 that $\N$ is infinite.    
\smallskip

In the final section of the paper, we prove that the Church counting
axiom and the Rosser counting axiom are equivalent. The proof uses
our main result that the Church counting axiom implies $\N$ is infinite
and Church successor is one-to-one.   This equivalence is proved
in INF, without any assumptions.   

We refer to \cite{inf-basics} for notation, axioms, and the basic theorems of 
\INF,  including the properties of finite sets and  finite cardinals.
In particular we use $\F$ for the set of finite (Frege) 
cardinals.  In this paper
we make no use of arithmetic on the Frege cardinals, not even addition, let
alone multiplication and exponentiation.  What we mostly require from 
\cite{inf-basics} is the intuitionistic theory of finite sets.    Lemmas,
theorems, and definitions from that \cite{inf-basics} will be referenced  
like this:  Lemma~\sref{lemma:finitedecidable}.

Thanks to Thomas Forster for asking me (once a year for twenty years)
about the strength of \INF.  
Thanks to Randall Holmes for his attention to my first draft,
and for the idea of the proof of Lemma~\ref{lemma:additionfunction}. 
Thanks to Albert Visser for his 
careful reading of an earlier version; many errors were
thus corrected.   Thanks to the creators of the proof assistant Lean \cite{LeanSystem},  which 
has enabled me to state with high confidence that there are no errors in this 
paper.  Thanks to the users of Lean who helped me acquire sufficient 
expertise in using Lean by answering my questions, especially Mario Carneiro.%
\footnote{Although the proofs in this paper have been computer-checked 
for correctness using Lean,  they are presented here in human-readable form
with detailed proofs.  Issues concerning the notation in Lean will not be discussed
here.}

There are many lemmas in this paper, and the intention is that each 
of those lemmas is provable in \INF.  Inductions are stratified and proofs are 
intuitionistically valid.  When the counting axiom is used, it is explicitly 
mentioned as a hypothesis.

\section{The Church numbers} 
We  define the class of single-valued relations:
\begin{definition} \label{definition:FUNC} 
$$ \FUNC = \{ f : \forall x,y,z\, (\langle x,y \rangle \in f \imp \langle x,z\rangle \in f \imp y = z)\}.$$
\end{definition}
The definition does not rule out the possibility that $f$ might contain some members that are not ordered pairs.
A function is a single-valued relation, i.e. $f \in \FUNC \ \land \ Rel(f)$.
\smallskip

If   $f \in \FUNC$  and $\langle x,y \rangle \in f$, then
informally we write $y = f(x)$.  Formally this is $y = Ap(f,x)$, where $Ap$ is a 
function symbol defined using stratified comprehension.  For details about $Ap$ see
Definition~\sref{definition:Ap} and Lemma~\sref{lemma:Ap}.
We will often suppress mention of the symbol $Ap$, as there is no other way to interpret $f(x)$.  
In fact, we will informally follow the $\lambda$-calculus convention of writing $(f x)$ or just $fx$
for function application, with association to the left, so $fxy$ means $((f x) y)$. 

The Church successor 
function $\s $ is defined in $\lambda$-calculus by
\begin{eqnarray*}
\s (z) &=& \lambda f \lambda x f(zfx) 
\end{eqnarray*}
Imitating this definition in NF we wish we could define
\begin{eqnarray*}
\s  &=& \{ \langle z,\{ \langle f, \{ \langle x, f(zfx)\rangle \} \rangle\}\rangle\}
\end{eqnarray*}
Expanding the formula on the right, it is equivalent, at least 
for functions $f$ and $z$, to  the formula in the following definition.
To explain the relation between the two formulas,
$t = zf$, $q = tx = zfx$, $w = fq = f(zfx)$.
But we emphasize, everything about Church successor up to this 
point is merely motivation for the definition below.  

\begin{definition} \label{definition:sucdef}
Church successor is defined by
\begin{eqnarray*}
&& \s = \{ \langle z, \{ \langle f,p \rangle : f \in \FUNC \ \land \\
&& \forall u\, (u \in p \iff \exists x,w,t,q \,(u = \langle x,w\rangle \ \land \ t \in \FUNC \ \land \\
&& \langle f,t\rangle \in z \ \land \ \langle x,q \rangle \in t \ \land \ \langle q,w \rangle \in f)) \}
         : z \in \FUNC \}
\end{eqnarray*}
\end{definition}

\begin{lemma} \label{lemma:sucdef} The definition of  Church successor can be given in \INF\ using
stratified comprehension; that is, the graph of Church successor is definable in \INF. 
\end{lemma}

\noindent{\em Proof}.
To stratify the formula in Definition~\ref{definition:sucdef}, 
we assign indices as follows:

\begin{center}
\begin{tabular}{rr}
$x$, $q$, and $w$  & 0 \\
$u, \langle x,q \rangle$, and $\langle q,w \rangle$ & 2 \\
$p$, $t$ and  $f$ & 3 \\
$\langle f,t \rangle$ & 5 \\
$z$             & 6  
\end{tabular}
\end{center}

\noindent
With this assignment,  the left and right members of ordered pairs get the same index, and the formula is stratified.
$t \in \FUNC$ can be stratified assigning $t$ any desired index $\ge 2$, and that condition is satisfied by 
the assignments in the table.
That completes the proof. 

\begin{lemma}\label{lemma:ChurchSuccessor3} Let $f \in \FUNC$ and 
$z \in \FUNC$.  Then $\s zf$ is a relation (contains only ordered pairs).
\end{lemma}

\noindent{\em Proof}.  Suppose $t \in \s zf$.  We have to prove $t$ is an ordered pair.
Officially $\s z f$ is $Ap(\s z, f)$.  By definition of $Ap$, there exists $y$ such that
$\langle f,y \rangle \in \s z$ and $t \in y$.  Then by definition of Church successor,
$z \in \FUNC$ and $f \in \FUNC$ and for every $u$,
$$ u \in y \iff \exists x,w,r,q\, (u = \langle x,w\rangle \ \land\ r \in \FUNC \ \land \ 
\langle f,r \rangle \in z \ \land \ \langle x,q \rangle \in r \ \land \ \langle q,w\rangle \in f).$$
Instantiate the quantified $u$ to $t$; then $t = \langle x,w \rangle$ for some $x$ and $w$.
That completes the proof of the lemma. 

\begin{lemma} \label{lemma:ChurchSuccessor2} Let $f \in \FUNC$ and 
$z \in \FUNC$.  Then
\begin{eqnarray*}
\s zf &=&   \{ \langle x, w\rangle :
             \exists t,q\,
                     ( t \in \FUNC \land \langle f,t \rangle \in z \land
                       \langle x,q \rangle \in t \land
                       \langle q,w \rangle \in f 
                      )
                \} 
 \end{eqnarray*}
 \end{lemma}
 
 \noindent{\em Proof}.  Let $f \in \FUNC$ and $z \in \FUNC$.  By Lemma~\ref{lemma:ChurchSuccessor3},
 it suffices to prove 
 \begin{eqnarray*}
  \langle x,w \rangle \in Ap(\s z, f) \iff 
  \exists t,q\,(t \in \FUNC \ \land \ \langle f,t \rangle \in z \ \land \ \langle x,q \rangle \in t \ \land \ \langle q,w\rangle \in f)
 \end{eqnarray*}
 Left to right:  Assume $\langle x,w \rangle\in Ap(\s z,f)$.
    By Definition~\ref{definition:Ap},  that assumption is equivalent to 
 $$\exists y\, (\langle \langle f,y \rangle, z \rangle \in \s \ \land \ \langle x,y \rangle \in y.$$
 it suffices to prove 
 $$ t \in y \iff \exists y\, (\langle x,y \rangle \in f \ \land \ t \in y.$$
 Applying the definition of Church successor, in a few steps we obtain
 $$ t \in \FUNC \ \land \ \langle f,t\rangle \in z \ \land \ \langle x,q \rangle \in f  \ \land \ \langle q,w \rangle \in f.$$
 That completes the left-to-right direction.
 \smallskip
 
 Right to left:   Assume 
 $$ t \in \FUNC \ \land \ \langle f,t\rangle \in z \ \land \ \langle x,q \rangle \in t  \ \land \ \langle q,w \rangle \in f.$$
 We have to prove $\langle x,w\rangle \in Ap(\s z,f)$. Applying the definitions of $Ap$ and $\s$, we 
 find that it suffices to prove
 \begin{eqnarray*}
 && \exists y\,(z \in \FUNC \land \ \exists g,p\, (\langle f,y \rangle = \langle g,p \rangle \ \land \ 
 g \in \FUNC \ \land \  \\
 &&\forall x,w\,(\langle x,w\rangle \in p \iff 
 \exists t,q\, (t \in \FUNC \ \land \ \langle g,t\rangle \in z \ \land \ \langle x,q \rangle \in t \ \land \ \langle q,w\rangle \in g))) \\
 &&\land \ \langle x,w \rangle \in y).
 \end{eqnarray*}
 Now we choose 
 $$ Y := \{ \langle x,w \rangle: \exists t,q\,(t \in \FUNC \ \land \ \langle f,t\rangle \in z \ \land \ \langle x,q \rangle \in t  \ \land \ \langle q,w \rangle \in f.\}$$
 (It is important to quantify over $t$ and $q$ even though there are free variables $t$ and $q$ in scope here.)
 The formula is stratified giving $x,w,t,q$ index 0 and $z,f$ index 3; $\FUNC$ is a parameter.  Hence
 the definition can be given in \INF.  Using $Y$ to instantiate $\exists y$, we have to prove
 \begin{eqnarray*}
 &&   z \in \FUNC \land \ \exists g,p\, (\langle f,Y \rangle = \langle g,p \rangle \ \land \ 
 g \in \FUNC \ \land \  \\
 &&\forall x,w\,(\langle x,w\rangle \in p \iff 
 \exists t,q\, (t \in \FUNC \ \land \ \langle g,t\rangle \in z \ \land \ \langle x,q \rangle \in t \ \land \ \langle q,w\rangle \in g))) \\
 &&\land \ \langle x,w \rangle \in Y .
 \end{eqnarray*}
 We have $z \in \FUNC$ and $t\in \FUNC$;  take $g = f$ and $p = Y$; then it suffices to prove
 \begin{eqnarray*}
 &&\forall x,w\,(\langle x,w\rangle \in Y \iff 
 \exists t,q\, (t \in \FUNC \ \land \ \langle f,t\rangle \in z \ \land \ \langle x,q \rangle \in t \ \land \ \langle q,w\rangle \in f))) \\
 &&\land \ \langle x,w \rangle \in Y).
 \end{eqnarray*}
 The first line follows from the definition of $Y$.  It remains to prove $\langle x,w\rangle\in Y$. 
 Note that in the last line, $x$ and $w$ are free variables.   We have by assumption
 $$\langle f,t\rangle \in z \ \land \ \langle x,q \rangle \in t  \ \land \ \langle q,w \rangle \in f.$$
 Then by definition of $Y$ we have $\langle x,w\rangle\in Y$, as desired.  
 That completes the proof of the lemma.

  \begin{lemma}\label{lemma:sfunction2} For all $z$,  $z \in \FUNC \imp \s z \in \FUNC$.
 \end{lemma}
 
 \noindent{\em Proof}.  Immediate from Definition~\ref{definition:sucdef}.
 
 \begin{lemma}\label{lemma:srelation2} For all $z$, $\s z$ is a relation.
 \end{lemma}
 
 \noindent{\em Remark}. Since our definition of $\FUNC$ does not require a function
 to be a relation (i.e., contain only ordered pairs), this lemma adds something to the
 previous lemma.
 \medskip
 
\noindent{\em Proof}.  Immediate from Definition~\ref{definition:sucdef}.
 
\begin{definition} \label{definition:function}
$f: \N \to \N$ means   $f \in \FUNC$  and 
for each $n \in \N$ there is a 
unique $m \in \N$ such that $\langle n,m \rangle \in f$.
\end{definition}
The concept just defined does not prevent the domain of $f$ 
from being larger than $\N$.  

\begin{definition}\label{definition:ChurchFreged}
${\id}$ is the identity function, $\{ \langle x,x\rangle: x =x \}$.
\end{definition}

\begin{definition}\label{definition:ChurchZero}
$\ChurchZero$ is the function $\lambda f\,x.x$, which as a set
of ordered pairs is 
$$ \{ \langle f,  \id \,\rangle  : f=f \},$$
\end{definition}

\noindent{\em Remark}.  Then $\ChurchZero: X \to X$ for any set $X$,
but the domain of $\ChurchZero$ is the whole universe, usually larger than $X$.
\smallskip

\begin{lemma} \label{lemma:zeroFUNC} 
 $\ChurchZero \in \FUNC$.
\end{lemma}

\noindent{\em Proof}.  Immediate from the definition of $\ChurchZero$.

\begin{lemma} \label{lemma:ApZero} For all $x$, $\ChurchZero x = \id $, the identity function.
\end{lemma}

\noindent{\em Proof}.  By definition of $\ChurchZero$,  $\langle f,u \rangle \in \ChurchZero$
if and only if $u = \id $.    Then for any $f$,  $\langle f, \id \rangle \in 
\ChurchZero$,  so by Lemma~\sref{lemma:Ap}, we have $\id  = Ap\, (\ChurchZero, f) = \ChurchZero f$
as desired.  That completes the proof. 

\begin{lemma} \label{lemma:zeroAp}
For all $x$,  $\ChurchZero f x = x$.
\end{lemma}

\noindent{\em Proof}. By Lemma~\ref{lemma:ApZero},   $\ChurchZero f$ is the identity
function.  It follows from 
Lemma~\sref{lemma:Ap} that $x = Ap ((\ChurchZero f), x)$.  Suppressing explicit mention of $Ap$,
that is $\ChurchZero f x = x$.  That completes the proof of the lemma.

\begin{lemma} \label{lemma:ApOne} For all $f$, $f \in \FUNC \imp Rel(f) \imp \s \ChurchZero f = f$.
\end{lemma}

\noindent{\em Proof}
Lemma~\ref{lemma:ChurchSuccessor2}, with $z$ in the lemma set to $\ChurchZero$, 
 \begin{eqnarray*}
\s \ChurchZero f &=&   \{ \langle x, w\rangle :
             \exists t,q\,
                     ( t \in \FUNC \land \langle f,t \rangle \in \ChurchZero \land
                       \langle x,q \rangle \in t \land
                       \langle q,w \rangle \in f 
                      )
                \} 
 \end{eqnarray*}
 By the definition of $\ChurchZero$,  $\langle f,t \rangle \in \ChurchZero$ is 
 equivalent to $t = \id $; then $\langle x,q \rangle \in t$ is equivalent 
 to $x = q$, and we have 
 \begin{eqnarray*}
\s \ChurchZero f &=&   \{ \langle x, w\rangle : \langle x,w \rangle \in f \}  \\
                 &=& f \mbox{\qquad by extensionality and $Rel(f)$} 
 \end{eqnarray*}

\begin{definition}\label{definition:N}
The set of Church numbers is the least set containing $\ChurchZero$
and closed under Church successor $\s $.  That is, it is the intersection
of all sets containing $\ChurchZero$ and closed under $\s $. 
\end{definition}

\begin{theorem} \label{theorem:ChurchNumbers}
The set of Church numbers is definable (by stratified formulas) in \INF.
\end{theorem}

\noindent{\em Proof}.
By Lemma~\ref{lemma:sucdef}, $\s$ is definable in NF. Then 
$$ \N = \{\ x : \forall X(0 \in X \land \forall u \in X( \s u \in X)) \imp x \in X\ \} $$
In other words,
$$ \N = \{\ x : \exists z (z = 0) \land
 \forall X(z \in X \land \forall u \in X \exists v( v= \s u \land v\in X)) \imp x \in X \ \} $$
To check that this definition is legal in NF, we stratify the 
formula on the right, giving $v$ and $u$ the same index, say 6,
since $\s$ is a function.  Then $X$ gets index 7 and $x$ and $z$ get index 6.
As discussed above we can stratify $z = 0$ giving $z$ any
index $\ge 3$, so 6 is OK.  $\N$ is a parameter and does not need an index.
This stratification shows that $\N$ is well-defined in NF.
That completes the proof of the theorem.

\begin{lemma}[Proof by induction] \label{lemma:Ninduction}
 
$$ \ChurchZero \in x \land \forall n \in \N(n \in x \imp \s n \in x) \imp \forall k \in \N\,( k \in x).$$
\end{lemma}

\noindent{\em Proof}.
 By Definition~\ref{definition:N} and 
 Theorem~\ref{theorem:ChurchNumbers},
 $\N$ is the 
intersection of all sets closed under successor.  There is at least
one such set, since $\V$ is closed under successor and contains $0$,
so $\N$ is not empty.    Then $x$ contains $0$ and is closed under
successor. Hence $\N \subseteq x$.  That completes the proof.
\smallskip

{\em Remark}.  If we wish to prove a stratified formula $\phi$  ``by induction on $n$'', we use
stratified comprehension to define $x = \{ n: \phi(n)\}$, and then prove the ``base case'' that $\Phi(0)$and the ``induction step'' that $\phi(n) \imp \phi(\s n)$.  Then $\ChurchZero \in x$ and 
$n \in x \imp \s n \in x$.  Then Lemma~\ref{lemma:Ninduction} can
be used to conclude that $\N \subseteq x$.  Hence $\forall n \in \N\, \phi(n)$. 
\smallskip 

\begin{lemma}\label{lemma:zeroN}
 $\ChurchZero \in \N$.
\end{lemma}

\begin{lemma}\label{lemma:successorN}
 $\s : \N \to \N$
\end{lemma}

\noindent{\em Proof}. By definition $\N$ is the 
intersection of all sets containing $\ChurchZero$ and closed under successor.   
Therefore $\ChurchZero \in \N$.   That completes the proof of the lemma.

\noindent{\em Proof}. By definition $\N$ is the 
intersection of all sets containing $\ChurchZero$ and closed under successor.   
Let $n \in \N$.  Then $n$ belongs to every set $X$ containing $\ChurchZero$ and 
closed under successor.  Hence $\s n$ belongs to every such set $X$.
Hence $\s n$ belongs to $\N$.  That completes the proof of the lemma.

\begin{lemma} \label{lemma:Churchnumbersarefunctions}
Every Church number $n$ is a function.
\end{lemma}

\noindent{\em Proof}.  By induction on $n$.
\smallskip

Base case: $\ChurchZero \in \FUNC$ by Lemma~\ref{lemma:zeroFUNC}.
\smallskip

Induction step:  Suppose $n \in \N$ and $n \in \FUNC$. 
By Lemma~\ref{lemma:sfunction2}, $\s n \in \FUNC$.
That completes the induction step.  That completes the proof of the lemma.

\begin{lemma} \label{lemma:Churchnumbersarerelations}
Every Church number $n$ is a relation.
\end{lemma}

\noindent{\em Proof}.  By induction on $n$, similar to Lemma~\ref{lemma:Churchnumbersarefunctions},
but appealing to Lemma~\ref{lemma:srelation2} in 
the induction step. 

 \begin{lemma} \label{lemma:ChurchSuccessor4} Let $n \in \FUNC$ and $f \in \FUNC$.  
 Then there exists $y$ such that $\langle f,y \rangle \in \s n$.  
 \end{lemma}
 
 \noindent{\em Proof}.  Let $n$ and $f$ be given, with $n \in \FUNC$ and $f \in \FUNC$.
 Define
 $$ y := \{ \langle x, z\rangle : \langle f,p \rangle \in n \ \land \ \langle x,q \rangle \in p
 \ \land \ \langle q,z \rangle \in f \ \land \ p \in \FUNC\}.$$
 The formula defining $y$ is stratified, giving $x$, $z$, and $q$ index 0; then $\langle x, q\rangle$
 gets index 2, so $p$ gets index 3.  Then $f$ gets index 3 and $\langle f,p \rangle$ gets
 index 5, so $n$ gets index 6.  Therefore the formula is stratified, and the definition 
 of $y$ is legal.
 
 The verification that $\langle f, y \rangle \in \s n$ then proceeds by unfolding the 
 definitions of $\s n$ and $y$.  We omit the 65 routine steps of this verification.   
 
\section{Iteration of a function}
If we have a mapping $f:X\to X$, we can iterate it $j$ times.  Often mathematicians
write the $j$-times iterated mapping as $f^j$, or if there is danger of confusion, as $f^{(j)}$.
Formally it is just $jf$, where $j$ is a Church number.  In treating this subject rigorously
one has to distinguish the relevant concepts precisely.  Namely, we have
\begin{eqnarray*}
f:X \to Y \\
Rel(f) \\
f \in \FUNC \\
oneone(f,X,Y) 
\end{eqnarray*}
$Rel(f)$ means that all the members of $f$ are ordered pairs.  $f \in \FUNC$ means that 
two ordered pairs in $f$ with the same first member have the same second member.  (Nothing
is said about possible members of $f$ that are not ordered pairs.)  $f:X \to Y$ means that 
if $x \in X$, there is a unique $y$ such that $\langle x,y \rangle \in f$ and that $y$ is in $Y$.
(But nothing is said about $\langle x, y\rangle \in f$ with $x \not \in X$.)  ``$f$ is one-to-one
from $X$ to $Y$'',  or $oneone(f,X,Y)$,  means $f:X \to Y$ and in addition, if $\langle x,y \rangle \in f$
and $\langle u,y\rangle \in f$ then $x = u$,  and if $y \in Y$ then $x \in X$.  (So $x=u$ does
not require $y \in Y$ or $x \in X$.)  In particular, $f:X \to Y$ does not require $dom X \subseteq X$,
so the identity function maps $X$ to $X$ for every $X$; but the identity function (on the universe)
has to be restricted to $X$ before it is one-to-one.  

We shall be mostly concerned with iterations of a map $f$ from some set $X$ to that same set.
In that setting the following concept is useful.

\begin{definition}\label{definition:permutation2} 
$f$ is a {\bf permutation} of a finite set $X$ if and only if $f:X \to X$,
and $Rel(f)$ and $f \in \FUNC$, and $dom(f) \subseteq X$,
and $f$ is both one-to-one and onto from $X$ to $X$. 
\end{definition} 

But for some purposes, we don't need $f$ to be onto,  but we still need it 
to be a relation and a function and to control its range and domain.  Therefore we define

\begin{definition}\label{definition:ChurchFregenjection}
$f$ is an {\bf injection} of a set $X$ into $Y$ if and only if $f:X \to Y$,
and $Rel(f)$ and $f \in \FUNC$, and $dom(f) \subseteq X$,
and $f$ is  one-to-one  from $X$ to $X$.
\end{definition} 

Note that the definition does not require $X$ to be finite.
\smallskip

Any function can be iterated, even if it doesn't map some $X$ to itself:

\begin{lemma} \label{lemma:nf_defined} Let $n \in \N$ and $f \in \FUNC$.  Then there
exists $y$ such that $\langle f,y\rangle \in n$.
\end{lemma}

\noindent{\em Proof}.  The formula is stratified, so we may prove it by induction.
\smallskip

Base case: By the definition of Church zero, we have $\langle f, \id  \rangle \in \ChurchZero$.
\smallskip

Induction step: By Lemma~\ref{lemma:ChurchSuccessor4}.  That completes the induction step.
That completes the proof of the lemma.

\begin{lemma} \label{lemma:nfFUNC} Let $n \in \N$ and suppose $f \in \FUNC$ and $Rel(f)$. 
Then $nf \in \FUNC$  and $Rel(nf)$.
\end{lemma}

\noindent{\em Proof}. By induction on $n$, which is legal since the formula is stratified.
(Although $Ap(n,f)$ gets the same type as $f$, that observation is not even needed here, as $\FUNC$
is just a parameter, so $f$ can be given any type and it doesn't matter what type $nf$ gets.)
\smallskip

Base case: $\ChurchZero f$ is the identity function, by definition of $\ChurchZero$.  Since the 
identity function is also a a relation, that 
completes the base case (though it requires 24 steps, here omitted, to spell out the details). 
\smallskip

Induction step: Suppose $f \in \FUNC$ and $n \in \N$ and $nf \in \FUNC$. 
 By Lemma~\ref{lemma:ChurchSuccessor2},
we have 
\begin{eqnarray}
\s nf &=&   \{ \langle x, w\rangle :
             \exists t,q\,
                     ( t \in \FUNC \land \langle f,t \rangle \in n \land
                       \langle x,q \rangle \in t \land
                       \langle q,w \rangle \in f 
                      )
                \}  \label{eq:E527}
 \end{eqnarray}
 Then $\s nf$ is a relation.  We next will prove $\s nf \in \FUNC$. 
 Suppose $\langle x,y\rangle \in \s nf$ and $\langle x,z\rangle \in \s nf$. We must prove $y=z$.
 By (\ref{eq:E527}) there exist $t_1,q_1$ and $t_2,q_2$ such that 
 \begin{eqnarray*}
 \FUNC(t_1) &&\\
 \FUNC(t_2) &&\\
 \langle f,t_1\rangle \in n && \\
 \langle f,t_2\rangle \in n &&\\
 \langle x,q_1 \rangle \in t_1 &&\\
 \langle x,q_2 \rangle \in t_2 &&\\
 \langle q_1,y \rangle \in f && \\
 \langle q_2,z \rangle \in f &&
 \end{eqnarray*}
 Using the definition of $\FUNC$ several times we obtain, in order, $t_1 = t_2$, then $q_1 = q_2$,
 and finally $y=z$.   
That completes the induction step.
   That completes the proof of the lemma.

\begin{lemma} \label{lemma:iteration}
Let $X$ be any set.  Suppose $f: X \to X$
and $f \in \FUNC$  and $Rel(f)$). 
Then for all
$n\in \N$ and $x \in X$, 
$$ \langle f, nf \rangle \in n \ \land \ nf : X \to X$$
\end{lemma}

\noindent{\em Proof}.   Let $f \in \FUNC$ and $Rel(f)$ and $f:\N \to \N$. 
By Lemma~\sref{lemma:Ap}, $\langle f, nf \rangle \in n$ is equivalent 
to $\exists y\, (\langle f,y \rangle \in n)$.
\smallskip

The formula is stratified, giving $X$ index 1, $f$ index 3 (since 
the members of $f$ are ordered pairs of members of $X$),  $x$ index 0;
$nf$ gets index 3, since $n$ is a function by Lemma~\ref{lemma:Churchnumbersarefunctions};
so we have to give $n$ index 6, since its members are pairs of objects of type 3.

Base case: 
By Lemma~\ref{lemma:zeroFUNC}, $\ChurchZero$ is a function, and by definition of $\ChurchZero$,
$\ChurchZero f$ is the identity function,  so $\ChurchZero f: X \to X$.  That completes the 
base case.
\smallskip

Induction step:  We first have to show that $\langle f, \s nf \rangle \in \s n$.
By Lemma~\ref{lemma:ChurchSuccessor4}, we have $\exists y\,( \langle f, y \in \s n)$.
Then by Lemma~\sref{lemma:Ap}, we have $\langle f, \s n f \rangle \in \s n$, as claimed.
\smallskip

We turn to the proof that $\s nf : X \to X$.
 By Lemma~\ref{lemma:Churchnumbersarefunctions}, 
  $n$ is a function, and by Lemma~\ref{lemma:nfFUNC},
$nf$ is a function.   Then according to 
  Lemma~\ref{lemma:ChurchSuccessor2},
we have 
\begin{eqnarray}
\s nf &=&   \{ \langle x, w\rangle :
             \exists t,q\,
                     ( t \in \FUNC \land \langle f,t \rangle \in n \land
                       \langle x,q \rangle \in t \land
                       \langle q,w \rangle \in f 
                      )
                \}  \label{eq:658}
 \end{eqnarray}
 Let $x \in X$.  By Lemma~\ref{lemma:nf_defined},
$\langle f,t\rangle \in n$ for some $t$; by Lemma~\sref{lemma:Ap}, $t = nf$. 
By the induction hypothesis, $nf:X \to X$, so there exists $q$ with $\langle x,q\rangle \in t$
and $q \in X$.  Then since $f:X \to X$, there exists $w$ with $\langle q,w \rangle \in f$.
Then by (\ref{eq:658}), we have $\langle x,w \rangle \in \s nf$.  Since $x$ was arbitrary,
we have proved $\s nf : X \to X$. 
That completes the induction step.  That completes the proof of the lemma.
 
\begin{theorem}[successor equation] \label{theorem:successorequation}
Let $X$ be any set and $f$ any function ($f \in \FUNC$ and $Rel(f)$) 
with $f:X \to X$. 
Then for all $n\in \N$ and $x \in X$, 
$$ \s  n f x = f(n f x).$$ 
\end{theorem}

\noindent{\em Proof}.   Let $f:X \to X$ and $n \in \N$ and $x \in X$.  
 By Lemma~\ref{lemma:Churchnumbersarefunctions}, 
  $n$ is a function, and by Lemma~\ref{lemma:nfFUNC},
$nf$ is a function.   Then according to 
  Lemma~\ref{lemma:ChurchSuccessor2},
we have 
\begin{eqnarray*}
\s nf &=&   \{ \langle x, w\rangle :
             \exists t,q\,
                     ( t \in \FUNC  \land \langle f,t \rangle \in n \land
                       \langle x,q \rangle \in t \land
                       \langle q,w \rangle \in f 
                      )
                \}  
 \end{eqnarray*}
By Lemma~\ref{lemma:iteration}, $nf:X \to X$.  
Then as in the proof of that lemma, we have $t = nf$, $q = tx = nfx$, and $w = fq$,
with $\langle x,w \rangle \in \s nf$.  That completes the proof of the theorem.

\begin{lemma} \label{lemma:Church1notequal0} 
Define $1 := \s 0$.  Then $1 \neq 0$.
\end{lemma}

\noindent{\em Proof}.
Let $f$ and $z$ be functional relations.  Then by Lemma~\ref{lemma:ChurchSuccessor3},
$\s zf$ is a relation, and by 
 Lemma~\ref{lemma:ChurchSuccessor2},
\begin{eqnarray*}
\langle x,w \rangle \in \s zf &\iff& 
           \exists t,q(t \in \FUNC  \land \langle f,t\rangle \in z
                       \land \langle x,q \rangle \in t 
                       \land \langle q,w \rangle \in f) 
\end{eqnarray*}
Take $z = 0$.  Then on the right, $t = zf$ and $q = tx = zfx = x$
since $z=0$.   Then $w = fx$.  Thus
\begin{eqnarray*}
\langle x,w \rangle \in \s 0f &=&   \langle x,w\rangle \in   f
\end{eqnarray*}
By Lemma~\ref{lemma:ChurchSuccessor4}, $\s 0f$ is a relation, and 
by hypothesis $f$ is a relation.  Therefore 
$$ \forall u\, (u \in \s 0 f \iff u \in f).$$
Then by extensionality, 
$1 f = f$, for all functional relations $f$.

  Now suppose,
for proof by contradiction, that $1 = 0$.  Then 
on the one hand, $1f =f$, and on the other hand 
$1f  = 0f  = \id  $.  Now we can get a contradiction
by exhibiting some (any) functional relation $f$ that
is not the identity. For example, we can use $f = \{ \langle \empty, \{\empty\}\}$.
It is easily verified that $f$ is a functional relation and is not
equal to $\id $.  That completes the proof of the lemma. 
\smallskip

\noindent{\em Remark}. The proof does not follow immediately
from Theorem~\ref{theorem:successorequation}, it seems.
For if we assume $1=0$, that equation says
$$  1 f x\  = \ \s 0 fx \ = \ f(0f x)$$
and since we have assumed 1 = 0, also $1 fx = 0fx$.
But by definition of $0$, we have $0 f x = x$.  
Thus $x = 0fx = 1fx = f(0fx) = fx$.  Hence $f$ is 
the identity function on $\N$.  That is, however,
not yet a contradiction.

\begin{theorem} \label{theorem:successoromitszero}
The Church successor function 
does not take the value $\ChurchZero$ on $\N$; that is, $n \in \N \imp \s n \neq \ChurchZero$.
\end{theorem}

\noindent{\em Proof}.  Let $a$ and $b$ be any two unequal members of $\N$;
by Lemma~\ref{lemma:Church1notequal0} there do exist two unequal members of $\N$.
Let $f$ be
the constant function with value $a$.  Then $f:\N \to \N$.
 Suppose, for proof by 
contradiction, that $\s (z) = 0$.  Applying both sides to 
$f$ and $b$ we have, by the definitions of $\s $ and $0$, 
\begin{eqnarray*}
\s zfb &=& 0fb\\
f(zfb) &=& 0fb \mbox{\qquad by Theorem~\ref{theorem:successorequation}}\\
f(zfb) &=& b  \mbox{\qquad since $0fb = b$ by definition of $0$}\\
a &=& b \mbox{\qquad since $f(u) = a$ for all $u$}
\end{eqnarray*}
But that contradicts $a \neq b$.  That completes the proof of the 
theorem.
\smallskip

{\em Remark}.  All we needed to prove that successor omits the value $0$
is that there is some function $f:\N \to \N$ that omits some value; 
and we can construct such a function if 
there are two distinct members of $\N$.

\begin{lemma}[Predecessor] \label{lemma:predecessor}
If $n \in \N$ and $n \neq \ChurchZero$,
then $n = \s m$ for some $m \in \N$.
\end{lemma}

\noindent{\em Remark}.  The predecessor is, of course, not asserted
to be unique.
\medskip

\noindent{\em Proof}. By induction on $x$ we prove 
that $x \in \N \imp x \neq 0 \imp \exists y\in \N\,(\s y = x)$. 
The base case and induction step are both
immediate.  

\begin{lemma}\label{lemma:decidable0} $\forall n \in \N\,(n = \ChurchZero \ \lor \ n \neq \ChurchZero)$.
\end{lemma}

\noindent{\em Proof}.  By induction on $n$, which is legal since the formula is stratified. 
The base case is immediate; and the induction step is 
immediate from Theorem~\ref{theorem:successoromitszero}.  That completes the proof of the lemma.

\begin{lemma} \label{lemma:iterationFUNC}
  Suppose $f: X \to X$, and $f \in \FUNC$ and $Rel(f)$,
and $dom(f) \subseteq X$.  Let $n \in \N$ and suppose $n \neq \ChurchZero$.  Then 
$dom(nf) \subseteq X$.
\end{lemma}

\noindent{\em Remark}.  When $n = \ChurchZero$, $nf$ is the identity function, whose
domain is $\V$.  Hence the restriction $n \neq \ChurchZero$ is necessary.
\medskip

\noindent{\em Proof}.  By induction on $n$, which is legal since the formula is stratified.
Base case:  there is nothing to prove because of the hypothesis $n \neq \ChurchZero$.
\smallskip

Induction step: Assume $\s n \neq \ChurchZero$.  
 By Lemma~\ref{lemma:decidable0}, $n = \ChurchZero \ \lor \ n \neq \ChurchZero$.  We argue
by cases accordingly.
\smallskip

Case 1, $n = \ChurchZero$.  We have to show $dom (\s \ChurchZero f) \subseteq X$.
  Suppose $x \in dom (\s nf)$.  It suffices to show $x \in X$.
    By Lemma~\ref{lemma:ApOne}, $\s\ChurchZero f = f$.
Therefore $x \in dom(f)$.  Since $dom(f) \subseteq X$, we have $x \in X$ as desired.
That completes Case~1.   
\smallskip

Case 2, $n \neq \ChurchZero$.  By Lemma~\ref{lemma:nfFUNC}, we have $nf \in \FUNC$ and $Rel(nf)$,
and also $\s n f \in \FUNC$ and $Rel (\s n f)$.
  Suppose
$\langle x,y \rangle \in \s nf$.  We must show $x \in X$.  Since $\s nf \in \FUNC$, 
we have $y =  (\s n f) x$.  By Lemma~\ref{lemma:ChurchSuccessor2}, we have 
\begin{eqnarray*}
x \in dom (\s nf) \iff \exists w,t,q\,(t \in \FUNC \ \land \ \langle f,t \rangle \in n 
\ \land \ \langle x,q \rangle \in t \ \land \ \langle q,w\rangle \in f) 
\end{eqnarray*}
Suppose $x \in dom (\s nf)$.  Then for some $w,t,q$ we have 
\begin{eqnarray*}
 t \in \FUNC \ \land \ \langle f,t \rangle \in n 
\ \land \ \langle x,q \rangle \in t \ \land \ \langle q,w\rangle \in f 
\end{eqnarray*}
Then by Lemma~\sref{lemma:Ap}, $ t = nf$.   By the induction hypothesis, 
$dom(nf) \subseteq X$.  Since $\langle x,q \rangle \in t$ and $t = nf$, we have 
$x \in dom(nf)$ and hence $x \in X$.  That completes Case~2.  That completes
the induction step.  That completes the proof of the lemma.

\begin{lemma}\label{lemma:oneoneiteration_helper}
Suppose  $f \in \FUNC$ and $Rel(f)$, and $f:X \to X$ and $m \in \N$ 
and $ mf$ is one-to-one
from $X$ to $X$,  and $dom(f) \subseteq X$.  Then $\s m f :X \to X$,
and $\s m f$ is  
one-to-one.
\end{lemma}

\noindent{\em Proof}.  
Suppose $\s m x = \s m z$.  I say $x = z$.  We have
\begin{eqnarray*}
\s m f x = f (m f x))  && \mbox{\qquad by Theorem~\ref{theorem:successorequation}}\\
\s m f z = f (m f z))  && \mbox{\qquad Theorem~\ref{theorem:successorequation}}\\
f( m f x) = f( m f z)  && \mbox{\qquad by the preceding two lines}\\
x = z  &&   \mbox{\qquad since $mf$ is one-to-one}
\end{eqnarray*}
Technically, however, the definition of one-to-one involves more than 
just $x=z$.  We have
\begin{eqnarray*}
\s m f \in \FUNC && \mbox{\qquad by Lemma~\ref{lemma:nfFUNC}}
\end{eqnarray*}
We also have to show that $m f$ is a relation,
that its domain is a subset of $X$, and that its range is a subset of $X$.
These verifications require about 100 proof steps (here omitted),
using for example Lemmas~\ref{lemma:Churchnumbersarefunctions} and
~\ref{lemma:Churchnumbersarefunctions} 
and \ref{lemma:iterationFUNC}.
That completes the proof of the lemma.

\begin{lemma} \label{lemma:oneoneiteration}
Let $X$ be any set, and let $f:X \to X$ be a permutation.
 Let $m \in \N$.  Then $mf:X \to X$  for  $m \neq \ChurchZero$, $mf$ is a permutation.   
\end{lemma}

\noindent{\em Remark}.  When $ m = \ChurchZero$, $m f$ is the identity function,
which has domain $\V$, so it is not a permutation of $X$ (unless $X = \V$).  
\medskip

\noindent{\em Proof}.  By induction on $m$.  The formula is stratified, as we 
have already checked that ``one-to-one'' and ``$f:X \to X$'' are stratified.
\smallskip

Base case, when $m = \ChurchZero$ we have $mf = \id $, by definition of $\ChurchZero$,
and the identity function maps $X$ to $X$, and the identity
function is one-to-one.  That completes the base case.
\smallskip

Induction step.  Suppose $mf:X\to X$ is one-to-one, and   
$$ m \neq \ChurchZero \imp range(mf) \subseteq X.$$
We must prove $\s mf:X \to X$ and $\s m$ is one-to-one.  We have $m = \ChurchZero \ \lor \ m \neq \ChurchZero$,
by Lemma~\ref{lemma:decidable0}. We argue by cases.
\smallskip

Case~1: $m = \ChurchZero$.  Then $mf = f$, so by hypothesis, $mf$ is a permutation
of $X$.
\smallskip

Case~2: $m \neq \ChurchZero$.  
Assume $x \in X$.  Then by Lemma~\ref{lemma:oneoneiteration_helper},
$\s m f:X \to X$  and $\s m f$ is one-to-one.   That completes Case~2. 
That completes the induction step.  That completes the proof of the lemma.

\section{Definition of addition on $\N$}
In this section we define addition on $\N$ and prove some of 
its properties.   To define the graph of a binary function 
we use ordered triples, which are defined in Definition~\sref{definition:triples}.

\begin{lemma} \label{lemma:sum} \INF\ can define a 
 set $Sum$ such that for $x,n,y \in \N$
\smallskip

(i) $\langle x,0,x  \rangle \in Sum$, and 
\smallskip

(ii) $\langle x,n,y \rangle\in Sum \imp
 \langle x, \s n, \s y \rangle \in Sum$, and
 \smallskip
 
(iii) $Sum$ is the intersection of all sets $X$
satisfying those two conditions; 

\end{lemma}

\noindent{\em Proof}.  $Sum$ is the intersection of 
sets $X$ satisfying  conditions (i) and (ii) 
with $Sum$ replaced by $X$.  Specifically those  conditions are 
\begin{eqnarray*}
 x \in \N &\imp& \langle x,0,x \rangle \in X  \\
x \in X \land y \in \N \land n \in \N &\imp& 
(\langle x,n,y \rangle \in X \imp
 \langle x, \s n, \s y \rangle \in X)
 \end{eqnarray*}
 These formulas can be stratified by assigning  $x$, $y$, and $n$
 all index 0, $\N$ index 1, and $X$ index 5.

 Then the conjunction of these two 
 conditions, preceded by $\forall X$, is also stratifiable, and 
 it defines $Sum$. 
 
 Now we must prove that $Sum$ so defined satisfies the 
three conditions itself.  Suppose $x \in \N$.  Then 
 $\langle x,0,x \rangle$ belongs to 
 every $X$ satisfying the conditions. Hence it belongs 
 to $Sum$.  Hence $Sum$ satisfies the first condition.
 
 Suppose $x\in \N$ and $y \in \N$ and $n \in \N$.   Suppose 
 $\langle x,n,y \rangle \in Sum$.
 Then for every $X$ satisfying the conditions,
  $\langle x,n,y \rangle \in X$. 
 Then for every $X$ satisfying the conditions,
 $\langle x,\s n, \s y\rangle \in X$.
 Then 
  $\langle x, \s n, \s y \rangle \in Sum$.
 That verifies that $Sum$ satisfies the second condition.
 That completes the proof of the lemma.
\smallskip

\begin{lemma} \label{lemma:additionfunction} \INF\ proves that
for each $x,n \in \N$, there is a unique $y \in N$ such that 
$\langle x,n,y \rangle \in Sum$.  
\end{lemma}

\noindent{\em Proof}. (Holmes) First, by induction on $n$, there is 
some $y$ such that $\langle x, n, y \rangle \in Sum$; the
clauses (i) and (ii) in the definition yield the base case
and induction steps, respectively.  So it suffices to prove
by induction on $x$ that 
$$ \forall y \in \N\,(\langle x,y,z \rangle \in Sum \land \langle x,y,w \rangle\in Sum 
    \imp z = w).$$

Base case: Since $\langle 0,y,y \rangle \in Sum$, it suffices to prove that
$\langle 0,y,z\rangle \in Sum$
implies $z = y$.  To that end define 
$$ X = \{ \langle x,y,z \rangle : x = 0 \imp y = z \}.$$
That is legal as the formula is stratifiable.
I say that $X$ satisfies the closure conditions in the definition
of $Sum$.  Ad~(i): We have $\langle x, 0, x \rangle \in X$,
since $x = 0 \imp 0 = x$.  
Ad~(ii):  Suppose $\langle x, y, z \rangle \in X$.
We must show $\langle x, \s y, \s z \rangle \in X$. 
Since $\langle x, y, z \rangle \in X$, we have $x = 0 \imp y = z$.
Then also $x=0 \imp \s y = \s z$.  Therefore 
$\langle x, \s y, \s z \rangle \in X$.  Therefore $X$ satisfies 
both closure conditions.  Therefore $Sum$ is a subset of $X$. 
Therefore $\langle 0,y,z \rangle \in Sum$  implies $y=z$,
as desired.   That completes the base case.
\smallskip

Induction step:  The induction hypothesis is that (with $x$ fixed) 
for every $y \in \N$, there
is a unique $p \in \N$ such 
that $\langle x, y, p \rangle \in Sum$. We denote that unique $p$
by $x\oplus y$.  
Suppose that 
\begin{eqnarray}
&& \langle \s x, y, z \rangle \in Sum  \label{eq:707} \\
&&\langle \s x,y, w \rangle \in Sum \label{eq:708}
\end{eqnarray}
  We must prove $z=w$.
\smallskip

We define a set $X$ (depending on $x$, which is now fixed until 
we finish the induction step):
$$X = \{ \langle u,v,z \rangle : u = \s x \imp z = \s(x\oplus v)\}.$$
That formula is stratifiable, as ``$z = \s(x\oplus v)$'' can be replaced
by $\exists p\,(Sum(x,v,p) \land \s p = z$'', and all the variables
can be given the same type.  Hence the definition of $X$ is legal.
I say that $X$ satisfies the closure conditions in the definition 
of $Sum$.  Ad~(i): We must show $\langle x, 0,x \rangle \in X$.
That holds if and only if $x = \s x \imp x = \s(x\oplus 0)$.
But $x\oplus 0 = 0$, so the condition is $x = \s x \imp x = \s x$,
which is indeed valid. 
\smallskip

Ad (ii):  Suppose $\langle u,v,z \rangle \in X$.  
We must show
$\langle u, \s v, \s z \rangle \in X$. 
We have
\begin{eqnarray*}
u = \s x \imp z = \s(x\oplus v)  && \mbox{\qquad since $\langle u,v,z \rangle \in X$}
\end{eqnarray*}
By the definition of $Sum$ we have $\s(x\oplus v) = x \oplus  \s v$.  Therefore
$$ u = \s x \imp z = x \oplus  \s v.$$
Taking the successor of both sides of the equation after the implication,
$$ u = \s x \imp \s z = \s(x\oplus \s v).$$
By the definition of $X$, this is equivalent to 
$$\langle u, \s v, \s z \rangle \in X.$$
That completes the verification that $X$ satisfies (ii).
Hence $Sum$ is a subset of $X$. Then $\langle u,v,z \rangle \in Sum$
implies $\langle u,v,z\rangle \in X$.  Take $u = \s x$.  Then by 
definition of $X$, we have 
$$ \langle \s x, v, z \rangle \in Sum \imp z = x\oplus v.$$
Applying this to (\ref{eq:707}) and (\ref{eq:708}) we have 
$ z = x \oplus v$ and $w = x\oplus v$.  Therefore $z=w$ as desired. That completes
the induction step.  That completes the proof of the lemma.

Lemma~\ref{lemma:additionfunction} allows us to make the
following definition.

\begin{definition} \label{definition:Churchaddition}
We henceforth write $x \oplus  n = y$ instead 
of $\langle x,n,y \rangle \in Sum$, 
and when $x,n \in \N$, we write
 $x\oplus n$ for the unique $y$ such that $x\oplus n = y$. 
\end{definition}

\noindent{\em Remark}.  We already used ``$x+y$'' for addition
of Frege numerals in \cite{inf-basics}.  While we never need
addition of Frege numerals in this paper, we have chosen to 
keep the notation consistent between the two papers, by using
a different symbol for addition of Church numbers.

\begin{lemma}\label{lemma:ChurchZero_equation} 
$\forall x\in \N\,(x\oplus 0 = x)$.
\end{lemma}

\noindent{\em Proof}.  By Definition~\ref{definition:Churchaddition},
this formula can be expressed in terms of $Sum$ 
as $\langle x,0, x\rangle \in Sum$, which is
 proved in Lemma~\ref{lemma:sum}.
 
\begin{lemma}\label{lemma:ChurchAddition_equation}
$\forall x,n \in \N\,  (x \oplus  \s n = \s (x\oplus n))$.
\end{lemma}

\noindent{\em Proof}.  By Definition~\ref{definition:Churchaddition},
this formula can be expressed in terms of $Sum$ 
as $\langle x,   n, y\rangle \in Sum \imp  \langle x, \s n, \s y \rangle \in Sum$, which is
 proved in Lemma~\ref{lemma:sum}.

\section{Alternate definitions of addition}
In this section we discuss two definitions that we do not use,
and the reasons we do not use them.
\subsection{Addition as iterated successor}
We could consider defining addition by 
\begin{eqnarray}
 x \oplus  y &:=& y\s x  \label{eq:825}
\end{eqnarray}
Technically we have  defined the ``add $y$'' function
$y \s$,  which takes an argument $x$ and adds $y$ to it. 

This definition of addition as iterated successor makes it 
immediate that addition is single-valued, but the 
defining formula (\ref{eq:825}) cannot be stratified 
giving $x$ and $y$ the same type, for if we give $x$ index 0,
then $s$ has index 3 and $y$ has to get index 6.  So 
this definition does not make addition a function of 
the ordered pair $\langle x,y \rangle$.  

The 
laws of addition follow from the definition of successor:
\begin{eqnarray*}
x \oplus  \s y &=& (\s y) \s x \\
         &=& \s y \s x \\
         &=& \s( y \s x) \mbox{\qquad by definition of $\s$} \\
         &=& \s(x\oplus y)  \mbox{\qquad by (\ref{eq:825})}
\end{eqnarray*}
\begin{eqnarray*}
x \oplus  0 &=& 0\s x \\
      &=& 0 \mbox {\qquad by definition of $0$}
\end{eqnarray*}

With this definition of addition, the formula $x\oplus y=z$
is $y \s x = z$, which is stratified since it has only one occurrence
of each variable, but not homogeneous.  For example the formula
$0 \oplus  x = x$ is $ x \s 0 = x$, which is not stratified.  Hence, with 
this definition, we would not see how to prove $0\oplus x = x$.  Similarly,
the formula asserting the equivalence of the two definitions is 
$$ x \oplus  y = y \s x,$$
where $\oplus $ means the first definition.  This is not a stratified formula, 
since $y$ on the right must get a greater index than on the left.
Hence we cannot prove,  at least not by induction on $y$, that
the two definitions are equivalent.   This gives us a 
second reason not to use this 
definition.
 
\subsection{Addition via composition}\label{section:churchaddition}
Church and Kleene (in \cite{church1941} and \cite{kleene1935})
define addition to satisfy this formula:
\begin{eqnarray}
 (x\oplus y)fz &=& xf(yfz)\label{eq:861}
\end{eqnarray}
This formula is stratified giving $z$ index 0, $x$ and $y$ both index 6, and
$f$ index 3, so it is possible to give this definition in INF.   
The set-theoretical definition of addition given in Definition~\ref{definition:Churchaddition} produces
an addition function defined only on the Church numbers;  the more
general definition here can add any two functions mapping some set 
into itself, not just mapping $\N$ into $\N$. 

Unlike the definition of addition by iterated successor, there is 
no compelling reason not to use the Church-Kleene definition.  But
there are several details to attend to in translating from the 
$\lambda$-calculus to NF, for example, just to go from the definition 
above to the set of ordered triples that is really the function $\oplus $.
We wrote out all the details required to reach the basic properties
of addition, and found it required twice as much 
space as the set-theoretic details using Definition~\ref{definition:Churchaddition}.
 We shall see in 
 Lemma~\ref{lemma:doubleiteration} that (\ref{eq:861}) is satisfied  
by the addition of Definition~\ref{definition:Churchaddition}.

From (\ref{eq:861}) and the equation for successor we find 
$ \s(x\oplus y) = \s x \oplus  y$ and $0\oplus y = y$.  
From these we can prove the equivalence of this definition to 
the one given in Definition~\ref{definition:Churchaddition}, when 
restricted to Church numbers $x$ and $y$.

\section{Stratification}
Let $L$ be the fragment of the language of Peano 
arithmetic that does not involve the symbol for multiplication; 
thus $L$ has a constant $0$ and function symbols for 
successor and addition, from which compound terms can be 
built up.

Now that we have defined addition on $\N$, it is possible 
to define an interpretation of (the language of) $L$
into NF (which does not have terms, constant symbols, or 
function symbols).  Namely, for each term $t$ in $L$
with free variables $x$ 
there is a formula of NF with free variables 
$x$ and one additional variable $y$
expressing $t = y$.  This formula contains many fresh 
existentially quantified variables; rather than give a 
recursive definition, or a program for computing it, we 
illustrate with an example.  If $t$ is $x \oplus  s(z)$,
then the formula in question is 
$$\exists u,v\, (\langle x,u, y\rangle \in Sum
   \land \langle z,u \rangle \in \s)
 $$ 
where $w \in \s$ abbreviates the formula in Definition~\ref{definition:sucdef}, and 
$w \in Sum$ stands for the formula defining $Sum$.
Similarly, the formula $y = 0$ is expressed by the 
formula in Definition~\ref{definition:sucdef}.

\begin{lemma} \label{lemma:stratification} 
Any formula in the language 
of Peano arithmetic without multiplication is 
interpreted by a formula of NF that can be stratified
by giving all the variables the same type.
\end{lemma}

\noindent{\em Proof}. By induction on the complexity 
of the formula $\phi$.  Since all the variables are to be 
given the same type, no conflict can arise between different
occurrences of a variable; hence we need consider only
atomic formulae $\phi$.  These have the form $p=q$ for 
terms $p$ and $q$.  We can replace $p=q$ by 
$\exists u\,(p=u \land q=u)$, so we need only atomic
formulae $p=u$.  These we prove stratifiable by induction
on the complexity of the term, which is either $q \oplus  r$
or $\s(q)$ (often written $q^{\, \prime}$ in PA). 
We omit the details, which are technical but typical of 
interpretation proofs.
\smallskip

\noindent{\em Examples}.  In the rest of this
 paper we have occasion to prove several theorems or lemmas by 
 induction in NF.  To prove something by induction in NF
we have to check that the formula being proved is 
stratified. The  theorems are all special cases of 
the preceding lemma.  
Some formulas to which we apply 
Lemma~\ref{lemma:stratification} to obtain these 
formulas are as follows:
\begin{eqnarray*}
0\oplus x = x && \mbox{\qquad Lemma~\ref{lemma:zeroplusx}} \\ 
\forall x\,(x \oplus  \s(n) = \s (x) \oplus  n)&& \mbox{\qquad 
Lemma~\ref{lemma:ChurchSuccessorShift}} \\
x \neq 0 \imp 0 < x&& \mbox{\qquad Lemma~\ref{lemma:zerosmall} }\\
x < y \imp \s  x < y \ \lor \ \s  x = y && \mbox{\qquad Lemma~\ref{lemma:order2}} 
\end{eqnarray*}

\section{Properties of addition}\label{section:additionproperties}

For the rest of the paper, it does not matter how addition was defined;
we use only that it is defined by a stratified homogeneous formula and  satisfies the two formulas in Lemmas~\ref{lemma:ChurchAddition_equation} and \ref{lemma:ChurchZero_equation}, namely 
\begin{eqnarray*}
x\oplus 0 &=& x \mbox{ \qquad and} \\
x \oplus  \s n &=& \s(x\oplus n).
\end{eqnarray*}

Indeed one can easily prove that if $x \oplus  y$ is another 
function satisfying these properties then $x \oplus  y = x \oplus  y$
on Church numbers $x,y$.  Above we gave a set-theoretical definition
of addition, in Definition~\ref{definition:Churchaddition};  and 
a definition closer to $\lambda$-calculus in spirit, in 
\sect\ref{section:churchaddition}.
  The former is defined
only on Church numbers, while the latter can add any two functions; 
but as just remarked, they {\em necessarily} agree on Church numbers.
In this section we develop further properties of addition, using 
only the two properties listed above.

\begin{lemma} \label{lemma:zeroplusx} 
For $x \in \N$, $0 \oplus  x = x$.
\end{lemma}

\noindent{\em Proof}.  By
induction on $x$.  The base case is $0\oplus 0=0$, which follows
from $x\oplus 0=0$, which is Lemma~\ref{lemma:ChurchZero_equation} part (i). 
For the induction step, assume $0\oplus x=x$.  Applying successor
to both sides, we have $\s(0\oplus x) = \s x$.  By 
Lemma~\ref{lemma:ChurchAddition_equation}   we have 
$\s(0\oplus x)= 0 \oplus  \s x$.  Therefore $\s x = 0 \oplus  \s x$.
That completes the induction step,
and that completes the proof of the lemma.

\begin{lemma} \label{lemma:ChurchSuccessorShift}  
For  $x,n \in \N$,
 $$ x \oplus  \s n = \s x \oplus  n .$$
\end{lemma}

\noindent{\em Proof}.  
We quantify universally over Church numbers $x$, 
obtaining
$$ \forall x \in \N \,(x \oplus  \s n = \s x \oplus  n ),$$
and prove that by induction on $n$.  The formula
to be proved can be stratified by giving all variables
type 0.  
\smallskip

Base case:  $x \oplus  \s 0 = \s x \oplus  0$.  
\begin{eqnarray*}
x \oplus  \s 0 &=& \s(x\oplus 0)  \mbox{\qquad by Lemma~\ref{lemma:ChurchAddition_equation}}\\
&=& \ s x    \mbox{\qquad \qquad \ by Lemma~\ref{lemma:ChurchZero_equation}} \\
&=& \s(x) \oplus  0   \mbox{\qquad by Lemma~\ref{lemma:ChurchZero_equation}}
\end{eqnarray*}
That completes the base case.
\smallskip

Induction step: 
\begin{eqnarray*}
x \oplus  \s \s n &=& \s (x \oplus  \s n) \mbox{\qquad
 by Lemma~\ref{lemma:ChurchAddition_equation}}\\
&=& \s (\s x \oplus  n)   \mbox{\qquad by the induction hypothesis}\\
&=&  \s x \oplus  \s n     \mbox{\qquad\  by 
Lemma~\ref{lemma:ChurchAddition_equation},  with $x$ replaced by $\s x$}
\end{eqnarray*}
The replacement of $x$ by $\s x$ in the last step is legal, because
the statement being proved by induction is universally quantified over $x$.
That completes the proof of the theorem. 

\begin{lemma} \label{lemma:ChurchAdditionMaps}
$\forall x,y \in \N\ (x \oplus  y \in \N).$
\end{lemma}

\noindent{\em Proof}.  By induction on $y$, which is legal since 
the formula is stratified.   We omit the straightforward proof.

\begin{lemma}[Associativity] \label{lemma:ChurchAdditionAssociative}
$\forall x,y,z \in \N,((x\oplus y)\oplus z = x \oplus  (y\oplus z))$. 
\end{lemma}

\noindent{\em Proof}.  
 By induction on $y$, which is legal since the formula is stratified.
 \smallskip
 
  Base case: $(x\oplus 0)\oplus  z = x \oplus  z$
and $x \oplus  (0\oplus z) = x \oplus  z$, by Lemma~\ref{lemma:zeroplusx}.  Hence
$(x\oplus 0)\oplus z = x \oplus  (0 \oplus  z)$, completing the base case.
\smallskip

Induction step: 
\begin{eqnarray*}
(x \oplus  \s y) \oplus  z = \s(x\oplus y) \oplus  z &&\mbox{\qquad by Lemma \ref{lemma:ChurchAddition_equation}} \\
= (x\oplus y) \oplus  \s z &&\mbox{\qquad by Lemma~\ref{lemma:ChurchSuccessorShift}}\\
= \s( (x\oplus y)\oplus z) &&\mbox{\qquad by 
         Lemmas \ref{lemma:ChurchAddition_equation} and \ref{lemma:ChurchAdditionMaps}}\\
= \s (x \oplus  (y\oplus z)) &&\mbox{\qquad by the induction hypothesis} \\
= x \oplus  \s(y\oplus z) &&\mbox{\qquad by Lemmas~\ref{lemma:ChurchSuccessorShift}} \\
= x \oplus  (y \oplus  \s z) &&\mbox{\qquad by Lemma \ref{lemma:ChurchAddition_equation} and \ref{lemma:ChurchAdditionMaps}}\\
= x \oplus  (\s y \oplus  z ) &&\mbox{\qquad  by Lemma~\ref{lemma:ChurchSuccessorShift}} 
\end{eqnarray*}
That completes the proof of the lemma.
\smallskip

\begin{lemma}[Commutativity] \label{lemma:ChurchAdditionCommutative}
$\forall x,y  \in \N\,(x \oplus  y = y \oplus  x)$.
\end{lemma}

\noindent{\em Proof}. By induction on $y$, which is legal since the formula is 
stratified.
\smallskip

Base case,  $x\oplus \ChurchZero = \ChurchZero \oplus x$.  We have
\begin{eqnarray*}
 x\oplus \ChurchZero=x  && \mbox{\qquad by  Lemma~\ref{lemma:ChurchZero_equation}} \\
\ChurchZero\oplus x=x   && \mbox{\qquad  by Lemma~\ref{lemma:zeroplusx}} \\  
x \oplus \ChurchZero  = \ChurchZero  \oplus x &&\mbox{\qquad by the previous two lines}
\end{eqnarray*}
\smallskip

Induction step:
\begin{eqnarray*}
x \oplus  \s y &=& \s(x\oplus y) \mbox{\qquad\ by  Lemma~\ref{lemma:ChurchAddition_equation}} \\
&=& \s(y\oplus x) \mbox{\qquad\ by the induction hypothesis} \\
&=& y \oplus  \s x \mbox {\qquad\quad  by  Lemma~\ref{lemma:ChurchAddition_equation}}\\
&=& \s y \oplus  x \mbox{\qquad\quad  by Lemma~\ref{lemma:ChurchSuccessorShift}}
\end{eqnarray*}
That completes the induction step, and the proof of the lemma.

\begin{lemma} \label{lemma:doubleiteration} Let $f \in \FUNC$ and $f:X \to X$.
Then for Church numbers $j$ and $\ell$, and $x \in X$, we have
$$(jf)(\ell f x) = (j\oplus \ell)f x.$$
\end{lemma}

\noindent{\em Proof}.  The formula to be proved is stratified, 
so we may prove it by induction on $j$.  Base case: 
\begin{eqnarray*}
 0f(\ell fx) &=& \ell fx \mbox{\qquad\qquad by definition of 0} \\
 &=& (0 \oplus  \ell) fx \mbox{\quad \,  by Lemma~\ref{lemma:zeroplusx}}
 \end{eqnarray*}
 
 Induction step: 
\begin{eqnarray*}
(\s j f)(\ell fx) &=& f(j f(\ell fx)) \mbox{\quad\ \ by Theorem~\ref{theorem:successorequation}}\\
&=& f((j\oplus \ell)f x)  \mbox{\quad by the induction hypothesis} \\
&=& \s(j\oplus \ell)f x  \mbox{\qquad by Theorem~\ref{theorem:successorequation}}\\
&=& (j \oplus  \s \ell) f x \mbox{\qquad by Lemma~\ref{lemma:ChurchAddition_equation}}\\
&=& (\s j \oplus  \ell) \mbox{\qquad \ by Lemma~\ref{lemma:ChurchSuccessorShift}}
\end{eqnarray*}
That completes the proof of the lemma.

\section{Order on $\N$}

\begin{definition}\label{definition:orderN}
Order on the Church numbers is defined by 
$$ x < y \iff  \exists n \in \N\,(x\oplus \s n = y).$$
$$ x \le y \iff \exists n \in \N\,(x\oplus  n = y).$$
\end{definition}
These formulas are stratifiable, giving
$x$, $y$, and $n$ all index 0. ($\N$ is a parameter.) 
Therefore the relations
$x < y$  and $x \le y$ are definable in INF as sets of ordered pairs. 
\smallskip

\noindent{\em Remark.} We use the same symbols for these relations as
are used in \cite{inf-basics} for order on finite Frege cardinals; in our
formalization, we used different symbols, but for human readers, we think it 
 better not to introduce a new symbol.%
\footnote{Life is short, but the alphabet is shorter.  And the alphabet
of binary ordering relations is even shorter, and we later need $\preceq$ 
and $\prec$ for something else.}
\smallskip

 \begin{lemma} \label{lemma:notless0}  
For all $x\in \N$,  $ x \not < \ChurchZero$.
\end{lemma}

\noindent{\em Proof}.  Suppose $x \in \N$ and $x < \ChurchZero$.  Then
\begin{eqnarray*}
x \oplus \s n = \ChurchZero && \mbox{\qquad by definition of $<$}\\
\s(x \oplus n) = \ChurchZero && \mbox{\qquad by Lemma~\ref{lemma:ChurchSuccessorShift}}\end{eqnarray*}
But that contradicts Theorem~\ref{theorem:successoromitszero}, which
says that $\ChurchZero$ is not a successor.

\begin{lemma} \label{lemma:order1} For all $x$ and $y$ in $\N$,
$$x = y \ \lor \ x < y \imp x < \s y $$
\end{lemma}

\noindent{\em Proof}. 
Suppose $x = y \lor x < y$.  
\smallskip

Case~1,  $x=y$.  We must prove $x < \s x$. 
\begin{eqnarray*}
x \oplus  \s 0 = \s x \oplus  0  &&\mbox{\qquad by Lemma~\ref{lemma:ChurchAddition_equation}}\\
 = \s x &&\mbox{\qquad by Lemma~\ref{lemma:ChurchZero_equation}}\\
 x \oplus  \s 0 = \s x &&\mbox{\qquad by the preceding two lines}\\
 x < \s x  && \mbox{\qquad by Definition~\ref{definition:orderN}}
 \end{eqnarray*}
That completes Case~1. 
\smallskip 
 
Case~2, $x < y$.  We have to prove $x < \s y$.
\begin{eqnarray*}
x \oplus  \s t = y  && \mbox{\qquad for some $t \in \N$, by  Definition~\ref{definition:orderN}}\\
 \s(x\oplus  \s t) = \s y  && \mbox{\qquad by the preceding line}\\
 x \oplus  \s ( \s t) = \s y && \mbox{\qquad by Lemma~\ref{lemma:ChurchAddition_equation}}\\ 
x < \s y  &&\mbox{\qquad by   Definition~\ref{definition:orderN}}
\end{eqnarray*}
That completes Case~2.  
 That completes the proof of the
lemma.

\begin{corollary} \label{lemma:xlessthansx}
For all $x \in \N$, $x < \s x$.
\end{corollary}

\noindent{\em Remark}.  This does not guarantee $x \neq \s x$
since we do not have trichotomy.
\medskip

\noindent{\em Proof}. Take $x=y$ in Lemma~\ref{lemma:order1}.

\begin{lemma} \label{lemma:zerosmall} For all $x \in \N$, 
$x \neq \ChurchZero \imp \ChurchZero < x$. 
\end{lemma}

\noindent{\em Proof}.   By induction on $x$.  The formula to be proved
is stratifiable, by Lemma~\ref{lemma:stratification}.

We proceed with the induction.  The 
base case is immediate (since $\ChurchZero \neq \ChurchZero$ implies anything).
  To prove the induction step, we have to prove 
$\s  x \neq \ChurchZero \imp \ChurchZero < \s  x$.  Suppose $\s x \neq \ChurchZero$;
we have to prove $\ChurchZero < \s x$.   By Lemma~\ref{lemma:order1},
it suffices to prove $\ChurchZero < x \lor \ChurchZero = x$. But that follows from 
the induction hypothesis $x \neq \ChurchZero \imp \ChurchZero < x$, even with 
intuitionistic logic, because by Lemma~\ref{lemma:decidable0},
 $x \neq \ChurchZero \lor x = \ChurchZero$, and
if $x \neq \ChurchZero$ then $\ChurchZero < x$, while if $x=\ChurchZero$ then $\ChurchZero = x$.  That completes
the proof of the lemma.

\begin{lemma} \label{lemma:order2}
For all $x,y \in \N$,
$$ x < y \imp \s  x < y \ \lor \ \s  x = y $$
\end{lemma}

\noindent{\em Proof}. 
\begin{eqnarray*}
x < y && \mbox{\qquad assumption}\\
x\oplus  \s p = y  && \mbox{\qquad for some $p \in \N$, by definition of $<$}\\
x \oplus  \s p =  \s x \oplus  p = y && \mbox{\qquad  by Lemma~\ref{lemma:ChurchSuccessorShift}}\\
p = \ChurchZero \ \lor \ p \neq \ChurchZero && \mbox{\qquad by Lemma~\ref{lemma:decidable0}}
\end{eqnarray*}
 If $p = \ChurchZero$
then $\s x = y$ and we are done. If $ p \neq \ChurchZero$ then 
\begin{eqnarray*}
p = \s \ell &&\mbox{\qquad for some $\ell \in \N$, by Lemma~\ref{lemma:predecessor}}\\
x \oplus \s (\s \ell) = y && \mbox{\qquad since $x \oplus \s p = y$}\\
 \s x \oplus  \s \ell = y  &&\mbox{\qquad by Lemma~\ref{lemma:ChurchSuccessorShift}}\\
 \s x < y && \mbox{\qquad by definition of $<$}
 \end{eqnarray*}
That completes the proof of the lemma.

\begin{lemma} \label{lemma:Churchletolessthan} For $x,y \in \N$, if $x \le y$ and $x \neq y$
then $x < y$.
\end{lemma}

\noindent{\em Proof}.  Suppose $x,y \in \N$ and $x \le y$. Then for some $m$ we have
$x \oplus  m = y$.  If $m = \ChurchZero$ then $ x= y$, by Lemma~\ref{lemma:ChurchZero_equation}.
Hence $m \neq \ChurchZero$.  Then by Lemma~\ref{lemma:predecessor},
$m = \s r$ for some $r \in \N$.  Then $x \oplus  \s r = y$.  Then by definition of $<$,
we have $x < y$.  That completes the proof of the lemma.

\begin{lemma} \label{lemma:Churchle} For $x,y \in \N$,   
$$x \le y  \iff x < y \ \lor \ x =y.$$
\end{lemma}

\noindent{\em Proof}. Left to right: Suppose $x \le y$.  Then $x \oplus  m = y$ for some $m \in \N$.
By Lemma~\ref{lemma:decidable0}, $m = \ChurchZero \ \lor \ m \neq \ChurchZero$.  If $m = \ChurchZero$,
then by Lemma~\ref{lemma:ChurchZero_equation}, $x = y$.  If $m \neq \ChurchZero$,
then by Lemma~\ref{lemma:predecessor}, $m = \s r$ for some $r \in \N$.
Then $x < y$ by the definition of $<$.
\smallskip

Right to left.  Suppose $x < y \ \lor x = y$.  If $x < y $ then $x \oplus  \s m = y$ for some $m \in \N$.
Then $\s m \in \N$ by Lemma~\ref{lemma:successorN}, so $x \le y$ by definition of $\le$.
If $x = y$ then $x \oplus  \ChurchZero = y$, by Lemma~\ref{lemma:ChurchZero_equation}, so $x \le y$.
That completes the proof of the lemma.

\begin{lemma}[transitivity] \label{lemma:transitivity}
$x < y$ is a transitive relation.  That is,
for $x,y,z \in \N$, 
$$ x < y \land y < z \imp x < z.$$
\end{lemma}

\noindent{\em Proof}.    Suppose $x < y$ and $y < z$.  
Then for some $p,q$ we have $x\oplus p = y$ and $y \oplus  q = z$. 
Then $(x\oplus p)\oplus q = z$.  By the associativity 
of addition we have $x\oplus  (p\oplus q) = z$.  Then $x < z$. 
That completes the proof
of the lemma. 

\begin{lemma}\label{lemma:trichotomy1}
For $x,y \in \N$,
$$x < y \ \lor \ x = y \ \lor \ y < x$$ 
\end{lemma}

\noindent{\em Remark}. We do not claim that exactly one
of the three alternatives holds.
\medskip

\noindent{\em Proof}. We proceed
by induction on $x$.
  When $x=\ChurchZero$
we have to prove 
$$\ChurchZero < y \ \lor \ \ChurchZero = y\  \lor\  y < \ChurchZero.$$  We have
\begin{eqnarray*}
y = \ChurchZero \ \lor \ y \neq \ChurchZero  &&\mbox{\qquad by Lemma~\ref{lemma:decidable0}}
\end{eqnarray*}
If $y = \ChurchZero$, we are done.  If $y \neq  \ChurchZero$ then
\begin{eqnarray*}
\ChurchZero < y\, \lor\, \ChurchZero = y &&\mbox{\qquad by  Lemma~\ref{lemma:zerosmall}}
\end{eqnarray*}
That completes the base case.
\smallskip 

For the induction step, we assume 
\begin{eqnarray*}
 x < y \ \lor \ x = y \ \lor \ y < x 
\end{eqnarray*}
and must prove
$$\s  x < y \ \lor \ \s  x = y \ \lor \ y < \s  x.$$
We argue by cases.
\smallskip

Case~1, $x < y$.  Then by Lemma~\ref{lemma:order2}, $\s  x < y$
or $\s  x = y$.  That completes Case~1.
\smallskip

Case~2, $x = y$. Then by Corollary~\ref{lemma:xlessthansx}, $y < \s  x$.
That completes Case~2.
\smallskip

Case~3, $y < x$.  
 Then
\begin{eqnarray*}
 x < \s  x  &&\mbox{\qquad by Corollary~\ref{lemma:xlessthansx}}\\
 y < \s x   && \mbox{\qquad by Lemma~\ref{lemma:transitivity}}
 \end{eqnarray*}
That completes Case~3. 
That completes the proof of the lemma.

\section{Structure of $\N$ under successor:  The picture}
In this  and the following sections, we explore the consequences of the assumption 
that Church successor is not one to one.  We first attempt to convey an intuitive 
picture of the situation.

\begin{figure}[ht]
\caption{The stem $\Stem$, the loop $\L$, and unique double successor}
\label{figure:rho2}
\FigureRho
\end{figure}
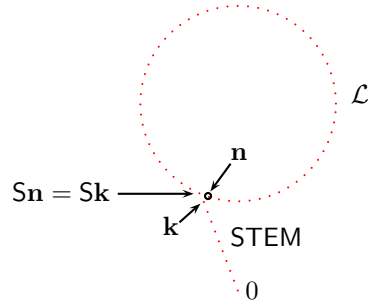
Figure~\ref{figure:rho2} (already shown in the introduction, but reprinted
here for convenience) illustrates the structure of $\N$ under
successor.  To arrive at this figure, imagine coloring $\ChurchZero$ red, and 
at each stage where you have just colored $x$ red, then color $\s x$
red unless $\s (\s x)$ is already red.  Then stop.  Let $\n$ be the last
number you encountered.  
You will have
colored every integer red except $\n$ (shown black in the figure).  The
reason you did not color $\n$ is that $\s \n = \s \k$, where $\k$ is some 
number that you already colored red.  We call $\s n$ a ``double successor.''
\smallskip

We emphasize that at this point we have not proved that this figure is 
accurate.  There might be many more double successors not shown; imagine
a gray spiderweb of mysterious Church numbers, merging at different
places into the red part of the figure.  But the red part, if it 
could be defined, contains 0 and is closed under successor, so it 
intuitively should be all of $\N$.  We shall prove in 
Theorem~\ref{theorem:rho} below  
that, at least if $\N$ is assumed to be finite, this picture 
is an accurate one. The part that you colored before reaching $\k$ 
(and including $\k$) is called the ``stem''.   The rest of the red numbers
(plus $\n$) comprise ``the loop.''   The next several sections will show in
detail that this picture is correct.  

\begin{definition} \label{definition:doublesuccessor}  $p$ is 
 not a double successor if 
$$ \forall a,b\, (a \in \N \ \land \ b \in \N \ \land \s a = \s b = p  \imp \ a = b).$$
\end{definition}

It might seem more natural to define the concept this way: 
 $p$ is a double successor if there exists $a,b \in \N$ with $a \neq b$ and  $\s a = \s b =p$. 
But negating this introduces a double negation, which we prefer not to have.  Hence
the definition above.   Soon we will be working under the hypothesis that $\N$ is finite,
which implies that $\N$ has decidable equality,  making this double negation irrelevant.
Also, we could strengthen the notion by dropping the condition $b \in \N$; we will do that
in one place below.  

\begin{lemma} \label{lemma:lessthansuccessorN}
If $x,y \in \N$ and 
 $\s x$ is not a double successor, then 
 $$ y < \s x \iff y< x \mspace{1mu} \lor \mspace{2mu} y = x.$$
Explicitly this means
$$ \forall x,y \in \N\,(( \forall u \in \N\, (\s u = \s x \imp u = x)) \imp 
y < \s x \iff y< x \mspace{1mu} \lor \mspace{2mu} y = x).$$
\end{lemma}

\noindent{\em Proof}.  
  Left to right:
\begin{eqnarray*}
y < \s x  && \mbox{\qquad assumption}\\
y \oplus \s p = \s x && \mbox{\qquad for some $p \in \N$, by definition of $<$}\\
\s(y \oplus p) = \s x && \mbox{\qquad by Lemma~\ref{lemma:ChurchAddition_equation}}\\
y \oplus p = x && \mbox{\qquad since $\s x$ is not a double successor}\\
p = \ChurchZero \ \lor \ p \neq \ChurchZero &&\mbox{\qquad by Lemma~\ref{lemma:decidable0}}
\end{eqnarray*}
Case~1, $p = \ChurchZero$.  Then $y = x$.  That completes Case~1.
\smallskip

Case~2, $p \neq \ChurchZero$.  Then
\begin{eqnarray*}
p = \s m && \mbox{\qquad for some $m \in \N$, by Lemma~\ref{lemma:predecessor}}\\
y \oplus \s m = x && \mbox{\qquad since $y \oplus p = x$}\\
y < x  && \mbox{\qquad by definition of $<$}
\end{eqnarray*}  
That completes Case~2. 
  That completes the left-to-right implication.

Right-to-left:  Suppose $y < x \lor y = x$, and $\s x$
is not a double successor.  We have to prove $y < \s x$.
We argue by cases.
\smallskip

Case~1, $y < x$.  Then 
\begin{eqnarray*}
x < \s x && \mbox{\qquad by Corollary~\ref{lemma:xlessthansx}}\\
y < \s x && \mbox{\qquad by transitivity}
\end{eqnarray*}
That completes Case~1.
\smallskip

Case~2, $y = x$. 
Then $y < \s x$ by Corollary~\ref{lemma:xlessthansx}.
That completes Case~2.  That completes the right-to-left direction.
That completes the proof of the lemma.

\section{Structure of $\N$ under successor:  The stem}

\begin{definition} \label{definition:stem}
The set $\Stem$  is the intersection of all subsets of $\N$ 
containing 0 and closed under successors that are not double successors.
More precisely, $ \Stem$ is the intersection of $\N$ and all $X$ such that 
\begin{eqnarray*}
 0 \in X \ \land\ \forall u \in \N\,(u \in X  \ \land \  \forall v \in \N\,
(\s v = \s  u \imp v = u) \imp \s u \in X). 
\end{eqnarray*}
\end{definition}
The intention of the definition is that
 $\Stem$ should contain everything from $\ChurchZero$ up to but not including
the first double successor.

\begin{lemma} \label{lemma:SN} $\Stem \subseteq \N$.
\end{lemma}

\noindent{\em Proof}. Immediate from the definition of $\Stem$ as the intersection
of $\N$ with some other sets.

\begin{lemma} \label{lemma:S1} $\Stem$  is one of the sets
used to define $\Stem$.  That is,  
$$\ChurchZero \in \Stem \ \land\ \forall u \in \N\,(u \in \Stem \land \forall v \in \N\,
(\s v = \s  u \imp v = u) \imp \s u \in \Stem).$$
\end{lemma}

\noindent{\em Proof}. 
Let $X$ satisfy the formula in the lemma (with $\Stem$ replaced by $X$).
  Then $\ChurchZero \in X$.
Since $X$ was arbitrary, $\ChurchZero \in \Stem$.  Now suppose $u \in \N$ and $u \in \Stem$
and $\forall v\in \N\, (\s v = \s u \imp v = u)$.  Then   $\s u \in X$.
Since $X$ was arbitrary, $\s u \in \Stem$, by Definition~\ref{definition:stem}.
  That completes the proof of the lemma.

\begin{lemma} \label{lemma:Soneone} 
Church successor is one-to-one on $\Stem$.  What is more, 
 $$ \forall u \in \Stem\, ( \s u \in \Stem \imp \forall v\in \N\,( \s u = \s v \imp u = v)).$$
 That is, there are no double successors in $\Stem$.
 \end{lemma} 
\medskip

\noindent{\em Remark.} ``What is more'' because $v$ is not required to be in $\Stem$.
\medskip
 
\noindent{\em Proof}.  Define $X$ to be 
\begin{eqnarray}
 X:= \{   p \in \N : \forall u, v\in \N\,(  p =  \s u \imp \s u = \s v  \imp u = v)\}. \label{eq:E1661}
 \end{eqnarray}
 The formula is stratified, giving all the variables index 0; $\N$ is a parameter.
 Hence the definition can be given in INF. 
Then $X \subseteq \N$.  By Theorem~\ref{theorem:successoromitszero},
$\ChurchZero \in X$.   $X$ is closed under successors except double successors; 
that is, if $x \in X$ and $\forall v \in \N\,(\s x= \s v \imp x = v)$,
then $\s x \in X$, as we see by putting $p = \s x$ in the definition of $X$.
(By the hypothesis $x \in X$, we have $x \in \N$.) 
Therefore, by the definition of $\Stem$, we have $\Stem \subseteq X$.
\smallskip

Suppose $S u \in \Stem$ and $\s u  = \s v$.   Then $\s u \in X$, since $\Stem \subseteq X$.
Therefore $u=v$.
That completes the proof of the lemma.

\begin{lemma} \label{lemma:Spred}
  If   $y \in \N$ and  $\s y \in \Stem$, then $y \in \Stem$.
\end{lemma}

\noindent{\em Proof}.
By Lemma~\ref{lemma:Soneone} there are no double successors in $\Stem$, so
it suffices to show that every nonzero element of $\Stem$ is the successor of 
something in $\Stem$.  Let $X$  be the set of elements of $\Stem$
that are equal to $\ChurchZero$ or are successors of something in $\Stem$.  Explicitly
$$ X = \{ x: x \in \Stem \ \land \ x = \ChurchZero \ \lor\ \exists y\,( \s y = x \ \land \ y \in \Stem)\}.$$
The formula is stratified, giving $x$ and $y$ index 0, with $\Stem$ as a parameter.
I say that $X$  is closed under successors that are not double successors.   Let
$x \in X$ and suppose 
\begin{eqnarray}
\forall v\in \N\,(\s x = \s v \imp x = v) \label{eq:E1693}
\end{eqnarray}
 (informally, $\s x$ is 
not a double successor).   We must show 
$\s x \in X$.  Since $x \in X$,  $x \in \Stem$.  To show $\s x \in X$ we must show two things:
\begin{eqnarray}
&& \s x \in \Stem  \label{eq:E1695}  \\
&&  \exists y\, (\s y = \s x \ \land \ y \in \Stem)  \label{eq:E1696}
\end{eqnarray}
(\ref{eq:E1696}) is immediate, taking $y=x$.  To verify (\ref{eq:E1695}) we use
that $x \in \Stem$ and $\s x$ is not a double successor (\ref{eq:E1693}).
By Lemma~\ref{lemma:S1}, $\Stem$ is closed under successors except double successors,
so $\s x \in \Stem$ as desired.  That completes the proof that $X$ is closed under
successors except double successors.   Then by the definition of $\Stem$, we have 
$\Stem \subseteq X$.   
\smallskip

Now suppose $\s x \in \Stem$ and $x \in \N$; we must prove $x \in \Stem$.  Since $\Stem \subseteq X$,
we have $\s x \in X$.  By definition of $X$, 
$$ \s x = \ChurchZero \ \lor \ \exists y\,(y \in \Stem \ \land \ \s y = \s x).$$
By Theorem~\ref{theorem:successoromitszero}, and the hypothesis $x \in \N$,
we have $\s x \neq \ChurchZero$.   Therefore, for some $y \in \Stem$, we have $\s y = \s x$.
By Lemma~\ref{lemma:Soneone}, we have $y = x$.  Since $y \in \Stem$ and $y = x$, we have 
$x \in \Stem$ as desired. 
That completes the proof of the lemma.

\begin{lemma} \label{lemma:Sdecidable} $\Stem$ has decidable equality.  In fact,
\begin{eqnarray*}
 \forall x \in \N\, \forall y\,( x \in \Stem \imp  y \in \N \imp x = y \ \lor\  x \neq y).
\end{eqnarray*}
\end{lemma}

\noindent{\em Remark}. It is not necessary to assume $y \in \Stem$.
\medskip 

\noindent{\em Proof}. 
 We prove by induction on $x$ that 
\begin{eqnarray}
 \forall y( x \in \Stem \imp y \in \N \imp x = y \ \lor\  x \neq y).
 \label{eq:E2132}
\end{eqnarray}
That formula is stratified, so it is legal to prove it by induction.
The base case follows from Lemma~\ref{lemma:decidable0}.
For the induction step, suppose $\s x \in \Stem$ and $y \in \N$; 
we have to prove $\s x = y \ \lor\  \s x \neq y$. 
By Lemma~\ref{lemma:decidable0}, we may argue by cases
according as $y = \ChurchZero$ or not.  If $y = \ChurchZero$, we 
are done by Lemma~\ref{lemma:decidable0}.  If $y \neq \ChurchZero$,
then $y = \s q$ for some $q$.  Then 
\begin{eqnarray*}
x \in \Stem && \mbox{\qquad by Lemma~\ref{lemma:Spred}}\\
\s x = y \ \lor \s x \neq y   \iff   \s x = \s q \lor \s x \neq \s q 
        && \mbox{\qquad since $y = \s q$} \\
 \iff  x = q \ \lor x \neq q &&\mbox{\qquad by Lemma~\ref{lemma:Soneone}, since $\s x \in \Stem$} 
\end{eqnarray*}
and that  follows from the induction hypothesis (\ref{eq:E2132}).
That completes the proof  of the lemma.

\begin{lemma} \label{lemma:Sinit} Suppose $y \in \Stem$ and $x \in \N$ and $x < y$.  Then 
$x \in \Stem$. 
\end{lemma}

\noindent{\em Proof}.  By induction on $y$ we prove
$$ y \in \N \imp y \in \Stem \imp \forall x\in \N\, (x < y \imp x \in \Stem).$$
That formula is stratified, so induction is legal.
\smallskip

Base case:  When $y = \ChurchZero$, it is impossible that $x < y$, by 
Lemma~\ref{lemma:notless0}. Therefore $x < \ChurchZero \imp x \in \Stem$.
That completes the base case.
\smallskip

Induction step: Suppose $\s y \in \N$ and $\s y \in \Stem$ and $x < \s y$.
We must prove $x \in \Stem$.  We have
\begin{eqnarray*}
y \in \Stem  && \mbox{\qquad  by Lemma~\ref{lemma:Spred}}\\
\s y \in \Stem && \mbox{\qquad by hypothesis}\\
 \s y \mbox{\ is not a double successor} && \mbox{\qquad by Lemma~\ref{lemma:Soneone}}\\
x < y \ \lor \ x = y  && \mbox{\qquad by  Lemma~\ref{lemma:lessthansuccessorN}}
\end{eqnarray*}
If $x < y$, then by the induction hypothesis, $x \in \Stem$.   If $x = y$
then $x \in \Stem$ because $y \in \Stem$.  That completes the induction step.
That completes the proof of the lemma.
   
\begin{lemma} \label{lemma:Smax}
Suppose  $\s k = \s n$ with $k \in \Stem$ and $n \in \N$ and $k \neq n$.
Then $k$ is a maximal element of $\Stem$; more precisely,  
$$ \Stem = \{ x\in \N: x < k \ \lor \ x = k\}.$$
\end{lemma}

\noindent{\em Proof}.
Define 
$$Z = \{x \in \Stem : x < k \ \lor\ x = k\}.$$ 
The formula is stratified, since $<$ is definable as a relation in INF.
I say that $Z$ contains $\ChurchZero$ and is closed under successor except
double successors.
\smallskip

To prove $\ChurchZero \in Z$: 
\begin{eqnarray*} 
 \ChurchZero \in \Stem  &&\mbox{\qquad by Lemma~\ref{lemma:S1}}\\
k = \ChurchZero = k \ \lor \ k \neq \ChurchZero  &&\mbox{\qquad by Lemma~\ref{lemma:decidable0}}
\end{eqnarray*}
We argue by cases.  
\smallskip

Case 1, $k = \ChurchZero$.  Then $\ChurchZero \in Z$, by definition of $Z$.
\smallskip

Case 2, $\k \neq \ChurchZero$.  Then 
\begin{eqnarray*}
k = \s m && \mbox{\qquad for some $m \in \N$, by Lemma~\ref{lemma:predecessor}}\\
\ChurchZero \oplus  \s m = \s m && \mbox{\qquad by Lemma~\ref{lemma:zeroplusx}}\\
\ChurchZero \oplus  \s m = k  && \mbox{\qquad since $k = \s m$}\\
\ChurchZero < k && \mbox{\qquad by definition of $<$}
\end{eqnarray*}

  To prove $Z$ is 
closed under successor except double successor:
Suppose $x \in Z$ and $\s x$ is not a double successor.
Since $x \in Z$, $x < k \ \lor\ x = k$.  But $\s k$ is
a double successor, so $x \neq  k$.  Therefore $x < k$.   
Then $\s x = k \lor \s x < k$,  by Lemma~\ref{lemma:order2}.
Hence $\s x \in Z$.  
Therefore $Z$ contains $\ChurchZero$ and is closed 
under successor except double successors.  Therefore $\Stem \subseteq Z$.  
Now I say
$$\Stem = \{ x \in \N : x < k \lor x = k\}.$$
It suffices to prove
$$ x \in \Stem \iff x \in \N \ \land \ ( x < k \lor x = k).$$
Left to right:  Suppose $x \in \Stem$.  By Lemma~\ref{lemma:SN}, $\Stem \subseteq \N$,
so $x \in \N$.  Since
$\Stem \subseteq Z$, we
have $x \in Z$.  Then  $x < k \lor x = k$, by definition of $Z$.
\smallskip

Right to left:  Suppose $x \in \N \ \land \ ( x < k \lor x = k)$.
Since $k \in \Stem$, we have $x \in \Stem$ by Lemma~\ref{lemma:Sinit}.
That completes the proof of the lemma.

\begin{lemma} \label{lemma:trichotomyonS}
Let $P$ be any subset of $\N$ satisfying the following two conditions:
\smallskip

(i)  $\forall x \in \N\,( S x \in P \imp x \in P)$.
\smallskip

(ii) $\forall x,y \in P (y < \s x \iff y< x \mspace{1mu} \lor \mspace{2mu} y = x)$.
\smallskip

Then trichotomy holds on $P$. That is, for $x,y \in P$, exactly 
one of $x < y$, $x = y$, or $y < x$ holds.
\end{lemma}

\noindent{\em Proof}.  Assume (i) and (ii). 
By Lemma~\ref{lemma:trichotomy1}, at least one of the 
three alternatives (of trichotomy) holds.  We prove by induction on $y$ that 
$$ y\in P \imp \forall x \in P\, ( \neg(x < y \ \land \ y < x) \ \land \ x \not< x).$$
The formula is stratified, giving $x$ and $y$ index 0 and $P$ index 1, so 
we may proceed by induction.
\smallskip

 Base case, $y=\ChurchZero$.  Suppose $x \in P$ and $\ChurchZero \in P$.   We do not have
$x < \ChurchZero$, by Lemma~\ref{lemma:notless0}. Suppose $x=\ChurchZero \ \land \ \ChurchZero < x$.
Then $\ChurchZero < \ChurchZero$, contradicting Lemma~\ref{lemma:notless0}.  That 
completes the base case. 
\smallskip

Induction step:  Suppose $\s y \in P$.  Assume $x\in P$.  We have to prove
\begin{eqnarray*}
 \neg\,(x < \s y \ \land \s y < x) \ \land \ \s y \not< \s y. \label{eq:E2113}
 \end{eqnarray*}
We have
\begin{eqnarray*}
y \in P  && \mbox{\qquad by (i), since $\s y \in P$}\\
x < \s y  && \mbox{\qquad assumption}\\
x < y \lor x = y  && \mbox{\qquad by (ii), with $x$ and $y$ switched}\\
x < y \imp x \oplus  j = y   && \mbox{ for some $j\in \N$, by the definition of $<$}\\
x = y \imp x \oplus j = y    && \mbox{ for $j = \ChurchZero$, by Lemma~\ref{lemma:ChurchZero_equation}}\\
x \oplus j = y && \mbox{ for some $j\in \N$, by the preceding three lines}\\
\end{eqnarray*}
Now assume that also $\s y < x$. Then arguing as above, but switching $x$ and $y$,
we have
\begin{eqnarray*}
  \s y \oplus  \ell = x && \mbox{\qquad for some $\ell \in \N$} \\
 \s y \oplus  \ell \oplus  j = x\oplus j  && \mbox{\qquad by the preceding line}\\
 \s y \oplus  \ell \oplus  j  = y  &&\mbox{\qquad since $x \oplus j = y$}\\
 y \oplus  \s(\ell \oplus  j) = y && \mbox{\qquad by Lemma~\ref{lemma:ChurchSuccessorShift}}\\
 y < y && \mbox{\qquad by the definition of $<$}
\end{eqnarray*}
But that contradicts the induction hypothesis, since $y \in P$.
We have now proved the first half of (\ref{eq:E2113}), namely
$$  \neg\,(x < \s y \ \land \s y < x).$$
Then by Lemma~\ref{lemma:trichotomy1} we have $\s y = x$.
It remains to prove $\s y \not< \s y$.  Suppose $\s y < \s y$.  Then
\begin{eqnarray*}
y < \s y && \mbox{\qquad by Corollary~\ref{lemma:xlessthansx}}\\
\s y \in P  && \mbox{\qquad by hypothesis}\\
\s y < y \ \lor \ \s y  = y && \mbox{\qquad by hypothesis (ii)} 
\end{eqnarray*}
We argue by cases accordingly.
\smallskip

Case~1, $\s y < y$.  Since $y < \s y$ we have $y < y$ by 
transitivity, contradicting the induction hypothesis.
\smallskip

Case~2, $\s y = y$.  Then since $y < \s y$ we again have $y < y$,
contradicting the induction hypothesis.
That completes the proof of (\ref{eq:E2113}.
That completes the induction step.  That 
completes  the proof of the lemma.

\begin{lemma}\label{lemma:knotlessthank}
Suppose there is a double successor
$\s \k = \s  \n$ with $\k \in \Stem$ and $ \n \neq \k$.
Then $\k \not < \k$.
\end{lemma}

\noindent{\em Proof}.   Let
$$ X:= \Stem - \{ x \in \N: \k < x\}.$$
I say that $X$ is closed under non-double successors.  Suppose $x \in X$ and 
$\s x$ is not a double successor.  Then
\begin{eqnarray*}
\s x \in \Stem  &&\mbox{\qquad  by Lemma~\ref{lemma:S1}}\\
\k \not< x  &&\mbox{\qquad since $x \in X$}
\end{eqnarray*}
 I say $  \k \not < \s x$.  Suppose $\k < \s x$.
Then by Lemma~\ref{lemma:lessthansuccessorN}, we have $\k < x \ \lor \ \k = x$.
We do not have $\k = x$, since $\s \k$ is a double successor but $\s x$ is not.  
Therefore $\k < x$, contradiction.  Therefore $X$ is closed under non-double 
successors, as claimed.  Therefore $\Stem \subseteq X$, by definition of $\Stem$.
But $\k < \k$ implies $\k \not \in X$,  while by hypothesis, $\k \in \Stem$.
That completes the proof of the lemma.

\begin{lemma} \label{lemma:kisunique}
Suppose there is a double successor
$\s \k = \s  \n$ with $\k \in \Stem$ and $ \n \neq \k$.
Let $x \in \Stem$ with $x \neq \k$.  Then $\s x$ is 
not a double successor; that is, $\forall u\in \N( \s u = \s x \imp u = x).$
\end{lemma}

\noindent{\em Proof}.
Suppose $\s \k$ is a double 
successor, and $x \in \Stem$ and $x \neq \k$. Then by Lemma~\ref{lemma:Smax},
we have $x< \k$.  If $\s x$ is a double successor, then by Lemma~\ref{lemma:Smax},
  $\k < x$.   Then by Lemma~\ref{lemma:transitivity}, $\k < \k$.
But that contradicts Lemma~\ref{lemma:knotlessthank}.  Hence,
$\s x$ is not a double successor.  That is, 
$$ \neg \exists u \in \N\,(\s u = \s x \ \land \ u \neq x).$$
By logic,
$$ \forall u \in \N\,(\s u = \s x \imp \neg\neg\,(u = x).$$
Since $x \in \Stem$ and $y \in \N$,  we have $\neg\neg\,u = x \imp u = x$,
by Lemma~\ref{lemma:Sdecidable}.  Therefore we can drop the double negation:
$$ \forall u \in \N\,(\s u = \s u \imp  u = x).$$
That completes the proof of the lemma.

\begin{lemma}\label{lemma:Strichotomy}
Suppose there is a double successor
$\s \k = \s \n$ with $\k \in \Stem$ and $\n \neq \k$.  Then
 for $y \in \Stem$,  we have $\neg(x < y \ \land y < x) \ \land \ \neg\,(y < y)$.
\end{lemma}

\noindent{\em Proof}.
We intend to apply Lemma~\ref{lemma:trichotomyonS},
with $P$ replaced by $\Stem$.   To do that, it suffices to verify 
the  hypotheses of Lemma~\ref{lemma:trichotomyonS}, namely
\begin{eqnarray*}
&&(i) \qquad   x \in \N \imp \s x \in \Stem  \imp x \in \Stem \\ 
&&(ii) \qquad x,y \in \Stem  \imp (y < \s x \iff y< x \mspace{1mu} \lor \mspace{2mu} y = x)  
 \end{eqnarray*}
Ad (i):  This is Lemma~\ref{lemma:Spred}.
\smallskip

Ad (ii):  Assume $x,y \in \Stem$.  By Lemma~\ref{lemma:Sdecidable}, we have $x = \k \ \lor \ x \neq \k$.  We
argue by cases.
\smallskip

Case 1, $x = \k$.  We have to prove $y < \s \k \iff y < \k \ \lor \ y = \k$.  
\smallskip

Left to right: By
Lemma~\ref{lemma:Smax}, the right side is equivalent to $y \in \Stem$, which we have assumed.
\smallskip

Right to left: Assume $y < \k \ \lor \ y = \k$; we have to prove $y < \s \k$.  We have
$\k < \s \k$ by Lemma~\ref{lemma:xlessthansx}.  If $y = \k$ we are done;  if 
$y < \k$ then by transitivity (Lemma~\ref{lemma:transitivity}) we have $y < \s \k$.
That completes Case~1.
\smallskip 

Case 2, $x \neq \k$.  By Lemma~\ref{lemma:kisunique}, $x$ is not a double successor.
Then by Lemma~\ref{lemma:lessthansuccessorN}, we have (ii). 
That completes Case~2.  That completes the proof of the lemma.

\begin{lemma}\label{lemma:nneqzero}
Suppose there is a double successor
$\s \k = \s \n$ with $\k \in \Stem$ and $\n \neq \k$.  Then $\n \neq \ChurchZero$.
\end{lemma}

\noindent{\em Proof}. Assume  $\n = \ChurchZero$.
 Then $\n \in \Stem$, by Lemma~\ref{lemma:S1}.
By Lemma~\ref{lemma:Smax},
$\ChurchZero$ is the maximal element of $\Stem$.  Then $\ChurchZero$ is the only 
element of $\Stem$, by Lemma~\ref{lemma:notless0}.
  Then $\k = \ChurchZero$, since $\k \in \Stem$. But 
that contradicts $\n \neq \k$. That completes the proof of the lemma.

\section{Structure of $\N$ under successor: the loop}

\begin{definition} \label{definition:loop}
Suppose there is a double successor
$\s \k = \s \n$ with $\k \in \Stem$ and $\n \neq \k$.  Then the {\bf loop} 
  $\L(\n)$ is the intersection of all sets $X$ containing $\n$ and 
closed under successor. 
\end{definition}

The formula is stratified, giving $\n$ and $\k$ both index 0.  $\Stem$ is a parameter.
Hence the definition is legal in INF.

\begin{lemma}\label{lemma:L1} Suppose there is a double successor
$\s \k = \s \n$ with $\k \in \Stem$ and $\n \neq \k$.   Then $\n \in \L(\n)$ and 
$\L(\n)$ is closed under Church successor.
\end{lemma}

\noindent{\em Proof}. Follows from the definition of $\L(\n)$ as the intersection of 
all sets $w$ that 
  contain $\n$ and are closed under successor.  Since $\n$ belongs to every such $w$,
  it belongs to their intersection.  Suppose $x \in \L(\n)$;  then $x\in w$, so $\s x \in w$.
  Then $\s x$ belongs to the intersection of all such $w$, i.e., $\s x \in \L(\n)$.
  That completes the proof of the lemma.
  
\begin{lemma} \label{lemma:LN} Suppose there is a double successor
$\s \k = \s \n$ with $\k \in \Stem$ and $\n \neq \k$.   Then $\L (\n) \subseteq \N$.
\end{lemma}

\noindent{\em Proof}. $\n \in \L n$ by Lemma~\ref{lemma:L1}.  Then $\N$ is a set containing
$\n$ and closed under successor.  Then by definition of $\L (\n)$, $\L (\n) \subseteq \N$.
That completes the proof of the lemma.

\begin{lemma} \label{lemma:LcapS}  
Suppose there is a double successor
$\s \k = \s \n$ with $\k \in \Stem$ and $\n \neq \k$. 
Then
 $\L(\n) \cap \Stem = \phi$.
 \end{lemma}
 
\noindent{\em Proof}.   Assume $\s \k = \s \n$ and $\k \in \Stem$ and $\n \neq \k$.
We will prove by induction on $j$ that 
 $$j \in \Stem \imp j \not \in \L(\n).$$
 The formula is stratified, giving $j$ index 0; $\Stem$ and $\L(\n)$ are parameters.
 \smallskip
 
  Base case: We must 
show $0 \not\in\L$.  Let $Z = \L(\n)- \{0\}$.  Then $Z$ is closed
under successor, by Theorem~\ref{theorem:successoromitszero}.
And $Z$ contains $\n$, since $\n \neq 0$ by Lemma~\ref{lemma:nneqzero}.
Therefore $\L(\n) \subseteq Z$.  Therefore $0 \not\in\L$, as desired.
\smallskip

Induction step:  Suppose $\s j \in \Stem$. 
 We have $j \in \Stem$ by Lemma~\ref{lemma:Spred}.
By Lemma~\ref{lemma:Soneone},
 $\s j$ is not a double successor.  Therefore,
if $\s j = \s \n$, $j= \n$.  But $j \neq \n$, since
$\n \in \L(\n)$ by definition of $\L(\n)$, but $j \not\in\L$ by
the induction hypothesis.   Therefore,   
$\s j \neq \s \n$.
\smallskip
  
Define $Z = \L(\n) - \{\s j\}$.  I say $Z$ is closed under 
successor.  Let $x \in Z$; then $x \in \L(\n)$, so $\s x \in \L$.
By induction hypothesis $j \not\in\L(\n)$, but $x \in \L$;
therefore $x \neq j$.
If $\s x = \s j$ then $\s j$ is a double successor, 
contradicting Lemma~\ref{lemma:Soneone}, since $\s j \in \Stem$.  Hence $\s x \in Z$
as claimed. 
\smallskip

Now I say $\n \in Z$.  Since $\n \in \L(\n)$ it suffices to show 
that $\n \neq \s j$.  Suppose to the contrary that $\n = \s j$.
Then $\n \in \Stem$, since $\s j \in \Stem$.
  Since $\s\n$ is a double successor, by 
Lemma~\ref{lemma:Smax}, $\n$ is the maximal element of $\Stem$. 
But $\k$ is the maximal element of $\Stem$, by definition of $\k$.
Therefore $\k = \n$, contradiction.  Hence $\n \neq \s j$.
Hence $\n \in Z$, as claimed.
\smallskip

 Therefore $Z$ satisfies the conditions defining $\L(\n)$.
Therefore $\L(\n) \subseteq Z$. Therefore $\s j \not \in \L(\n)$, as 
desired.  That completes the induction step. 
That completes the proof of the lemma.

\begin{lemma}\label{lemma:LcupS}
Suppose there is a double successor
$\s \k = \s \n$ with $\k \in \Stem$ and $\n \neq \k$. 
Then $\N = \L(\n) \cup \Stem$.
\end{lemma}

\noindent{\em Proof}.  Assume $\s \k = \s \n$ and $\k \in \Stem$ and $\n \neq \k$.
We will prove by induction on $x$ that 
\begin{eqnarray}
x \in \N \imp x \in \L(\n) \cup \Stem \label{eq:2027}
\end{eqnarray}
The formula is stratified, giving $x$ index 0; $\L(\n)$ is a parameter.
\smallskip

Base case: $\ChurchZero \in \Stem$, by Lemma~\ref{lemma:SN}.  Therefore
$\ChurchZero \in \L(\n) \cup \Stem$.  That completes the base case.
\smallskip

Induction step: 
Let $x \in \L(\n) \cup \Stem$. Then $x \in \L(\n) \ \lor \ x \in \Stem$. 
\smallskip

Case~1:  
$x \in \L(\n)$. Then by Lemma~\ref{lemma:L1}, $\s x \in \L(\n)$, so $\s x \in \L(\n) \cup \Stem$.
\smallskip

Case~2: $x \in \Stem$.  We have 
\begin{eqnarray*}
\k \in \Stem && \mbox{\qquad by hypothesis}\\
\k \in \N  && \mbox{\qquad by Lemma~\ref{lemma:SN}}\\
x = \k \lor x \neq \k &&\mbox{\qquad by Lemma~\ref{lemma:Sdecidable}}
\end{eqnarray*}
 Therefore we may 
argue by cases according as $x = \k$ or not. 
\smallskip

Case~2a: $x \neq \k$.  Then 
\begin{eqnarray*}
\forall u\,(u \in \N \imp \s u = \s x \imp u = x) &&\mbox{\qquad by Lemma~\ref{lemma:kisunique}}\\
 \s x \in \Stem && \mbox{\qquad by Lemma~\ref{lemma:S1}}\\
x \in \L(\n) \cup \Stem &&\mbox{\qquad by definition of union}
\end{eqnarray*}

Case~2b: $x = \k$.
Then 
\begin{eqnarray*}
\n \in \N && \mbox{\qquad by Lemma~\ref{lemma:LN}}\\
 \s \n \in \L(\n) && \mbox{\qquad by Lemma\ref{lemma:L1}}\\
\s x = \s \k = \s n && \mbox{\qquad by hypothesis}\\
\s x \in \Stem \cup \L(\n) &&\mbox{\qquad by definition of union}
\end{eqnarray*}
That completes the induction step.
 That completes the proof of (\ref{eq:2027}).
\smallskip

By (\ref{eq:2027}), $\N \subseteq \L(\n) \cup \Stem$.  It remains
to prove $\L(\n) \cup \Stem \subseteq \N$.  Suppose $x \in \L(\n) \cup \Stem$.
Then $x \in \L(n)$ or $x \in \Stem$.  If $x \in \L(n)$, then $x \in \N$ by Lemma~\ref{lemma:LN}.
If $x \in \Stem$, then $x \in \N$ by Lemma~\ref{lemma:SN}.  That completes
the proof of the lemma. 

\begin{lemma} \label{lemma:Lindependent}
Suppose there is a double successor
$\s \k = \s \n$ with $\k \in \Stem$ and $\n \neq \k$.   
Then  $\L(\n) = \N- \Stem$.   Consequently $\L(n)$ does not 
depend on the choice of $\n$.
\end{lemma}

\noindent{\em Remark}.  This lemma is never used; we include it only to 
clarify why we keep writing in English ``the loop'', while in formulas
we keep writing $\L(\n)$ as if ``the loop'' depended on $\n$.
\medskip
 
\noindent{\em Proof}.  By Lemma~\ref{lemma:LcupS}, we have $\N = \L (\n) \cup \Stem$.
Therefore it suffices to prove
\begin{eqnarray}
\L(\n) = (\L (n) \cup \Stem) - \Stem  \label{eq:2163}
\end{eqnarray}
Left to right:  Suppose $t \in \L(\n)$.  By Lemma~\ref{lemma:LcapS}, $t \not\in\Stem$.
Therefore $t \in (\L (n) \cup \Stem) - \Stem$, as desired.
\smallskip

Right to left: Suppose $t \in (\L (n) \cup \Stem) - \Stem$.  Then $t \in \L(\n)$. 
That completes the proof of the lemma. 

\begin{lemma} \label{lemma:nissuccessor}
Suppose there is a double successor
$\s \k = \s \n$ with $\k \in \Stem$ and $\n \neq \k$. 
Then $\exists p \in \N\, ( p \in \L(\n) \ \land \ \s p = \n)$.
\end{lemma}

\noindent{\em Proof}.   
We have
\begin{eqnarray*}
\n \not\in \Stem &&\mbox{\qquad by Lemma~\ref{lemma:kisunique}}\\
\ChurchZero \in \Stem  &&\mbox{\qquad by Lemma~\ref{lemma:S1}} \\
 \n \neq \ChurchZero &&\mbox{\qquad by the preceding two lines}\\
\n = \s p  &&\mbox{\qquad for  some $p$, by Lemma~\ref{lemma:predecessor}}
\end{eqnarray*}
 By Lemma~\ref{lemma:LcupS}, $p \in \L(\n) \ \lor \ p \in \Stem$.  We argue
 by cases accordingly.
 \smallskip
 
 Case 1,  $p \in \L(\n)$.  Then we use $p$ to instantiatiate $\exists p$.
 That completes Case~1.
 \smallskip
 
 Case~2,  $p \in \Stem$.
 Since we have decidable equality on $\Stem$.
By Lemma~\ref{lemma:Sdecidable},
we have $p = \k \ \lor \ p \neq \k$.
We argue by cases accordingly.
\smallskip

Case~2a, $p = \k$.
Then $\n = \s p = \s \k = \s \n$,
so $\n = \s \n$.   We use 
$n$ to instantiate $\exists p$.  We have
$\n \in \L(\n)$ by Lemma~\ref{lemma:L1}. 
That completes Case~2a.
\smallskip

Case~2b, $p \neq \k$.  Then 
\begin{eqnarray*}
p < \k && \mbox{\qquad by Lemma~\ref{lemma:Smax}} \\
 \s p \in \Stem &&\mbox{\qquad by Lemma~\ref{lemma:order2}} \\
\Stem \cap \L(\n) = \empty &&\mbox{\qquad by Lemma~\ref{lemma:LcapS}}\\
 \n \in \L(\n)  &&\mbox{\qquad by Lemma~\ref{lemma:L1}}\\
 \s p \in \L(\n) && \mbox{\qquad since $\s p = \n$}\\
 \s p \not\in \Stem && \mbox{\qquad since $\Stem \cap \L(\n) = \phi$}
\end{eqnarray*}
But that contradicts $\s p \in \Stem$.  
That completes Case~2b.
That completes the proof of the lemma.

\begin{theorem} \label{theorem:looponto} 
Suppose there is a double successor
$\s \k = \s \n$ with $\k \in \Stem$ and $\n \neq \k$. 
Then $\s: \L(\n) \to \L(n)$ is onto.
\end{theorem}

\noindent{\em Proof}.  By Lemma~\ref{lemma:L1}, $\L(\n)$ is closed under
successor, so $\s: \L(\n) \to \L(n)$.  Define
$$ Z:= \{ x \in \L(\n) : \exists y\, (y \in \L(\n) \ \land \ \s y = x) \}.$$
The formula is stratified, giving $x$ and $y$ index 0; $\L(\n)$ is a parameter. 
We have
\begin{eqnarray*}
\n \in \L(\n) && \mbox{\qquad by Lemma~\ref{lemma:L1}}\\
\exists p\in \L(\n)\,(\s p = \n) && \mbox{\qquad by Lemma~\ref{lemma:nissuccessor}}\\
 \n \in Z &&\mbox{\qquad by the definition of $Z$}
\end{eqnarray*}
I say that $Z$ is closed under successor.  Suppose $x \in Z$.  
Then
\begin{eqnarray*}
x \in \L(\n) &&\mbox{\qquad by definition of $Z$}\\
\s x \in \L(\n) && \mbox{\qquad by Lemma~\ref{lemma:L1}}\\
\s x \in \Z && \mbox {\qquad  by the definition of $\Z$}
\end{eqnarray*}
We have shown that $Z$ contains $\n$ and is closed under successor. Then by 
the definition of $\L(\n)$, we have $\L(n) \subseteq Z$.  That completes
the proof of the theorem.

\begin{lemma} \label{lemma:kinN}
Suppose there is a double successor
$\s \k = \s \n$ with $\k \in \Stem$ and $\n \neq \k$.
Then $\k \in \N$.
\end{lemma}

\noindent{\em Proof}.  We have $\Stem \subseteq \N$, since $\N = \Stem \ \cup \ \L(n)$
by Lemma~\ref{lemma:LcupS}. Since $\k \in \Stem$, we have $\k \in \N$.
That completes the proof of the lemma.

 \section{The Annihilation Theorem}

\begin{theorem} [Annihilation Theorem]
\label{theorem:annihilation} 
 Suppose $u$ and $m$ are Church numbers such that
$u \oplus  m = u$.  
Let $X$ be any set and let $f:X\to X$ be an injection. 
 Then $f$, iterated $m$ times, is the 
identity on $X$.  In symbols, $m fx = x$ for all $x \in X$.
\end{theorem} 

\noindent{\em Remarks}.  This theorem is proved for any set $X$,
not just for any finite set,  and we do not assume $\N$ is 
finite.  Definition~\ref{definition:ChurchFregenjection} defines ``injection''. 
We also do not need to know that there is only one 
double successor (and even if there is none, the theorem 
is still true, although then $m = 0$ is the only possibility.)
\medskip

\noindent{\em Proof}.   Let $f:X \to \N$, and 
assume $f \in \FUNC$ and $Rel(f)$.  By Lemma~\ref{lemma:iteration},
each iterate $uf$ of $f$ maps $X$ to $X$ and is also a functional relation.
Moreover, by Lemma~\ref{lemma:oneoneiteration}, each iterate of $f$ is also 
one-to-one from $X$ to $X$.
Then
 
\begin{eqnarray*}
u &=& u \oplus m \mbox{\qquad\qquad by hypothesis} \\
ufx &=& (u \oplus  m)fx  \mbox{\qquad since $u$, $u\oplus m$, $uf$, and $(u\oplus m)f$ are functions} \\
&=& uf(mfx) \mbox{\qquad \ \ by Lemma~\ref{lemma:doubleiteration}}\\
\end{eqnarray*}
Since $uf$ is one-to-one from $X$ to $X$,  this implies $x = mfx$.
 That completes the proof of the theorem.
\smallskip

\noindent{\em Remark}.
Nothing proved up to now rules out the possibility that $\n = \s \k$
and $\s \n = \n$.  That would make the loop $\L$ contain only one element,
and $\m$ would be 1.  The following corollary shows that $\m$ is 
much greater.

\begin{corollary} \label{lemma:mbig}  
Suppose  $x = x\oplus  m$ for some 
$x,m \in \N$.  
Then  $m$ is not equal to 1, 2, 3, $\ldots$, where by $1$ we mean $\s \ChurchZero$, etc.  
\end{corollary}

\noindent{\em Remark}.  Formally, this is a different theorem for each 
value of $m$.  We formalized the cases $m=1$ and $m=2$ in Lean, which was
sufficient for our application.
\smallskip

\noindent{\em Proof}. If there is 
a finite set $X$ with a permutation $f$ that is not the identity on $X$, but
  $d f$ is the identity on $X$, then $\m \neq d$.  (Here $d$ is not
a variable, but a specific named integer, with a different proof for each $d$.
) 
For example, when $d = \s \ChurchZero$, we have $\ChurchZero \neq \s \ChurchZero$
by Theorem~\ref{theorem:successoromitszero}. We can define
 a permutation $f$ of $\{ \ChurchZero, \s \ChurchZero\}$
 that interchanges  $\ChurchZero$ and $ \s \ChurchZero$. 
  (It takes about 400 steps to verify 
 that formally, as the definition of permutation has several clauses.) 
 Therefore $\m \neq \s \ChurchZero$.
\smallskip

Therefore there are at least three   elements in $\N$, namely 
$\ChurchZero$, $\n$, and $\s \n$.  I say these are distinct elements.
We have $\s \n \neq \n$ since we have just shown $\m \neq \ChurchZero$.
We have $\s n \neq \ChurchZero$ by Theorem~\ref{theorem:successoromitszero}.
And we have $\n \neq \ChurchZero$ by Lemma~\ref{lemma:nneqzero}.
We can then construct a permutation $f$ of $\N$ such that $jf$ is not
the identity for $j = \s \ChurchZero$ or $j = \s(\s \ChurchZero)$.
Therefore $\m \neq \s (\s \ChurchZero)$.  Then one can show that there
are three distinct members of $\L(\n)$, so there are four distinct 
members of $\n$, and we can construct a permutation of those members
to show that $\m \neq \s(\s(\s\ChurchZero))$.  Similarly we can continue
through any particular value of $\m$.   That is, 
 $\m$ is not equal to any integer with a name, as for such $d$ we can construct
the required permutation.  
\smallskip

\begin{corollary} \label{lemma:snneqn} 
Nothing
 is its
own successor.  That is, for $x \in \N$, we have $\s x \neq x$.
\end{corollary}

\noindent{\em Remark}.  This corollary shows that the loop
does not degenerate to a singleton, in that $\s \n \neq \n$,
but it applies more generally to any Church number $x$.   Thanks
to Albert Visser for pointing out that we can obtain this corollary
immediately 
for any $x$, not just for $\n$.  In fact we do not even need to 
assume that there is a double successor.
\medskip

\noindent{\em Proof}. We have $\s x = x\oplus \s 0 = x \oplus 1$.
Then if $\s x = x$ we have $x = x\oplus m$ with $m=1$, contradicting
Corollary~\ref{lemma:mbig}.  That completes the proof. 

\begin{corollary} \label{lemma:ssnneqn} 
For $x \in \N$, we have $\s(\s x) \neq x$. 
\end{corollary}

\noindent{\em Remark}. To formalize this result, we have to formalize
Lemma~\ref{lemma:mbig} for $\m = 2$, or more precisely, $\m = \s(\s \ChurchZero))$,
which involves constructing a permutation of three elements.  We define
$X = \{ a, b,c\}$ where $a = \ChurchZero$, $b = \s a$, and $c = \s b$.
Those three elements are distinct, by Lemma~\ref{lemma:snneqn} and 
Theorem~\ref{theorem:successoromitszero}. To prove that there is a 
permutation of $X$ requires about 700 steps,  which we omit here.
(There are several 
arguments by cases with nine cases.) Somewhat
surprisingly, one does not need to first prove $X$ is finite.  
\smallskip 

\noindent{\em Proof}.
Suppose $\s(\s x) = x$.  Then
\begin{eqnarray*}
x = x \oplus \ChurchZero && \mbox{\qquad by Lemma~\ref{lemma:ChurchZero_equation}}\\
\s x = \s(x \oplus \ChurchZero) &&\mbox{\qquad by the previous line}\\
 = x \oplus \s \ChurchZero && \mbox{\qquad by Lemma~\ref{lemma:ChurchSuccessorShift}} \\
\s (\s x) = \s (x \oplus \s \ChurchZero)  &&\mbox{\qquad by the previous line}\\
 = x \oplus \s (\s \ChurchZero) &&\mbox{\qquad by Lemma~\ref{lemma:ChurchSuccessorShift}} \\
x = x \oplus \s (\s \ChurchZero) && \mbox{\qquad since $\s (\s x) = x$}
\end{eqnarray*}
But that contradicts Lemma~\ref{lemma:mbig}.
That completes the proof.

\begin{corollary}\label{lemma:klessthann}
 If  $\s \k = \s \n$  with $\k \in \Stem$ and $\n \in \N$ and $\k \neq \n$, then $\k < \n$.
That is, there exists $m \in \N$ such that $\n = \k \oplus  m$.
\end{corollary}

\noindent{\em Proof}.   Define
$$Z = \{ x\in \Stem: x \le \n \}.$$
I say that $Z$ contains 0 
and is closed under successor except double successors.
We have
\begin{eqnarray*}
\ChurchZero \in \Stem  && \mbox{\qquad by Lemma~\ref{lemma:S1}}\\
\ChurchZero \oplus  \n = \n && \mbox{\qquad by Lemma~\ref{lemma:zeroplusx}}\\
\ChurchZero \le \n  && \mbox{\qquad by definition of $<$}\\
\ChurchZero \in Z && \mbox{\qquad by definition of $Z$}
\end{eqnarray*}
Now suppose $x \in Z$ and $\s x$ is not a double successor.  We 
must show $\s x \in Z$. We have
\begin{eqnarray*}
x \in \Stem \ \land \ x \le \n && \mbox{\qquad by the definition of $Z$}\\
\s x \in \Stem &&\mbox{\qquad by Lemma~\ref{lemma:S1}}\\
x \neq \n  && \mbox{\qquad by Lemma~\ref{lemma:kisunique}, since $\s x$ is not a double successor} \\
x < \n && \mbox{\qquad by Lemma~\ref{lemma:Churchletolessthan}}\\
 \s x \le \n &&\mbox{\qquad by Lemma~\ref{lemma:Churchle}} 
\end{eqnarray*}
 Hence $\Stem \subseteq Z$.
Then $\k \in Z$.  Therefore $\k \le \n$.  Since $\k \neq \n$ we 
have $\k < \n$, by Lemma~\ref{lemma:Churchletolessthan}.
  By definition of $<$, there exists $m \in \N$
such that $\n = \k \oplus  m$.  That completes the proof of the corollary.
\medskip

\noindent{\em Remarks}. $m$ is not asserted to be unique.  We do not 
know if $m$ has to be in the loop or has to be in the stem.
\smallskip

In order to apply the Annihilation Theorem (Theorem~\ref{theorem:annihilation}), we need to know that the 
iterates of $f$ still map $X$ to $X$.  That is the content of the next lemma.
\begin{lemma} \label{lemma:xfmaps}
Let $X$ be any set. Let $f:X \to X$, and suppose $f \in \FUNC$ and $Rel(f)$.
   Then 
$$ q \in \N \imp  x \in X \imp  q f x \in X.$$
\end{lemma}

\noindent{\em Proof}. The formula is stratified, giving $x$ index 0,
$f$ index 3, and $q$ index 6.  Therefore we may proceed by induction on $q$.
\smallskip

Base case, $q=\ChurchZero$.  Then
\begin{eqnarray*}
x \in X && \mbox{\qquad by hypothesis} \\
\ChurchZero f x = x && \mbox{\qquad by Lemma~\ref{lemma:zeroAp}}
\end{eqnarray*}
That completes the base case.
\smallskip

Induction step.  
\begin{eqnarray*}
q f x \in X && \mbox{\qquad by the induction hypothesis}\\
\s (q f x) \in X && \mbox{\qquad since $f:X \to X$} \\
\s q f x = \s (q f x) && \mbox{\qquad by Theorem~\ref{theorem:successorequation}}\\
\s q f x \in X && \mbox{\qquad by the preceding two lines}
\end{eqnarray*}
That completes the induction step.
That completes the proof of the lemma.
  
\section{Some consequences of assuming $\N$ is finite}
We take this opportunity to point out that ``$\N$ is not finite'' is,
on the face of it at least, 
a weaker assertion than ``$\N$ is infinite'',  where the latter is 
taken in Dedekind's sense, that the Church successor function is 
one-to-one.   Thus ``$\N$ is finite'' is a stronger assumption 
than ``$\N$ is not infinite''.   In this section we show that under the assumption 
that $\N$ is finite,  we rather quickly reach several important results:
  $\N$ has decidable equality, 
  successor is one-to-one on the loop $\L$, and there is a unique double successor.   
  
 That $\N$ has decidable equality 
is immediate if we assume $\N$ is finite, since according to Lemma~\sref{lemma:finitedecidable}, every finite set has decidable
equality.   
 
\begin{lemma} \label{lemma:loopfinite} 
If $\N$ is finite, and there is a double successor 
$\s \k = \s \n$ with $\k \in \Stem$, then  $\L(\n)$ is finite.
\end{lemma}

\noindent{\em Proof.} Assume $\N$ is finite 
and there is a double successor 
$\s \k = \s \n$ with $\k \in \Stem$.
  I say that 
$\L(\n)$ is a separable subset of $\N$. By Definition~\sref{definition:SSC},
that means that $\N = \L \cup (\N-\L)$.  
By Lemma~\ref{lemma:LcupS}, $\N-\L = \Stem$,
and $\N = \L \cup \Stem$, so $\L$ is a separable subset of $\N$, as claimed.
Then by Lemma~\sref{lemma:separablefinite}, $\L$ is finite.
That completes the proof of the lemma.

\begin{theorem} \label{theorem:looponeone}
If $\N$ is finite, and there is a double successor 
$\s \k = \s \n$ with $\k \in \Stem$, then  
Church successor restricted to $\L(\n)$ is one-to-one.
\end{theorem}

\noindent{\em Proof}.  Assume $\N$ is finite 
and there is a double successor 
$\s \k = \s \n$ with $\k \in \Stem$.
By Lemma~\ref{lemma:loopfinite}, $\L(n) \in \FINITE$.
By Theorem~\ref{theorem:looponto}, successor is onto as a map from $\L(\n)$ to $\L(\n)$.
By Theorem~\sref{theorem:dedekind2},   successor is 
 one-to-one as a map from $\L$ to $\L$. 
That completes
the proof of the theorem.

\begin{theorem} \label{theorem:rho}
Suppose  $\N$ is finite and there is a double successor 
$\s \k = \s \n$ with $\k \neq \n$ and $\n \in \N$ and  $\k \in \Stem$.
 Then there is exactly 
one double successor.  More precisely,  if $j, \ell \in \N$ and $j \neq \ell$
and $j < \ell$ and 
 $\s j = \s \ell$,   
 then $j = \k$ and $\ell = \n$. 
\end{theorem}

\noindent{\em Proof}.  Suppose $\N$ is finite and $\s \k = \s \n$
with $\k \neq \n$ and $\k \in \Stem$.
Suppose $\s j = \s \ell$ with $j < \ell$.   We have to prove $\ell = \n$.
\smallskip

 By Theorem~\ref{theorem:looponeone},
successor is one-to-one on $\L$, so not both $j$ and $\ell$ can belong to $\L(\n)$.
By Lemma~\ref{lemma:LcupS}, each of them belongs to $\L(\n)$ or to $\Stem$,
and by Lemma~\ref{lemma:LcapS}, $\L(n)$ and $\Stem$ are disjoint.
I say that 
\begin{eqnarray}
\ell \in \Stem \imp j \in \Stem \ \land \ \s j \in \Stem \label{eq:2631}
\end{eqnarray}
To prove that, assume $\ell \in \Stem$.  Then 
\begin{eqnarray*}
\ell \in \Stem && \mbox{\qquad by assumption} \\
j \in \Stem && \mbox{\qquad by Lemma~\ref{lemma:Sinit}, since $j < \ell$}\\
\end{eqnarray*}
By Lemma~\ref{lemma:S1}, to prove $\s j \in \Stem$ it suffices to prove 
that $\s j$ is not a double successor. To that end, assume $\s j = \s v$;
we must prove $j = v$.  I say that $j \neq \k$.  Here is the proof:
\begin{eqnarray*}
j = \k && \mbox{\qquad assumption}\\
\ell < \k \ \lor \ \ell = \k &&\mbox{\qquad by Lemma~\ref{lemma:Smax}, since $\ell \in \Stem$}\\
\ell < \k && \mbox{\qquad since $\k \neq \ell$}\\
\k < \k  && \mbox{\qquad by transitivity, since $\k < \ell$} \\
\k \not < \k  && \mbox{\qquad by Lemma~\ref{lemma:knotlessthank}}
\end{eqnarray*}
That contradiction completes the proof that $j \neq \k$. 
Then by Lemma~\ref{lemma:kisunique}, we have $v=j$ as desired.
That completes the proof of (\ref{eq:2631}).
\smallskip

Now I say that $\ell \not \in \Stem$.  To prove that:
\begin{eqnarray*}
\ell \in \Stem && \mbox{\qquad by assumption} \\
j \in \Stem && \mbox{\qquad by (\ref{eq:2631})}\\
\s j \in \Stem && \mbox{\qquad by (\ref{eq:2631})}\\
j = \ell  && \mbox{\qquad by Lemma~\ref{lemma:Soneone}}\\
j \neq \ell && \mbox{\qquad by Lemma~\ref{lemma:Strichotomy}, since $j < \ell$}
\end{eqnarray*}
That contradiction completes the proof that $\ell \not\in \Stem$.
\smallskip

Then
\begin{eqnarray*}
\ell \in \L (\n)  && \mbox{\qquad by Lemma~\ref{lemma:LcupS}, since $\ell \not \in \Stem$}\\
j \in \L(\n) \imp j = \ell && \mbox{\qquad by Theorem~\ref{theorem:looponeone} }\\
j \not\in \L(\n)  && \mbox{\qquad since $j \neq \ell$}\\
j \in \Stem  && \mbox{\qquad by Lemmas~\ref{lemma:LcupS} and \ref{lemma:LcapS}}\\
\Stem = \{ x \in \N : x < \k \lor x = \k\}  && \mbox{\qquad by Lemma~\ref{lemma:Smax}}\\
\s \k \not\in\Stem && \mbox{\qquad by Lemma~\ref{lemma:Soneone} applied to $\k,\n$} \\
 \k \neq \s \k  &&\mbox{\qquad since $\k \in \Stem$  but $\s \k \not\in \Stem$}\\
 \k \in \N && \mbox{\qquad by Lemma~\ref{lemma:kinN}}\\
 j = \k \ \lor \ j \neq \k && \mbox{\qquad by Lemma~\sref{lemma:finitedecidable}, since $\N \in \FINITE$}
\end{eqnarray*}
We argue by cases accordingly.
\smallskip

Case~1, $j=\k$. Then
\begin{eqnarray*}
 \s \k = \s \ell  &&\mbox{\qquad since $\s \k = \s j = \s \ell$}\\
 \s \n = \s \ell  && \mbox{\qquad since $\s \k = \s n$}\\
 \ell \in \L(\n) &&\mbox{\qquad by Lemma~\ref{lemma:LcupS}}\\
\s \ell \in \L(\n) && \mbox{\qquad by Lemma~\ref{lemma:L1}}\\
 \s \n \in \L(\n)  && \mbox{\qquad by Lemma~\ref{lemma:L1}}\\
 \n = \ell && \mbox{\qquad  by Theorem~\ref{theorem:looponeone}}
 \end{eqnarray*}
That completes Case~1.
\smallskip

Case~2, $j \neq \k$. Then 
\begin{eqnarray*}
\s j \mbox{\ is not a double successor} && \mbox{\qquad by Lemma~\ref{lemma:kisunique}}\\ \s j \in \Stem  && \mbox{\qquad by Lemma~\ref{lemma:S1}, since $j \in \Stem$}\\
j = \ell && \mbox{\qquad by Lemma~\ref{lemma:Soneone}} \\
j \neq \ell && \mbox{\qquad by hypothesis}
\end{eqnarray*}
That contradiction completes Case~2.
That completes the proof of the theorem.

\begin{corollary} \label{lemma:rho2}
Suppose  $\N$ is finite and there is a double successor 
$\s \k = \s \n$ with $\k \neq \n$ and $\n \in \N$ and  $\k \in \Stem$.
Suppose $j, \ell \in \N$ and $j \neq \ell$ and 
 $\s j = \s \ell$.  Then $\{j, \ell\} = \{ \k, \n\}$, 
 i.e., $j$ and $\ell$ are $\n$ and $\k$ or $\k$ and $\n$.  
\end{corollary}

\noindent{\em Proof}. By Lemma~\ref{lemma:trichotomy1}, 
$j < \ell \ \lor \ell < j$.  If $j < \ell$ then by 
Theorem~\ref{theorem:rho} we have $j = \k$ and $\ell = \n$.
If $\ell < j$ then (applying Theorem~\ref{theorem:rho} to $\ell$ and $j$ instead of to $j$ and $\ell$\,),
we have  $\ell = \k$ and $j = \n$.  That
completes the proof of the corollary.

\section{A linear order on $\N$}

In this section we introduce a certain linear ordering on $\N$, 
 which we write as $x \preceq y$,  or in its strict version, $x \prec y$.
The definition of $x \preceq y$ will be given in such a way that it 
does not presume that $\N$ is finite or that there is a double successor,
because we need it under those conditions near the end of the paper,
after we have proved $\N$ is not finite but still need to prove $\N$ is 
infinite. 
\smallskip

The intuitive meaning of $x \preceq y$ is that we come to $x$ before $y$ as we 
trace out the stem and then the loop (also allowing $x=y$).

\begin{definition} \label{definition:preceq_aux}
For $X \subseteq \N$, we say ``$X$ is closed under successors except greater double successors''  to mean
$$ \forall u\, (u \in X \imp \forall v \in \N\,(v < u \imp \s u = \s v \imp u = v) \imp \s u \in X)).$$
\end{definition}

\noindent{\em Remark}.  $\n$ does not actually appear in the definition of ``closed under 
successors except greater double successors,''  but the following lemma shows that,
if there is a double successor,  it really means ``closed under successors except $\n$.''
However, the definition does not {\em assume} that there is a double successor. 
\smallskip

\begin{lemma}\label{lemma:preceq_helper}
Suppose $\N$ is finite, and  $\s\k = \s\n$ with $\k \in \Stem$ and $\n \neq \k$.  Let $X \subseteq \N$.
 Then 
$$ \forall u\, (u \in X \imp \forall v \in \N\,( v < u \imp \s u = \s v \imp u = v) \imp \s u \in X))$$
(which is the formula in the preceding definition) is equivalent to  
$$ \forall u\, (u \in X  \imp u \neq \n  \imp \s u \in X).$$
\end{lemma}

\noindent{\em Proof}. Let $\s\k = \s\n$ with $\k \in \Stem$ and $\n \neq \k$,
and let $X$ be any set. 
We have to prove
\begin{eqnarray*}
&& \forall u\, (u \in X \imp \forall v\in \N\,(v < u \imp (\s u = \s v \imp u = v) \imp \s u \in X)) \\
\iff && 
 \forall u\, (u \in X  \imp u \neq \n  \imp \s u \in X).
\end{eqnarray*}
Left to right:  Assume $u \in X$ and $u \neq \n$.  Instantiating the left side to $u$,
we see that it suffices to prove
$$ \forall v \in \N\,(v < u \imp \s u = \s v \imp u = v).$$
Suppose $v \in \N$ and $v < u$ and $\s u = \s v$.  Since $u \in X$ and $X \subseteq \N$,
we have $u \in \N$.  Since $\N$ is finite, it has decidable equality, by Lemma~\sref{lemma:finitedecidable}.
Therefore $ u = v \ \lor \ u \neq v$.  If $u = v$, the desired conclusion is immediate,
so we may assume $u \neq v$. Then
by Theorem~\ref{theorem:rho}, since $u \neq \n$,
we have $u=v$ as desired.  That completes the left-to-right direction.
\smallskip

Right to left: Assume 
\begin{eqnarray}
 \forall u\, (u \in X  \imp u \neq \n  \imp \s u \in X)) \label{eq:E2793}
\end{eqnarray}
and suppose $u\in X$ and 
\begin{eqnarray}
\forall v < u\,(\s u = \s v \imp u = v). \label{eq:E2812}
\end{eqnarray}
  We must prove 
$\s u \in X$.  We have
\begin{eqnarray*}
\k \in \N && \mbox{\qquad by Lemma~\ref{lemma:kinN}}\\
\k < \n  && \mbox{\qquad by Corollary~\ref{lemma:klessthann}}
\end{eqnarray*}
 We have $u \neq \n$, since if $u = \n$ then taking $u = \n$ and $v = \k$ in (\ref{eq:E2812})  we have
$\s u = \s v$, so $u = v$, i.e., $\n = \k$, contradiction.   Then by (\ref{eq:E2793}),
we have $\s u \in X$, as desired.  That completes the right-to-left direction.
That completes the proof of the lemma. 

\begin{definition}\label{definition:preceq}
The relation $x \preceq y$,  
 means that $x\in \N$ and $y \in \N$ and $y$ belongs to 
  every separable
subset of $\N$ containing $x$ and closed under successors except greater
double successors.
\smallskip

Explicitly,
\begin{eqnarray*}
&& x \preceq y  \iff  \forall X\,( \N = X \cup (\N -X) \imp  
        x \in X  \\
 &&\imp
   \forall u\, (u \in X \imp \forall v\in \N\,(v < u \imp (\s u = \s v \imp u = v) \imp \s u \in X) \imp y \in X.
\end{eqnarray*}
\end{definition}

The definition is stratified, giving $x$ and $y$ index 0 and $X$ index 1.  $\N$ is a parameter.
Since $x$ and $y$ get the same index, the relation $x \preceq y$ is definable in INF.

\begin{lemma} \label{lemma:preceq} Suppose $\N \in \FINITE$ and $\s \k = \s \n$ and $  \k \neq \n$ and
$\k \in \Stem$ and $\n \in \N$.  Then  for all $x,y \in \N$, we have $ x \preceq y$ if and only if
$$   \forall w\,(\N = X \cup (\N -X) \imp x \in w \imp 
(\forall u\,(u \in w \imp u \neq \n \imp \s u \in w)     \imp y \in w).$$
\end{lemma}

\noindent{\em Remark}.  Although $\n$ appears in this lemma,
$\n$ does not appear in the definition of $\preceq$.  We can use therefore
use this lemma to express $\preceq$ in terms of any (hypothesized) double successor,
without it depending on the particular double successor. 
\medskip

\noindent{\em Proof}. Using Lemma~\ref{lemma:preceq_helper} (in the right-to-left direction)
to rewrite the closure condition in the lemma, we see that it suffices to prove
\begin{eqnarray*} x \preceq y \iff 
\forall w\,(\N = X \cup (\N -X) \imp x \in w \imp  \\
 \forall u\, (u \in w \imp \forall v \in \N\,( v < u \imp \s u = \s v \imp u = v) \imp \s u \in w) \\
 \imp y \in w)
\end{eqnarray*}
But that is just Definition~\ref{definition:preceq} (up to renaming a bound variable).
  That completes the proof of the lemma.

\begin{definition}\label{definition:prec}
We define $$x \prec y \iff x \preceq y \ \land x\ \neq y.$$
\end{definition}

\begin{lemma} [Transitivity of $\preceq$] \label{lemma:preceqtrans} 
For $x,y,z \in \N$ we have 
$$ x \preceq y \imp y \preceq z \imp x \preceq z.$$
\end{lemma}

\noindent{\em Proof}.  Suppose $x \preceq y$ and $y \preceq z$.
Let $X$ be a separable subset of $\N$  closed under successors except 
greater double successors.    Suppose $x \in X$.  Since $x \preceq y$
we have $y \in X$.  Since $y \preceq z$, we have $z \in X$.  Then 
by the definition of $\preceq$,  we have $x \preceq z$.  
That completes the proof of the lemma.

\begin{lemma} \label{lemma:preceqreflexive}
For $x \in \N$ we have $x \preceq x$.
\end{lemma}

\noindent{\em Proof}. $x$ belongs to every separable set $X$ containing $x$ and 
satisfying some condition; putting in the particular condition from the 
definition of $\preceq$ we have the desired result.  That completes the proof of the lemma.

\begin{lemma} \label{lemma:preceqzero} 
 Suppose $x \in \N$ and $x \preceq \ChurchZero$.  Then $x = \ChurchZero$.
\end{lemma}

\noindent{\em Proof}.   By Theorem~\ref{theorem:successoromitszero}, $Z:= \N - \{\ChurchZero\}$ is 
  closed under successor.   
By Lemma~\ref{lemma:decidable0},  $x = \ChurchZero \ \lor \ x \neq \ChurchZero$.
Therefore $Z$ is a separable subset of $\N$. 
If $x \neq \ChurchZero$, then $x \in Z$. 
Since $x \preceq \ChurchZero$,   then 
$\ChurchZero \in Z$.  But $\ChurchZero \not\in Z$.  Therefore 
$x = \ChurchZero$.   That completes the proof of the lemma.

\begin{lemma} \label{lemma:precmin2}  For $x \in \N$,  we have $\neg\, (x \prec \ChurchZero)$.
\end{lemma}

\noindent{\em Proof}.  Suppose $x \in \N$ and $x \prec \ChurchZero$.  By definition of $\prec$, $x \preceq \ChurchZero$
and $x \neq \ChurchZero$. By Lemma~\ref{lemma:preceqzero}, $x = \ChurchZero$, contradiction.
That completes the proof of the lemma.

\begin{lemma} \label{lemma:preceqsuccessor} 
 Suppose $\N \in \FINITE$ and $\s \k = \s \n$ and $  \k \neq \n$ and
$\k \in \Stem$ and $\n \in \N$.  Then for $x,y \in \N$ we have
$$ y \neq \n \imp  x \preceq \s y \iff  x \preceq y \ \lor \ x = \s y.$$
\end{lemma}

\noindent{\em Remark}.  The reader should refer to Fig.~\ref{figure:rho} 
to see why the condition $y \neq \n$ is needed.
\smallskip

\noindent{\em Proof}.  Suppose $y \neq \n$.  Left to right: suppose $x \preceq \s y$.
Since $\N$ is finite, it has decidable equality, so $x = \s y \ \lor \ x \neq \s y$.
If $x = \s y$ we are done, so we may suppose $x \neq \s y$. We
have to prove $x \preceq y$.  Let $X$ be a separable subset
of $\N$ closed under successor except $\n$ and containing $x$.  We have to prove 
$y \in X$.  Since $X$ is a separable subset of $\N$, we have $y \in X \ \lor \ y \not\in X$.
If $y \in X$, we are done, so we may assume $y \not\in X$. 
 Define $Z:= X- \{ \s y\}$. Since $X$ is a separable subset of $\N$ and 
$\N$ has decidable equality, $Z$ is a separable subset of $\N$ (80 steps omitted).
Since $x \neq \s y$, we have $x \in Z$. 
\smallskip

I say that $Z$ is closed under successor except $\n$.
To prove that, suppose $u \in \Z$ and $u \neq \n$; we must prove $\s u \in \Z$. We have
$\s u \in X$ since $X$ is closed under successor except $\n$.  Since $Z = X - \{ \s y\}$,
it suffices to prove $\s u \neq \s y$.  Suppose that $\s u = \s y$; we must 
derive a contradiction.   
\begin{eqnarray*}
u = y && \mbox{\qquad by Corollary ~\ref{lemma:rho2}, since $u \neq \n$ and $y \neq \n$}\\
u \not\in X && \mbox{\qquad since $y \not\in X$ and $u = y$}\\
u \in X   && \mbox{\qquad since $u \in Z = X - \{\s y\}$}
\end{eqnarray*}
  That completes the proof that 
$Z$ is closed under successor except $\n$.  Since $x \in Z$ and $x \preceq \s y$, 
we have $\s y \in Z$.  But that is a contradiction.  That completes the proof of the 
left-to-right direction of the lemma.
\smallskip

Right to left: Suppose $y \neq \n$ and  $x \preceq y \ \lor \ x = \s y$.
We must prove $x \preceq \s y$.   
Let $X$ be a separable subset of $\N$ closed
under successor except $\n$ and containing $x$; we must prove $\s y \in X$. 
Since $x \preceq y \ \lor \ x = \s y$, we may argue by cases.
\medskip

Case~1, $x \preceq y$.  Since $x \in X$ and $y \neq \n$, we have $\s y \in X$. 
That completes case~1.
\smallskip

Case~2, $x = \s y$.  Since $x \in X$ we have $\s y \in X$.
That completes Case~2.       
That completes the proof of the lemma.

\begin{corollary} \label{lemma:xpreceqsx} 
Suppose $\N \in \FINITE$ and $\s \k = \s \n$ and $  \k \neq \n$ and
$\k \in \Stem$ and $\n \in \N$.  Then for $x,y \in \N$ we have 
$$x \neq \n \imp x \preceq \s x.$$
\end{corollary}

\noindent{\em Proof}. By Lemma~\ref{lemma:preceqreflexive}, we have $x \preceq x$.
Taking $y = x$ in Lemma~\ref{lemma:preceqsuccessor}, we have $x \preceq \s x$ 
as desired.  That completes the proof of the corollary.

\begin{corollary}\label{lemma:xprecsx}
  Suppose $\N \in \FINITE$ and $\s \k = \s \n$ and $  \k \neq \n$ and
$\k \in \Stem$.  Then for all $x \in \N$, we have 
$$ x \neq \n \imp x \prec \s x.$$
\end{corollary}

\noindent{\em Proof}.  Suppose $x \neq \n$ and $x \in \N$.  We have
\begin{eqnarray*}
x \preceq \s x && \mbox{\qquad by Lemma~\ref{lemma:xpreceqsx}}\\
x \neq \s x && \mbox{\qquad by Lemma~\ref{lemma:snneqn}}\\
x \prec \s x && \mbox{\qquad by definition of $\prec$}
\end{eqnarray*}
That completes the proof of the lemma.

\begin{lemma} \label{lemma:leftpreceqsuccessor}
Suppose $\N \in \FINITE$ and $\s \k = \s \n$ and $  \k \neq \n$ and
$\k \in \Stem$ and $\n \in \N$.  Then for $x,y \in \N$ we have
$$   x \prec   y \imp  \s x \preceq y.$$
\end{lemma}

\noindent{\em Proof}.  Suppose $x \prec y$.  By definition of $\prec$, 
$x \preceq y$ and $x \neq y$.  We must prove $\s x \preceq y$. 
By Lemma~\ref{lemma:preceq}, it suffices to show that for every separable
subset $X$ of $\N$ that contains $\s x$ and is closed under successor except $\n$,
we have $y \in X$.   Let $X$ be such a set, and define
\begin{eqnarray}
 Z := X \cup \{ x\} \label{eq:E2991}
\end{eqnarray}
Since $\N$ is finite, it has decidable equality; hence $Z$ is a separable subset of $\N$.
Since $\s x \in X$, we have $\s x \in Z$, so $Z$ is closed under successor except $\n$.
We have $x \in Z$ by (\ref{eq:E2991}). 
Since $x \preceq y$, we have $y \in Z$ by Lemma~\ref{lemma:preceq}. 
Since $x \neq y$, we have $y \in X$ by (\ref{eq:E2991}). 
That completes the proof of the lemma.
 
\begin{lemma}\label{lemma:precmax}
Suppose $\N \in \FINITE$ and $\s \k = \s \n$ and $  \k \neq \n$ and
$\k \in \Stem$ and $\n \in \N$.  Then for all $x \in \N$ we have $x \preceq \n$.
\end{lemma}

\noindent{\em Proof}. 
Suppose $\N \in \FINITE$ and $\s \k = \s \n$ and $  \k \neq \n$ and
$\k \in \Stem$ and $\n \in \N$.  We begin by proving
\begin{eqnarray}
\s \k \preceq \n \label{eq:E3033}
\end{eqnarray}
To prove that, let $X$ be a separable subset of $\N$ containing
$\s \k$ and closed under successor except $\n$. We must prove $\n \in X$.
\smallskip

 Define $Z:= X \cup \{\n\}$.
I say that $Z$ is closed under successor.
To prove that: if $x \in Z$ then $x \in X \ \lor \ x = \n$.
If $x \in X$ and $x \neq \n$, then $\s x \in X$, since $X$ is closed under successor except $\n$.
But if $x = \n$ then $\s x = \s \n = \s k$, which is in $X$ by hypothesis, and hence in $Z$.
  Since $\N$
is finite, it has decidable equality, so these cases are exhaustive.  Hence $Z$ is 
closed under successor, as claimed.  
\smallskip

Now we can prove $\n \in X$:
\begin{eqnarray*}
\L(n) \subseteq Z && \mbox{\qquad by the definition of $\L(\n)$}\\
\n \in \L(\n)     && \mbox{\qquad by Lemma~\ref{lemma:L1}}\\
\ChurchZero \in \Stem  &&\mbox{\qquad by Lemma~\ref{lemma:S1}}\\
\L(n) \cap \Stem = \empty && \mbox{\qquad by Lemma~\ref{lemma:LcapS}}\\
\n \neq \ChurchZero && \mbox{\qquad by the preceding lines}\\
\n = \s r   && \mbox{\qquad for some $r \in \L(\n)$, by Theorem~\ref{theorem:looponto}}\\
r \in Z   && \mbox{\qquad since $\L(\n) \subseteq Z$}\\
r \neq n  && \mbox{\qquad by Lemma~\ref{lemma:snneqn}}\\
r \in X   && \mbox{\qquad since $Z = X \cup \{n\}$}\\
\s r \in X &&\mbox{\qquad since $X$ is closed under successor except $\n$, and $r \neq \n$}\\
\n \in X   && \mbox{\qquad since $\s r = \n$}
\end{eqnarray*}
That completes the proof that $\n \in X$.  That
completes the proof of (\ref{eq:E3033})
\smallskip

We must prove
$$ \forall x \in \N\,( x\preceq \n).$$
We will prove this by induction on $x$.
The formula to be proved is stratified, giving 
$x,\k,\n$ all index 0, since $\preceq$ is a definable relation, appearing
here as a parameter.     
Therefore we may proceed by induction on $x$.
\smallskip

Base case. $\ChurchZero \preceq \n$ by Lemma~\ref{lemma:preceqzero}.
\smallskip

Induction step.  The induction hypothesis is $x \preceq \n$.  
We have to prove $\s x \preceq n$.  Since $\N$ is finite,
it has decidable equality, so we have $x = \n \ \lor \ x \neq \n$.
If $x = \n$, we are done by (\ref{eq:E3033}), since $\s \n = \s \k$.
Therefore we may assume $x \neq \n$.
Then 
\begin{eqnarray*}
x \prec \n && \mbox{\qquad by definition of $\prec$}\\
\s x \preceq \n && \mbox{\qquad by Lemma~\ref{lemma:leftpreceqsuccessor} with $\n$ for $y$}
\end{eqnarray*} 
That completes the induction step.  That completes the proof of the lemma.

\begin{lemma}[Finite induction] \label{lemma:finiteinduction} 
Suppose $\N \in \FINITE$ and $\s \k = \s \n$ and $  \k \neq \n$ and
$\k \in \Stem$ and $\n \in \N$.
Suppose $\ChurchZero \in X$ and $\forall u\,(u \in X \imp u \neq \n \imp \s u \in X)$.
Then $\N \subseteq X$.
\end{lemma}

\noindent{\em Proof}. We will prove by  induction on $z$ that
\begin{eqnarray}
 \forall x \in \N\,(x \preceq z \imp x \in X). \label{eq:E3150}
\end{eqnarray}
The formula is stratified, giving $x$ and $z$ index 0, since $\preceq$ is a 
definable relation, so we may proceed by induction.
\smallskip

Base case, $z = \ChurchZero$. We must show $x \preceq \ChurchZero \imp x \in X$.
Suppose $x \preceq \ChurchZero$.  By Lemma~\ref{lemma:preceqzero}, $x = \ChurchZero$.
Then $x \in X$ by hypothesis.  That completes the base case.
\smallskip

Induction step.  Suppose 
\begin{eqnarray}
\forall x \in \N\,(x \preceq z \imp x \in X) \label{eq:E2940}
\end{eqnarray}
  and 
suppose $x \preceq \s z$.  We must prove $x \in X$. 
\begin{eqnarray*}
\N \in \DECIDABLE && \mbox{\qquad by Lemma~\sref{lemma:finitedecidable}}\\
z = \n \ \lor \ z \neq \n && \mbox{\qquad  since $\N \in \DECIDABLE$}
\end{eqnarray*}
We argue by cases.
\smallskip

Case~1, $z = \n$. Then 
\begin{eqnarray*}
\s z = \s \n = \s \k &&\mbox{\qquad since $z = \n$ and $\s \n = \s \k$}\\
x \preceq \n \imp x \in X && \mbox{\qquad by (\ref{eq:E2940})}\\
x \preceq \n  && \mbox{\qquad by Lemma~\ref{lemma:precmax}} \\
x \in X && \mbox{\qquad by the preceding two lines}
\end{eqnarray*}
That completes Case~1.
\smallskip

Case~2.  $z \neq \n$. 
\begin{eqnarray*}
x \preceq \s z  && \mbox{\qquad by hypothesis}\\
x \preceq z \ \lor \ x = \s z && \mbox{\qquad by Lemma~\ref{lemma:preceqsuccessor}, since $z \neq \n$} \\
\end{eqnarray*}
 If $x \preceq z$, we are done by (\ref{eq:E2940}),
so we can assume $x = \s z$.  Since $X$ is closed under successor except $\n$, 
and $z \neq \n$,  we have $x \in X$.  That completes the induction step.
That completes the proof of (\ref{eq:E3150}). 
\smallskip

Now under the assumptions of the lemma, we have to prove $\N \subseteq X$.  It suffices
to prove that for all $z$, $z \in \N \imp z \in X$.  Assume $z \in \N$.  Substituting
$z$ for the bound variable $x$ in (\ref{eq:E3150}),
we have $z \preceq z \imp z \in X$.  By Lemma~\ref{lemma:preceqreflexive}, we have 
$z \preceq z$.  Hence $z \in X$ as desired.  That completes the proof of the lemma.

\begin{theorem}\label{theorem:prectrichotomy1}
 Suppose $\N \in \FINITE$ and $\s \k = \s \n$ and $  \k \neq \n$ and
$\k \in \Stem$ and $\n \in \N$.  Then for $x,y \in \N$ we have 
$$  x \prec y \ \lor\ x = y \ \lor \ y \prec x.$$
\end{theorem}

\noindent{\em Proof}.  We will prove by finite induction on $x$ that 
\begin{eqnarray}
 \forall y\in \N\, (  x \preceq y \ \lor\  \ y \preceq x) 
 \label{eq:E3108}
\end{eqnarray}
Since $\N$ is finite, it has decidable equality, so $x \preceq y \iff x \prec y \ \lor \ x = y$.
Therefore (\ref{eq:E3108}) is equivalent to the lemma as stated.
\smallskip

The formula is stratified, giving $x$ and $y$ both index 0, since $\preceq$ is
a definable relation.  Hence induction is legal.
\smallskip

Base case.  By Lemma~\ref{lemma:preceqzero}, we have $\ChurchZero \preceq y$.   
 That completes the base case.
\smallskip

Induction step.  The induction hypothesis is (\ref{eq:E3108}).  Let $y$ be given. We must prove
\begin{eqnarray}
\s x \preceq y  \ \lor\  y \preceq \s x     \label{eq:E3124}  
\end{eqnarray}
Since we are using finite induction on $x$ (Lemma~\ref{lemma:finiteinduction}), we may assume 
\begin{eqnarray}
x \neq \n \label{eq:E3181}
\end{eqnarray}
 By (\ref{eq:E3108})
we have $x \preceq y   \ \lor \ y \preceq x$.  We argue by cases accordingly.
\smallskip

Case~1: $x \preceq y$.  By decidable equality and the definition of $\prec$,
we have $x \prec y$ or $x = y$.  If $x \prec y$, then $\s x \preceq y$,
 by Lemma~\ref{lemma:leftpreceqsuccessor}.  If $x = y$ then $y \preceq \s x$,
 by Lemma~\ref{lemma:preceqsuccessor}. That completes Case~1. 
\smallskip

Case~2:  $y \preceq x$.  Then $y \preceq \s x$, by Lemma~\ref{lemma:preceqsuccessor} and (\ref{eq:E3181}).
That completes Case~2.   That completes the proof of the lemma.

\begin{lemma} \label{lemma:successorprec}
   Suppose $\N \in \FINITE$ and $\s \k = \s \n$ and $  \k \neq \n$ and
$\k \in \Stem$ and $\n \in \N$.  Then for $x,y \in \N$ we have 
$$ y \neq \n \imp \s y \prec x \imp y \preceq x.$$
\end{lemma}

\noindent{\em Proof}.  By finite induction on $x$.  
\smallskip

Base case.  We must prove $\s y \prec \ChurchZero \imp y \preceq \ChurchZero$.
But $\s y \prec \ChurchZero$ can never hold, by Lemma~\ref{lemma:precmin2}. 
That completes the base case.
\smallskip

Induction step.  Suppose $\s y \prec \s x$ and $x \neq \n$.  We have to 
prove $y \preceq \s x$.  We have
\begin{eqnarray*}
\s y \preceq \s x  \ \land \ \s y \neq \s x && \mbox{\qquad by definition of $\prec$}\\
y \neq \n      && \mbox{\qquad by hypothesis}\\
\s y \preceq x  && \mbox{\qquad by Lemma~\ref{lemma:preceqsuccessor}, since $y \neq \n$}\\
y \preceq \s y  && \mbox{\qquad by Lemma~\ref{lemma:preceqsuccessor}}\\
y \preceq x    && \mbox{\qquad by Lemma~\ref{lemma:preceqtrans} and the preceding two lines}\\
y \preceq \s x  && \mbox{\qquad by Lemma~\ref{lemma:preceqsuccessor}}
\end{eqnarray*}
That completes the induction step.  That completes the proof of the lemma.

\begin{lemma} \label{lemma:sxnotpreceqx}
   Suppose $\N \in \FINITE$ and $\s \k = \s \n$ and $  \k \neq \n$ and
$\k \in \Stem$ and $\n \in \N$.  Then for $x \in \N$ we have  
$$x \neq \n \imp \neg\,(\s x \preceq x).$$
\end{lemma} 

\noindent{\em Proof}.  Suppose $x \in \N$ and $x \neq \n$ and 
\begin{eqnarray}
\s x \preceq x. \label{eq:E3369}
\end{eqnarray}
We must derive a contradiction. Define
$$ Z:= \{ u \in \N: \neg\, (\s x \preceq u \ \land \ u \preceq x)\}.$$
The formula is stratified, giving $x$ and $u$ index 0, so the definition
can be given in INF.  
\smallskip

I say that $\N \subseteq Z$. By the definition of $\subseteq$, that is equivalent to
\begin{eqnarray}
\forall u\in \N\,(u \in Z) \label{eq:E3398}
\end{eqnarray}
We will prove that by finite induction. 
\smallskip

Base case, $\ChurchZero \in Z$.  We have to prove $\neg\,( \s x \preceq \ChurchZero \ \land \ \ChurchZero \preceq x).$  It suffices to prove $\neg\,(\s x \preceq \ChurchZero).$  Suppose $\s x \preceq \ChurchZero$.
Then by Lemma~\ref{lemma:preceqzero}, $\s x = \ChurchZero$.  But that contradicts
Theorem~\ref{theorem:successoromitszero}.  That completes the base case.
\smallskip

Induction step.  We have to prove that 
$$  u \in Z \imp u \neq \n \imp \s u \in Z.$$
Using the definition of $Z$, that becomes
$$ \neg\, (\s x \preceq u \ \land \ u \preceq x)  \imp u \neq \n \imp \neg\,(\s x \preceq \s u \ \land \ \s u \preceq x).$$
It suffices to prove, assuming $u \neq \n$, that
$$ \s x \preceq \s u \ \land \ \s u \preceq x  \imp \s x \preceq u \ \land \ u \preceq x.$$
Suppose
\begin{eqnarray}
 && u \neq \n   \label{eq:E3373} \\
 && \s x \preceq \s u \ \land \ \s u \preceq x. \label{eq:E3408}
 \end{eqnarray}
We must prove
\begin{eqnarray}
\s x \preceq u \ \land \ u \preceq x \label{eq:E3412}
\end{eqnarray}
We have
\begin{eqnarray*}
\s x \preceq u \ \lor \ \s x = \s u  && \mbox{\qquad by Lemma~\ref{lemma:preceqsuccessor} and (\ref{eq:E3408}) and (\ref{eq:E3373})}
\end{eqnarray*}
We argue by cases accordingly to prove (\ref{eq:E3412}).
\smallskip

Case~1, $\s x \preceq  u$.  That is already the first half of (\ref{eq:E3412}); it 
remains to prove $u \preceq x$.  We have 
\begin{eqnarray*}
\s u \preceq x && \mbox{\qquad by (\ref{eq:E3408})}\\
u \preceq u && \mbox{\qquad by Lemma~\ref{lemma:preceqreflexive}}\\
u \preceq \s u && \mbox{\qquad by Lemma~\ref{lemma:preceqsuccessor}, since $u \neq \n$}\\ 
u \preceq x && \mbox{\qquad by Lemma~\ref{lemma:preceqtrans}}
\end{eqnarray*}
That completes Case~1.
\smallskip

Case~2, $ \s x = \s u$.   Since $\N$ has decidable equality, we have
$x = u \ \lor \ x \neq u$.  We argue by cases.
\smallskip

Case 2a,  $x = u$. Then (\ref{eq:E3412}) becomes 
$\s x \preceq x \ \land \ x \preceq x$, which follows from (\ref{eq:E3369})
and Lemma~\ref{lemma:preceqreflexive}.
\smallskip

Case 2b. $x \neq u$. 
  Then $\s x = \s u$
is a double successor.   Then
\begin{eqnarray*}
\{x,u\} = \{\k,\n\} &&\mbox{\qquad by Corollary~\ref{lemma:rho2}}\\
x = \n \ \land \ u = \k  && \mbox{\qquad since $u \neq \n$ by (\ref{eq:E3373})}\\
 x \neq \n && \mbox{\qquad by hypothesis}
 \end{eqnarray*}
The last two lines are contradictory.  
That contradiction completes the proof of (\ref{eq:E3412}).
That completes the proof that $\N \subseteq Z$.
\smallskip

Therefore $x \in Z$.  But by hypothesis we have $\s x \preceq x$,  and by Lemma~\ref{lemma:preceqreflexive}
we have $x \preceq x$. Hence $x \not\in Z$.  That contradiction completes the proof 
of the lemma.

\begin{lemma}\label{lemma:predecessornotn} 
   Suppose $\N \in \FINITE$ and $\s \k = \s \n$ and $  \k \neq \n$ and
$\k \in \Stem$ and $\n \in \N$.   Let $x \in \N$ with $x \neq \ChurchZero$.  Then there exists an $r \in \N$
with $\s r = x$ and $r \neq \n$.
\end{lemma} 

\noindent{\em Proof}.  Let $x \in \N$ be nonzero. By Lemma~\ref{lemma:predecessor},
there exists $u \in \N$ with $\s u = x$.  Since $\N$ is finite, it has decidable equality,
so $u = \n \ \lor \ u \neq \n$. 
 If $u \neq \n$, we may take $r = u$, and then we are done.
If $u  = \n$, then 
 $x = \s \n = \s \k$.  Since $\k \neq \n$ we may take $r = \k$.  That completes the proof of the lemma.

\begin{theorem} \label{theorem:prectrichotomy3}
  Suppose $\N \in \FINITE$ and $\s \k = \s \n$ and $  \k \neq \n$ and
$\k \in \Stem$ and $\n \in \N$.  Then for $x,y \in \N$ we have 
$$  x \preceq y  \imp  y \preceq x \imp x = y.$$
\end{theorem}

\noindent{\em Proof}. Since $\N$ is finite, it has decidable equality.
Therefore the lemma
as stated is equivalent to
$$ x \neq y \imp \neg (x \preceq y \ \land \ y \preceq x).$$
That formula is stratified, since $\preceq$ and $\prec$ are definable
relations.  We will prove it by finite induction on $y$.
\smallskip

Base case, $x \neq \ChurchZero \imp \neg\,(x \preceq \ChurchZero \ \land \ChurchZero \preceq x).$
is immediate from Lemma~\ref{lemma:preceqzero}.
\smallskip 

Induction step.  Suppose $y \neq \n$ and $x \neq \s y$.  We must prove
\begin{eqnarray}
 \neg\, (x \preceq \s y  \ \land \ \s y \preceq x). \label{eq:E3297}
\end{eqnarray}
To prove that, we must derive a contradiction from 
\begin{eqnarray}
x \preceq \s y  \label{eq:E3273} \\
\s y \preceq x  \label{eq:E3274}
\end{eqnarray}
We have
\begin{eqnarray}
x \preceq y  && \mbox{\qquad by Lemma~\ref{lemma:preceqsuccessor} and (\ref{eq:E3273}) and $y \neq \n$ and $x \neq \s y$} 
\label{eq:E3278}\\
\s y \prec x  && \mbox{\qquad by (\ref{eq:E3274}) and $x \neq \s y$ and the definition of $\prec$} \nonumber\\
x \neq y \label{eq:E3311} &&\mbox{\qquad by Lemma~\ref{lemma:sxnotpreceqx} and (\ref{eq:E3274})  } \nonumber\\
y \preceq x &&\mbox{\qquad by Lemma~\ref{lemma:successorprec}} \label{eq:E3309} \\
\neg\,(x \preceq y \ \land \ y \preceq x) && \mbox{\qquad by (\ref{eq:E3311}) and the induction hypothesis} \nonumber
\end{eqnarray}
But that is contradicted by (\ref{eq:E3278}) together with (\ref{eq:E3309}).  
That completes the induction step.  That completes the proof of the lemma.

\begin{theorem} \label{theorem:prectrichotomy2}
  Suppose $\N \in \FINITE$ and $\s \k = \s \n$ and $  \k \neq \n$ and
$\k \in \Stem$ and $\n \in \N$.  Then for $x,y \in \N$ we have 
$$ \neg (x \prec y \ \land\  y \prec x).$$
\end{theorem}

\noindent{\em Proof}.  Suppose $x \prec y$ and $y \prec x$. 
Then we have
\begin{eqnarray*}
x \preceq y  && \mbox{\qquad by definition of $\prec$}\\
y \preceq x  && \mbox{\qquad by definition of $\prec$}\\
x = y  && \mbox{ \qquad by Theorem~\ref{theorem:prectrichotomy3}}\\
x \prec x && \mbox{\qquad since $x \prec y$ and $x = y$}\\
x \neq x  && \mbox{\qquad by definition of $\prec$}
\end{eqnarray*}
That contradiction completes the proof of the theorem.

\begin{lemma}\label{lemma:precmax2} 
 Suppose $\N \in \FINITE$ and $\s \k = \s \n$ and $  \k \neq \n$ and
$\k \in \Stem$ and $\n \in \N$.  Let $x \in \N$. 
Then $\neg\,(\n \prec x).$
\end{lemma}

\noindent{\em Proof}. Suppose $x \in \N$.  Then 
\begin{eqnarray*}
x \preceq n && \mbox{\qquad by Lemma~\ref{lemma:precmax}}\\
 \n \prec x && \mbox{\qquad assumption, for proof by contradiction}\\
\n \preceq x \ \land \ n \neq x && \mbox{\qquad by definition of $\prec$}\\
x = \n   && \mbox{\qquad by Theorem~\ref{theorem:prectrichotomy2}}
\end{eqnarray*}
That contradiction completes the proof of the lemma.

\begin{lemma}\label{lemma:leastelement}
 Suppose $\N \in \FINITE$ and $\s \k = \s \n$ and $  \k \neq \n$ and
$\k \in \Stem$ and $\n \in \N$.  
Let $X$ be a finite nonempty subset of $\N$.  Then $X$ has a $\preceq$-least
element.  More formally,
$$\forall X \in \FINITE( X \subseteq \N \land X \neq \empty \imp
\exists p \in X\, \forall q \in X\,(p \preceq q)).$$
\end{lemma}

\noindent{\em Proof}.  The formula to be proved is stratified, giving 
$p$ and $q$ both index 0 and $X$ index 1.   Therefore we can proceed
by induction on finite sets $X$.  Because of the hypothesis that $X \neq \empty$,
the base case is immediate.  For the induction step, suppose
$X = Y \cup \{b\}$ with $b \not\in Y$ and $Y$ a finite set, and 
$X \subseteq \N$.  Then also $Y \subseteq \N$.  By Lemma~\sref{lemma:empty_or_inhabited},
$Y$ is empty or inhabited. If $Y = \empty$, then $b$ is the only element of 
$X$, and hence the least element of $X$.  So we may assume $Y$ is inhabited.
Then, by the induction hypothesis,
$Y$ has a $\preceq$-least element $r$.  By Theorem~\ref{theorem:prectrichotomy1},
we have $r \preceq b  \ \lor \ b \preceq r$.  We argue by cases.
\smallskip

Case~1:  $r \preceq b$.  Then $r$ is the desired $\preceq$-least member of $X$.
\smallskip

Case~2: $b \preceq r$.  Then $b$ is the desired $\preceq$-least member of $X$,
by Lemma~\ref{lemma:transitivity}.  
That completes the proof of the lemma.

\begin{lemma} [Transitivity of $\prec$] \label{lemma:prectrans} 
Suppose $\N \in \FINITE$ and $\s \k = \s \n$ and $  \k \neq \n$ and
$\k \in \Stem$ and $\n \in \N$.  Then for $x,y,z \in \N$ we have 
$$ x \prec  y \imp y \prec  z \imp x \prec  z.$$
\end{lemma}

\noindent{\em Proof}.  Suppose $x \prec y$ and $y \prec z$.  Then
\begin{eqnarray*}
x \preceq y  && \mbox{\qquad by the definition of $\prec$}\\
y \preceq z  && \mbox{\qquad by the definition of $\prec$}\\
x \preceq z  && \mbox{\qquad by Lemma~\ref{lemma:preceqtrans}}\\
x = z \imp z \preceq y  && \mbox{\qquad since $x \preceq y$}\\
x = z \imp y = z  && \mbox{\qquad by Theorem~\ref{theorem:prectrichotomy3}}\\
y \neq z  && \mbox{\qquad by the definition of $y \prec z$}\\
x \neq z  && \mbox{\qquad by the preceding two lines}\\
x \prec z && \mbox{\qquad by the definition of $\prec$}
\end{eqnarray*}
That completes the proof of the lemma.

\begin{lemma} \label{lemma:precpred} 
 Suppose $\N \in \FINITE$ and $\s \k = \s \n$ and $  \k \neq \n$ and
$\k \in \Stem$ and $\n \in \N$.   Let $x,y \in \N$ with $x \neq \n$.
Then $\s x \preceq \s y \imp x \preceq y$.
\end{lemma}

\noindent{\em Proof}. Assume all the hypotheses of the lemma, as well 
as $\s x \preceq \s y$.  We must prove $x \preceq y$.
Since $\N$ is finite, it has decidable equality.
Therefore $y = \n \ \lor \ y \neq \n$.  If $y = n$ then 
$x \preceq y$ by Lemma~\ref{lemma:precmax}, and we are done.  Therefore
we may assume $y \neq \n$.
\smallskip

We have 
\begin{eqnarray*}
x \preceq \s x  && \mbox{\qquad by Lemma~\ref{lemma:xpreceqsx}}\\
x \neq \s x   && \mbox{\qquad by Lemma~\ref{lemma:snneqn}}\\
\s x \preceq \s y && \mbox{\qquad by hypothesis}\\
x \preceq \s y && \mbox{\qquad by Lemma~\ref{lemma:preceqtrans}}\\
x\preceq \s y \iff x \preceq y \ \lor x = \s y && \mbox{\qquad by Lemma~\ref{lemma:preceqsuccessor}, since $y \neq \n$}\\
 x \preceq y \ \lor x = \s y  && \mbox{\qquad by the preceding two lines} 
\end{eqnarray*}
If $x \preceq y$ we are finished; so we may assume $x = \s y$.  Then 
$\s x \preceq \s y$ becomes $\s x \preceq x$.  Then
\begin{eqnarray*}
x = \s x && \mbox{\qquad by Theorem~\ref{theorem:prectrichotomy3}, since $x \preceq \s x \preceq x$}\\
x \neq \s x && \mbox{\qquad by Lemma~\ref{lemma:snneqn}}
\end{eqnarray*}
That contradiction completes the proof of the lemma.
 
\section{Multiplication of Church numbers}

In this section we define multiplication $x \ox y$  of Church numbers in 
such a way that
 it satisfies the ``defining'' laws $x \ox \ChurchZero = \ChurchZero$
 and $x \ox \s y = x \ox y \oplus  x$.  This requires knowing that either 
 successor is one-to-on on $\N$,  or that there is a double successor
 $\s \n = \s \k$ with $\k \neq \n$.  Then the further laws of 
 multiplication can be developed from those two, plus the decidability
 of equality on $\N$.  Hence, all the results of this section will be 
 valid when there is a double successor, and we use them later to show 
 that $\N$ cannot be finite.  But then, we still need multiplication to 
 interpret HA, so we need these results also in the case when successor is 
 one-to-one.   

\begin{lemma} \label{lemma:ChurchMultiplication}
Multiplication $x \ox y$ on $\N \times \N$ can be defined (as a function
of two variables) in INF and 
satisfies the following laws for all $x,y \in \N$:
\smallskip

(i) $x \ox \ChurchZero = \ChurchZero$
\smallskip

(ii) $x \ox  \s y = x \ox  y \oplus  x$.
\end{lemma}

\noindent{\em Remark}.  The lemma does not assume that $\N$ is finite
or that there is a double successor, or that there is no double successor.  
To prove it we have to give an ``agnostic'' definition of multiplication,
that works without any assumption of that sort.
\smallskip

\noindent{\em Proof}.
 We define multiplication as the intersection of
all sets
$Z$ of ordered triples $\langle x,y,z \rangle$ satisfying these 
conditions:
\begin{eqnarray*}
&& \forall y \in \N\,(\langle \ChurchZero,y,\ChurchZero \rangle \in Z) \\
&& \forall y \in \N\,(\langle x,y,z \rangle \in Z \land \neg \exists u\,(u \in \N \ \land \ u < y \ \land \ \s y = \s u) \imp
                             \langle x, \s y, z\oplus x \rangle \in Z)
\end{eqnarray*}
These formulas are stratified, giving $x,y,z$ index 0 and $Z$ 
index 5,  so the definition is legal in INF, and defines a
relation, which we write $x \ox  y = z$.  It remains to prove
that this relation is a function. 
 We will prove by induction on $y$ that for each 
$x$ there exists a unique $z$ such that $x \ox  y = z$. 
\smallskip

Base case: Existence: $x\ox  \ChurchZero = \ChurchZero$ by the first condition.
Uniqueness: $x \ox  \ChurchZero = z$ is only possibly by the first condition,
by Theorem~\ref{theorem:successoromitszero}.  That completes the base case.
\smallskip

Induction step: Existence: By the induction hypothesis, 
there exists $z$ such that $x \ox  y = z$.  Then 
by the second condition, $x \ox  \s y = z \oplus  x$.  
\smallskip

Uniqueness.  Suppose $x \ox  \s y = z$ and $x \ox  \s y = w$.
Then by the second condition, $z = x\ox  u \oplus  x$ and 
$w = x \ox  v \oplus  x$, where $\s u = \s v = \s y$.  Then 
$\neg\, (u < v)$ and $\neg\,(v < u)$, by the second condition.
Then by Lemma~\ref{lemma:trichotomy1} we have $u=v$.    Hence 
$z =  x \ox  v \oplus  x = x \ox  u \oplus  x =w$ as desired.  
That completes the induction step.
That completes the proof of the lemma.

\begin{lemma}\label{lemma:ChurchMultiplicationMaps}
$\forall x,y \in \N\,(x \ox y \in \N).$
\end{lemma}

\noindent{\em Proof}.  The formula is stratified, so we can prove it 
by induction on $y$.
\smallskip

Base case.  $x \ox \ChurchZero = \ChurchZero$ by Lemma~\ref{lemma:ChurchMultiplication}.
$\ChurchZero \in \N$ by Lemma~\ref{lemma:zeroN}.  That completes the base case.
\smallskip

Induction step.  Suppose 
\begin{eqnarray}
\forall x\,(x \ox y \in \N) \label{eq:E3688}
\end{eqnarray}  We have to prove
$$
\s y \in \N \imp\forall x\in \N\,( x \ox \s y \in \N)
$$
  Suppose $y,x \in \N$.  We have
\begin{eqnarray*}
x \ox \s y = x \ox y \oplus  x &&\mbox{\qquad By Lemma~\ref{lemma:ChurchMultiplication}}\\
x \ox y \in \N && \mbox{\qquad by (\ref{eq:E3688})}\\
\s y \in  && \mbox{\qquad by Lemma~\ref{lemma:successorN}}\\
x \ox y \oplus  x \in \N && \mbox{\qquad by Lemma~\ref{lemma:ChurchAdditionMaps}}\\
x \ox \s y \in \N && \mbox{\qquad by the preceding lines}
\end{eqnarray*}
That completes the induction step.  That completes the proof of the lemma.

\begin{lemma} \label{lemma:zerotimes}
For $y \in \N$ we have 
\begin{eqnarray*}
\ChurchZero \ox y = \ChurchZero  
\end{eqnarray*}
\end{lemma}

\noindent{\em Proof}.  By Lemma~\ref{lemma:ChurchMultiplication},
multiplication is well-defined and satisfies the laws in 
 Lemma~\ref{lemma:ChurchMultiplication}.  The formula is stratified, so we 
may use  induction on $y$.
\smallskip

Base case, $\ChurchZero \ox \ChurchZero = \ChurchZero$ by Lemma~\ref{lemma:ChurchMultiplication}.
\smallskip

Induction step, 
\begin{eqnarray*}
\ChurchZero \ox \s y = \ChurchZero \ox y \oplus  \ChurchZero  && \mbox{\qquad by Lemma~\ref{lemma:ChurchMultiplication}}\\
= \ChurchZero \ox y  && \mbox{\qquad  since $z \oplus  \ChurchZero = z$}\\
= \ChurchZero  && \mbox{\qquad by the induction hypthesis} 
\end{eqnarray*}
That completes the induction step.  That completes the proof of the lemma.

\begin{lemma} \label{lemma:successortimes}
For $x,y \in \N$ we have 
\begin{eqnarray*}
\s x \ox  y = x \ox  y \oplus  y.  
\end{eqnarray*}
\end{lemma}

\noindent{\em Proof}.  By Lemma~\ref{lemma:ChurchMultiplication},
multiplication is well-defined and satisfies the laws in 
 Lemma~\ref{lemma:ChurchMultiplication}.  The formula is stratified, so we 
may use induction on $y$.
\smallskip

Base case: 
\begin{eqnarray*}
\s x \ox  \ChurchZero = \ChurchZero && \mbox{\qquad by Lemma~\ref{lemma:ChurchMultiplication}}\\
\ChurchZero = x \ox \ChurchZero &&\mbox{\qquad by Lemma~\ref{lemma:ChurchMultiplication}}\\
\ChurchZero = x \ox \ChurchZero \oplus  \ChurchZero &&\mbox{\qquad by Lemma~\ref{lemma:ChurchZero_equation}}
\end{eqnarray*}
  That completes the base case.
\smallskip

Induction step:  Assume $\s x \ox  y = x \ox  y \oplus  y$.   We must prove
\begin{eqnarray*}
\s x \ox  \s y = x \ox  \s y \oplus  \s y. 
\end{eqnarray*} 
 
We have
\begin{eqnarray*}
\s x \ox  \s y &=&\s x \ox  y \oplus  \s x \mbox{\qquad\quad \ \ \  by definition of multiplication}\\
               &=& x \ox  y \oplus  y \oplus  \s x \mbox{\qquad\ \ by the induction hypothesis} \\
               &=& x \ox y \oplus  ( y \oplus  \s x) \mbox{\qquad by Lemma~\ref{lemma:ChurchAdditionAssociative}}\\
               &=& x \ox  y \oplus  (\s y \oplus  x)  \mbox{\qquad by Lemma~\ref{lemma:ChurchSuccessorShift}}\\
               &=& x \ox  y \oplus  (x \oplus  \s y)  \mbox{\qquad by Lemma~\ref{lemma:ChurchAdditionCommutative}}\\
               &=& x \ox  y \oplus   x \oplus  \s y  \mbox{\qquad\ \  by Lemma~\ref{lemma:ChurchAdditionAssociative}}\\
               &=& x \ox  \s y \oplus  \s y   \mbox{\qquad\quad\ \ \ by Lemma~\ref{lemma:ChurchMultiplication}}
\end{eqnarray*}
That completes the induction step. That completes the proof of the lemma.

\begin{lemma}\label{lemma:Church_leftdistrib}
For $x,y,z \in \N$ we have 
\begin{eqnarray*}
x \ox (y \oplus  z) &=& x \ox  y \oplus  x \ox  z 
\end{eqnarray*}
\end{lemma}
 
\noindent{\em Proof}. By induction on $x$.  
\smallskip

Base case. 
\begin{eqnarray*}
y \oplus  z \in \N && \mbox{\qquad by Lemma~\ref{lemma:ChurchAdditionMaps}}\\
\ChurchZero \ox (y \oplus  z) = \ChurchZero && \mbox{\qquad by Lemma~\ref{lemma:zerotimes}}\\
\ChurchZero \ox y = \ChurchZero && \mbox{\qquad by Lemma~\ref{lemma:zerotimes}}\\
\ChurchZero \ox z = \ChurchZero && \mbox{\qquad by Lemma~\ref{lemma:zerotimes}}\\
\ChurchZero \ox (y \oplus  z) = \ChurchZero \ox y \oplus  \ChurchZero \ox z &&\mbox{\qquad since $\ChurchZero\oplus \ChurchZero = \ChurchZero$}
\end{eqnarray*}
That completes the base case. 

Induction step.   
\begin{eqnarray*}
y \oplus  z \in \N && \mbox{\qquad by Lemma~\ref{lemma:ChurchAdditionMaps}}\\
x \ox y \in \N && \mbox{\qquad by Lemma~\ref{lemma:ChurchMultiplicationMaps}}\\
x \ox z \in \N && \mbox{\qquad by Lemma~\ref{lemma:ChurchMultiplicationMaps}}\\
\s x \ox (y\oplus z) = x \ox (y\oplus z) \oplus  (y\oplus z)  && \mbox{\qquad by Lemma~\ref{lemma:successortimes}}\\
                = x \ox  y \oplus  x \ox  z \oplus   y \oplus  z  &&\mbox{\qquad by the induction hypothesis}\\
                 = x \ox  y \oplus  (x \ox  z \oplus   y) \oplus  z  &&\mbox{\qquad by Lemma~\ref{lemma:ChurchAdditionAssociative}}\\                  = x \ox  y \oplus  (y \oplus  x \ox z)  \oplus  z  &&\mbox{\qquad by Lemma~\ref{lemma:ChurchAdditionCommutative}}\\              
               = (x \ox  y \oplus  y) \oplus  (x \ox  z \oplus  z) &&\mbox{\qquad by Lemma~\ref{lemma:ChurchAdditionAssociative}}\\         
                = \s x \ox  y \oplus  \s x \ox  z &&\mbox{\qquad by Lemma~\ref{lemma:successortimes}}
\end{eqnarray*}
That completes the proof of the lemma.

\begin{lemma}\label{lemma:Church_rightdistrib}
Suppose there is a double successor $\s \k = \s \n$ with $\k \in \Stem$ and $\k \neq \n$ and $\n \in \N$.
Then for $x,y \in \N$ we have 
\begin{eqnarray*}
(x \oplus  y) \ox  z &=& x \ox  z \oplus  y \ox  z 
\end{eqnarray*}
\end{lemma}
 
\noindent{\em Proof}. By induction on $z$.  The base case is immediate.  For the induction step we have
\begin{eqnarray*}
(x \oplus  y) \ox \s z  &=& (x\oplus y) \ox z  \oplus  (x \oplus  y) \mbox{\qquad\quad by Lemma~\ref{lemma:ChurchMultiplication}}\\
                &=& x \ox  z \oplus  y \ox  z \oplus   (x\oplus y) \mbox{\qquad by the induction hypothesis} \\
                &=& (x \ox  z \oplus  x) \oplus  (y \ox  z \oplus  y) \mbox{\quad by associativity and commutativity of $\oplus $} \\
                &=&  x \ox  \s z \oplus  y \ox  \s z  \mbox{\qquad\qquad\quad by Lemma~\ref{lemma:ChurchMultiplication}}
\end{eqnarray*}
That completes the proof of the lemma.

\begin{lemma}[Church multiplication associative]\label{lemma:ChurchMultiplicationAssociative} 
For $x,y,z \in \N$ we have 
\begin{eqnarray*}
   x \ox  ( y \ox z) &=& (x \ox y) \ox  z
\end{eqnarray*}
\end{lemma}

\noindent{\em Proof}.   
By induction on $y$.
\smallskip

Base case: 
\begin{eqnarray*}
x \ox  (\ChurchZero \ox  z) = x \ox  \ChurchZero && \mbox{\qquad by Lemma~\ref{lemma:zerotimes}}\\
= \ChurchZero  && \mbox{\qquad by Lemma~\ref{lemma:ChurchMultiplication}}\\
= (x \ox \ChurchZero) \ox  \ChurchZero && \mbox{\qquad by Lemma~\ref{lemma:ChurchMultiplication}}
\end{eqnarray*}
That completes the base case. 
\smallskip

Induction step:
\begin{eqnarray*}
x \ox  (\s y \ox  z)= x \ox  (y \ox  z \oplus  z) && \mbox{\qquad by Lemma~\ref{lemma:successortimes}} \\
                     =x \ox  (y \ox  z) \oplus  x \ox  z  && \mbox{\qquad by Lemma~\ref{lemma:Church_leftdistrib} } \\
                      = (x \ox  y) \ox  z) \oplus  x \ox  z && \mbox{\qquad by the induction hypothesis } \\
                       = (x \ox  y \oplus  x) \ox  z&& \mbox{\qquad  by Lemma~\ref{lemma:Church_rightdistrib} }  \\
                      =(x \ox  \s y) \ox  z && \mbox{\qquad by Lemma~\ref{lemma:ChurchMultiplication} } 
\end{eqnarray*}

\begin{lemma}[Church multiplication commutative]\label{lemma:ChurchMultiplicationCommutative} 
For $x,y \in \N$ we have 
\begin{eqnarray*}
   x \ox  y &=&  y \ox x
\end{eqnarray*}
\end{lemma}

\noindent{\em Proof}. By induction on $y$, which is legal since
the formula is stratified.
\smallskip

Base case. We have
\begin{eqnarray*}
x \ox \ChurchZero = \ChurchZero && \mbox{\qquad by Lemma~\ref{lemma:ChurchMultiplication}}\\
= \ChurchZero \ox x && \mbox{\qquad by Lemma~\ref{lemma:zerotimes}}
\end{eqnarray*}
That completes the base case.
\smallskip

Induction step. We have
\begin{eqnarray*}
\s x \ox y = x \ox y \oplus  y &&\mbox{\qquad by Lemma~\ref{lemma:successortimes}}\\
            = y \ox  x \oplus  y && \mbox{\qquad by the induction hypothesis}\\
 = y \ox  \s x.  && \mbox{\qquad by Lemma~\ref{lemma:ChurchMultiplication}}
\end{eqnarray*}
That completes the induction step. 
That completes the proof of the lemma.

\section{Successor and addition on the loop}

In this section we consider the map $f$ on the loop, defined by restricting 
Church successor to the loop.  We will show that $f$ is a permutation of the loop;
by the Annihilation Theorem then $\m f$ is the identity on the loop.  We then 
consider solutions $x$ of the equation $\n + x = \n$.  There is a solution,
and we show there is a $\preceq$-least solution $\m$.   If $\m$ were the 
order of $f$, we could reach a contradiction, proving that $\N$ is not 
finite. 
We assumed the Church counting axiom to reach that conclusion; but in the last
half of this section we prove, without the counting axiom, that the order
of successor at least exists.  That existence is not subsequently used, 
but we include it anyway.    

\begin{lemma}\label{lemma:loopaddition}
Suppose $\N$ is finite and  
 $\s \k = \s \n$ and $  \k \neq \n$ and
$\k \in \Stem$ and $\n \in \N$.  Then every element
of $\L(\n)$ has the form $\n \oplus  x$ for some $x \in \N$.
\end{lemma}

\noindent{\em Proof}.   Let 
$$X:= \{ \n \oplus  x : x \in \N\}.$$
The formula is stratified, so the definition is legal.
Then
\begin{eqnarray*}
\n \in X && \mbox{\qquad since $\n = \n \oplus  \ChurchZero$, by Lemma~\ref{lemma:ChurchZero_equation}}\\
z \in X \imp \s z \in X && \mbox{ \qquad since $\s(n\oplus x) = \n \oplus  \s x$, by Lemma~\ref{lemma:ChurchAddition_equation}} \\
\L(\n) \subseteq X && \mbox{\qquad by the definition of $\L(\n)$}
\end{eqnarray*}
That completes the proof of the lemma.

\begin{lemma}\label{lemma:mexists_helper}
  Suppose $\N \in \FINITE$ and $\s \k = \s \n$ and $  \k \neq \n$ and
$\k \in \Stem$ and $\n \in \N$.   Then there exists $\m \in \N$ such that $\n = \k \oplus  \m$. 
\end{lemma}

\noindent{\em Proof}.  
\begin{eqnarray*}
 \n \in \L(\n)   && \mbox{\qquad by Lemma~\ref{lemma:L1}}\\
 \n = \s p   && \mbox{\qquad for some $p \in \L(\n)$ by Theorem~\ref{theorem:looponto}}\\
p = \n\oplus  u && \mbox{\qquad for some $u \in \N$, by Lemma~\ref{lemma:loopaddition}}\\
\s p = \s (n \oplus  u) && \mbox{\qquad by the preceding line}\\
\n = \s(n\oplus u)  && \mbox{\qquad since $\s p = \n$}\\
\n = \n \oplus  \s u  && \mbox{\qquad by Lemma~\ref{lemma:ChurchAddition_equation}}\\
\n = \s \n \oplus  u  && \mbox{\qquad by Lemma~\ref{lemma:ChurchSuccessorShift}}\\
\n = \s \k \oplus  u  && \mbox{\qquad since $\s \n = \s \k$}\\
\n = \k \oplus  \s u  && \mbox{\qquad by Lemma~\ref{lemma:ChurchSuccessorShift}}
\end{eqnarray*}
Setting $\m:= \s u$ we have $\n = \k \oplus  \m$.  That completes the proof of the lemma.

\begin{lemma}\label{lemma:mexists}
  Suppose $\N \in \FINITE$ and $\s \k = \s \n$ and $  \k \neq \n$ and
$\k \in \Stem$ and $\n \in \N$.   Then there exists $\m \in \N$ such that $\n = \k \oplus  \m$ 
and $\m$ is the $\preceq$-least number with that property.  Explicitly,
$$ \forall p \in \N\,(\n = \k \oplus  p \imp \m \preceq p ).$$
\end{lemma}

\noindent{\em Proof}. Define
$$X = \{\ x : x \in \N \ \land \ \n = \k \oplus  x \ \}.$$
By Lemma~\ref{lemma:mexists_helper},  $X$ is inhabited.  
Since $\N$ is finite, it has decidable equality, by Lemma~\sref{lemma:finitedecidable}.
Therefore $X$ is a separable subset of $\N$.   By Lemma~\ref{lemma:separablefinite},
$X \in \FINITE$. 
By Lemma~\ref{lemma:leastelement}, $X$ has a $\preceq$-least element.
That completes the proof of the lemma.

\begin{lemma}\label{lemma:xsmapsloop}
 Suppose $\N \in \FINITE$ and $\s \k = \s \n$ and $  \k \neq \n$ and
$\k \in \Stem$ and $\n \in \N$.  Then $x\s : \L(\n) \to \L(\n)$.  That is,
$$ \forall x \in \N\, \forall y\,( y \in \L(\n) \imp (x \s y \in \L(\n)).$$
\end{lemma}

\noindent{\em Proof}.  To stratify the formula, we give $y$  index 0 and $x$ index 6.
Then $x\s y$ gets index 0, so the two occurrences of $\L(\n)$ could get the same index,
but since  
$\L(\n)$ is a parameter we do not even have to assign $\L(\n)$ an index. 
Since the formula is stratified, we may prove it by finite induction on $x$.
\smallskip

Base case, $x = \ChurchZero$.  Suppose $y \in \L(\n)$.  Then 
\begin{eqnarray*}
y \in \N && \mbox{\qquad by Lemma~\ref{lemma:LN}}\\
 \ChurchZero \s  y = y &&\mbox{\qquad by Lemma~\ref{lemma:ChurchZero_equation}}\\
  \ChurchZero \s   y \in \L(\n) && \mbox{\qquad by the preceding two lines}
 \end{eqnarray*}
 That completes the base case.
 \smallskip
 
 Induction step.  Suppose $x \neq \n$ and  $x \in \N$ and $y \in \L(\n)$.
 We must prove $ (\s x) \s y \in \L(\n)$.  We have
 \begin{eqnarray*}
 x \s  y \in \L(\n)  && \mbox{\qquad by the induction hypothesis} \\
 (\s x) \s y = \s (x \s y) && \mbox{\qquad by Theorem~\ref{theorem:successorequation}}\\
 \s (x \s y) \in \L(\n) && \mbox{\qquad by Lemma~\ref{lemma:L1}}\\
 (\s x) \s y \in \L(\n) && \mbox{\qquad by the preceding two lines}
 \end{eqnarray*}
 That completes the induction step.  That completes the proof of the lemma.
 
\begin{lemma} [Loop closed under addition]\label{lemma:additionmapsloop}
Suppose $\N$ is finite and  
 $\s \k = \s \n$ and $  \k \neq \n$ and
$\k \in \Stem$ and $\n \in \N$.  Then
$$ x \in \L(\n) \imp y \in \N \imp x \oplus y \in \L(\n).$$
\end{lemma}

\noindent{\em Proof}.  The displayed formula is stratified, giving $x$ and $\n$ index 
0, since $L(\n)$ and $\N$ are parameters.  So we may prove it by  induction
on $y$.
\smallskip

Base case, $y = \ChurchZero$. Then 
\begin{eqnarray*}
x \in \N && \mbox{\qquad by Lemma~\ref{lemma:LN}}\\
x \oplus \ChurchZero  = x && \mbox{\qquad by Lemma~\ref{lemma:ChurchZero_equation}} \\
x \oplus y  = x && \mbox{\qquad since $y = \ChurchZero$} \\
x \oplus y \in \L(\n) && \mbox{\qquad since $x \in \N$}
 \end{eqnarray*}
\smallskip

Induction step.  Suppose  $x \in \L(\n)$ and $\s y \in \L(\n)$ and $x \oplus y \in \L(\n)$.
We must prove $x \oplus \s y \in \L(\n)$. We have
\begin{eqnarray*}
x \oplus y \in \L(\n)  && \mbox{\qquad by the induction hypothesis}\\
\s (x \oplus y) \in \L(\n) && \mbox{\qquad by Lemma~\ref{lemma:L1}}\\
x \oplus \s y  \in \L(\n) && \mbox{\qquad by Lemma~\ref{lemma:ChurchAddition_equation}}
\end{eqnarray*}
That completes the induction step.
That completes the proof of the lemma.

\begin{lemma} \label{lemma:knm} 
Suppose $\N$ is finite and  
 $\s \k = \s \n$ and $  \k \neq \n$ and
$\k \in \Stem$ and $\n \in \N$.  Suppose $\k \oplus \m = \n$.
Then $\n \oplus \m = \n$.
\end{lemma}

\noindent{\em Proof}.   We have
\begin{eqnarray*}
\k \oplus \m = \n  && \mbox{\qquad by hypothesis}\\
\s (\k \oplus \m) = \s \n  && \mbox{\qquad by the previous line}\\
\k \in \N    && \mbox{\qquad by Lemma~\ref{lemma:SN}}\\
\k \oplus \s \m = \s \n && \mbox{\qquad by Lemma~\ref{lemma:ChurchAddition_equation}}\\
\s \k \oplus  \m = \s \n && \mbox{\qquad by Lemma~\ref{lemma:ChurchSuccessorShift}}\\
\s \k = \s \n  && \mbox{\qquad by hypothesis}\\
\s \n \oplus   \m = \s \n && \mbox{\qquad by the preceding two lines} \\
\n \oplus \s \m = \s \n && \mbox{\qquad by Lemma~\ref{lemma:ChurchSuccessorShift}}\\
\s (n \oplus \m) = \s \n && \mbox{\qquad by Lemma~\ref{lemma:ChurchAddition_equation}}\\
\n \oplus \m \in \L(\n) && \mbox{\qquad by Lemma~\ref{lemma:additionmapsloop}}\\
\n \oplus \m = \n  && \mbox{\qquad by Theorem~\ref{theorem:looponeone}}
\end{eqnarray*}
That completes the proof of the lemma.

\begin{lemma} \label{lemma:kplusm}
Suppose $\N \in \FINITE$ and $\s \k = \s \n$ and $  \k \neq \n$ and
$\k \in \Stem$ and $\n \in \N$.  Then for all $x \in \N$,
$$ x \neq \ChurchZero \imp \n = \n \oplus x \imp \n = \k \oplus x.$$
\end{lemma}

\noindent{\em Proof}.  Suppose $\n = \n \oplus x$. Then
\begin{eqnarray*}
\s \n = \s(\n \oplus x)  && \mbox{\qquad since $\n = \n \oplus x$}\\
\s \n = \n  \oplus  \s x  && \mbox{\qquad  by Lemma~\ref{lemma:ChurchAddition_equation}}\\
\s \n = \s \n \oplus x    && \mbox{\qquad by Lemma~\ref{lemma:ChurchSuccessorShift}}\\
\s \n = \s \k \oplus x     && \mbox{\qquad since $\s \n = \s \k$}\\
\k \in \N                 && \mbox{\qquad by Lemma~\ref{lemma:SN}, since $k \in \Stem$}\\
\s \n = \k \oplus \s x && \mbox{\qquad by Lemma~\ref{lemma:ChurchSuccessorShift}}\\
\s \n = \s(\k \oplus x)&& \mbox{\qquad  by Lemma~\ref{lemma:ChurchAddition_equation}}\\
\n \in \L(\n)   && \mbox{\qquad by Lemma~\ref{lemma:L1}}\\  
\s \n \in \L(\n)   && \mbox{\qquad by Lemma~\ref{lemma:L1}}\\  
\s \k \in \L(\n)   && \mbox{\qquad since $\s \k = \s \n$}\\
x \neq \ChurchZero  && \mbox{\qquad by hypothesis}\\
x = \s r            && \mbox{\qquad for some $r \in \N$, by Lemma~\ref{lemma:predecessor}}\\
\k \oplus x = \k \oplus \s r && \mbox{\qquad since $x = \s r$}\\
\k \oplus x = \s \k \oplus r  && \mbox{\qquad by Lemma~\ref{lemma:ChurchSuccessorShift}}\\
\s \k \oplus r  \in \L(\n)  && \mbox{\qquad by Lemma~\ref{lemma:additionmapsloop}}\\
\k = \n \oplus x        && \mbox{\qquad by Theorem~\ref{theorem:looponeone}, since $\s \n = \s(\k \oplus x)$}
\end{eqnarray*}
That completes the proof of the lemma.

\begin{definition} \label{definition:orderL}
Suppose $\N$ is finite and  
 $\s \k = \s \n$ and $  \k \neq \n$ and
$\k \in \Stem$ and $\n \in \N$.  Then we define
the {\bf order of successor on the loop} to be the $\preceq$-least Church number $\q$ such that
$\q \s$  is the identity on $\L(\n)$.  That is, 
\begin{eqnarray*} 
&&\forall x \in \L(\n)\,( \q \s x = x) \\
&& \forall r\in \N\,( r \prec \q \imp \neg\, \forall x \in \L(\n)\,(r \s x = x))
\end{eqnarray*}
For short we call $\q$ the ``order of $\L(\n)$'' or the ``order of the loop.''
\end{definition}
We shall show below that there actually exists such a number $\q$.   That, of course,
requires a proof, not just a definition.  Of course, assuming the Church
counting axiom, it is easy to prove that $\m$ is the order of the loop,
but we shall prove without the counting axiom that the order is well-defined.

\begin{lemma}\label{lemma:xsmapsN} Let $x\in \N$ and $y \in \N$. 
Then $x\s y \in \N$.
\end{lemma}

\noindent{\em Proof}.  The formula is stratified, giving $x$ index 6 
and $y$ index 0, with $\N$ as parameter.  Therefore we may prove it 
by induction on $x$.
\smallskip

Base case.  We  have 
\begin{eqnarray*}
\ChurchZero \s y = y  && \mbox{\qquad by Lemma~\ref{lemma:zeroAp}}\\
\ChurchZero \s y \in \N  && \mbox{\qquad since $y\in\N$}
\end{eqnarray*}
That completes the base case.
\smallskip

Induction step.  We have
\begin{eqnarray*}
x \s y \in \N  && \mbox{\qquad by the induction hypothesis}\\
x \in \FUNC  && \mbox{\qquad by Lemma~\ref{lemma:Churchnumbersarefunctions}}\\
\s(x \s y) \in \N && \mbox{\qquad by Lemma~\ref{lemma:successorN}, since $y \in \N$ }\\
\s x \s y = \s (x \s y) && \mbox{\qquad by Theorem~\ref{theorem:successorequation}}\\
\s x \s y \in \N && \mbox{\qquad by the preceding two lines}
\end{eqnarray*}
That completes the induction step.  That completes the proof of the theorem.

\begin{lemma} \label{lemma:successorflip} Let $t \in \N$ and $q \in \N$. 
Then 
$$ q \s (\s t) = \s (q \s t).$$
\end{lemma}

\noindent{\em Remark}.  Intuitively, both sides refer to successor applied
$q$ plus one times to $t$.
\smallskip

\noindent{\em Proof}.  We have 
\begin{eqnarray*}
\ChurchZero \s t = t  && \mbox{\qquad by Lemmas~\ref{lemma:zeroAp}}\\
\s (\ChurchZero \s t) = \s t && \mbox{\qquad by the preceding line}\\
q \s t \in \N && \mbox{\qquad by Lemma~\ref{lemma:xsmapsN}}\\
\s t \in \N && \mbox{\qquad by Lemma~\ref{lemma:successorN}}\\
\s t \in \FUNC && \mbox{\qquad by Lemma~\ref{lemma:Churchnumbersarefunctions}}\\
Rel(\s t)  && \mbox{\qquad by Lemma~\ref{lemma:Churchnumbersarerelations}}\\ 
\s \ChurchZero \s t = \s t && \mbox{\qquad by Lemma~\ref{lemma:ApOne}}\\
\s t = \s \ChurchZero \s t  && \mbox{\qquad by the preceding line}\\
q \s (\s t) = q \s (\s \ChurchZero \s t) && \mbox{\qquad by the preceding line}\\
    = (q \oplus \s \ChurchZero) \s t && \mbox{\qquad by Lemma~\ref{lemma:doubleiteration} with $X = \N$ and $f = \s$}\\
    = (\s q \oplus \ChurchZero) \s t && \mbox{\qquad by Lemma~\ref{lemma:ChurchSuccessorShift}}\\
 = \s q \s t  && \mbox{\qquad by Lemma~\ref{lemma:ChurchZero_equation}}\\
 = \s (q \s t) &&  \mbox{\qquad by Theorem~\ref{theorem:successorequation}}
\end{eqnarray*}
That completes the proof of the lemma.

\begin{lemma} \label{lemma:orderq_helper}
Suppose $\N$ is finite and  
 $\s \k = \s \n$ and $  \k \neq \n$ and
$\k \in \Stem$ and $\n \in \N$.  Suppose   $q \in \N$ and $q \neq \ChurchZero$ and $q\s \n = \n$. 
Then
  $q \s $ is the identity on $\L(\n)$.
\end{lemma}

\noindent{\em Proof}.   Define 
$X = \{ x \in \L(n):  q \s x = x \}$
The formula is stratified, giving $x$ index 0 and $q$ index 6, so the 
definition is legal.  By hypothesis, $\n \in X$.  I say that $X$ is closed
under successor.  Suppose $x \in X$.  Then $x \in \L(\n)$ and $q \s x = x$. 
\begin{eqnarray*}
 q \s ( \s x)      = \s (q \s x)  && \mbox{\qquad by Lemma~\ref{lemma:successorflip}}\\
 q \s (\s x)    = \s x   && \mbox{\qquad since $q \s x = x$} \\
 \s x \in \L(\n)  && \mbox{\qquad by Lemma~\ref{lemma:L1}}\\
 \s x \in X  && \mbox{\qquad by definition of $X$}
\end{eqnarray*}
That completes the proof that $X$ is closed under successor.
Then by definition of $\L(\n)$,  we have $\L(\n) \subseteq X$.  
That completes the proof of the lemma.

\begin{lemma} \label{lemma:orderq_helper2}
Suppose $\N$ is finite and  
 $\s \k = \s \n$ and $  \k \neq \n$ and
$\k \in \Stem$ and $\n \in \N$.  Suppose   $q \in \N$ and $q \neq \ChurchZero$ and 
 $t \in \L(\n)$ and $q\s t = t$. 
Then
  $q \s $ is the identity on $\L(\n)$.  That is,
  $$\forall x \in \L(\n)\, (q \s x = x).$$
\end{lemma}

\noindent{\em Proof}.  The formula is stratified, giving $x$ and $t$ index 0 and 
$q$ index 6.  We can therefore prove it by ``loop induction.'' That is,
we show that the set of $t$ for which the lemma holds contains $\n$ and 
is closed under successor.   That set, explicitly, is 
$$ Z:= \{\ t\in \L(\n):  q \s t = t \imp \forall x \in \L(\n)\,(q \s x = x)\ \}.$$
The formula defining $Z$ is stratified, giving $q$ index 6 and $x$ and $t$ 
index 0.  $\L(\n)$ is a parameter.  Therefore $Z$ can be defined in INF.
\smallskip

$Z$ contains $\n$, by Lemma ~\ref{lemma:orderq_helper}.  It remains to 
show $Z$ is closed under successor.  Suppose $t \in Z$.  We must show 
$\s t \in \Z$.  Suppose $q \s (\s t) =  \s t $, and let $x \in \L(\n)$ 
be given.  We must show $\q \s x = x$.  We have 
\begin{eqnarray*}
q \s (\s t) = \s (q \s t)  && \mbox{\qquad by Lemma~\ref{lemma:successorflip}}\\
\s (q \s t) = \s t  && \mbox{\qquad since $q \s (\s t) =  \s t $}\\
q \s t = t   && \mbox{\qquad by Theorem~\ref{theorem:looponeone}}\\
\forall x \in \L(\n)\, (q x = x) && \mbox{\qquad by the induction hypothesis}
\end{eqnarray*}
That completes the induction step.  That completes the proof of the lemma.

\begin{lemma} \label{lemma:successoronloop}
Suppose $\N$ is finite and  
 $\s \k = \s \n$ and $  \k \neq \n$ and
$\k \in \Stem$ and $\n \in \N$.  Let $f$ be Church successor
restricted to $\L(\n)$.  Then for every $q \in \N$ and $x \in \L(\n)$,
$$ qfx = q \s x.$$
\end{lemma}

\noindent{\em Proof}.  The displayed formula is stratified, giving $x$ 
index 0, $f$ index 3, and $q$ index 6.  Therefore we can prove it by induction 
on $q$.
\smallskip

Base case, $\ChurchZero f x = x$ and $\ChurchZero \s x = x$, by Lemma~\ref{lemma:zeroAp}.
Therefore $\ChurchZero f x = \ChurchZero \s x$.  That completes the base case.
\smallskip

Induction step.  We have 
\begin{eqnarray*}
f:\L(\n) \to \L(\n) && \mbox{\qquad by Lemma~\ref{lemma:L1}}\\
Rel(f) \ \land \ f \in \FUNC  && \mbox{\qquad since $f$ is a subset of the graph of Church successor} \\
\s q f x = f(q f x) && \mbox{\qquad by Theorem~\ref{theorem:successorequation}}\\
\s q f x  = f (q \s x)  && \mbox{\qquad by the induction hypothesis, $q f x = q \s x$}\\
 \s q \s x) =\s (q \s x)   && \mbox{\qquad by Theorem~\ref{theorem:successorequation}} \\
 q \s x \in \L(\n)  && \mbox{\qquad by Lemma~\ref{lemma:xsmapsloop}} \\
f(q \s x) = \s (q \s x) && \mbox{\qquad since $f$ is the restriction of $\s$ to $\L(\n)$} \\
\s q f x = \s (q \s x)  && \mbox{\qquad since $\s q f x = f(q f x)$}\\
\s q f x = \s q \s x   && \mbox{\qquad since $ \s q \s x) =\s (q \s x) $}
 \end{eqnarray*}
That completes the induction step.  That completes the proof of the lemma.
 
\begin{lemma} \label{lemma:orderq}  Suppose $\N$ is finite and  
 $\s \k = \s \n$ and $  \k \neq \n$ and
$\k \in \Stem$ and $\n \in \N$.   Then the order of successor on the loop exists.
\end{lemma}

\noindent{\em Proof}. 
\begin{eqnarray*}
\n = \k \oplus \m  &&\mbox{\qquad for some $\m \in\N$, by  Lemma~\ref{lemma:mexists}}\\
\n = \n \oplus \m  && \mbox{\qquad by Lemma~\ref{lemma:knm}}\\
 \m \s \mbox{  is the identity on $\L(\n)$ }&&\mbox{\qquad by the Annihilation Theorem}\\
 \m \s \n = \n && \mbox{\qquad by the previous line}
 \end{eqnarray*}
 Define
$$X := \{ q \in \N  : q \s \n = \n\}.$$
The formula defining $X$ is stratified, giving $q$ index 6 and $\n$ index 0, so 
$X$ can be defined in INF.   Define
$f$ to be Church successor restricted to $\L(\n)$.   Then 
\begin{eqnarray*}
f:\L(\n) \to \L(\n) && \mbox{\qquad by Lemma~\ref{lemma:L1}}\\
f \mbox{\ is an injection} && \mbox{\qquad by Theorem~\ref{theorem:looponeone}} \\
\m f \mbox{\ is the identity on $\L(\n)$} && \mbox{\qquad by the Annihilation Theorem} \\
x \in \L(\n)\imp \m f x = \m \s x && \mbox{\qquad by Lemma~\ref{lemma:successoronloop}}\\
\m \s \mbox{\ is the identity on $\L(\n)$} && \mbox{\qquad by the preceding two lines} \\
\m \in X  && \mbox{\qquad by the definition of $X$}\\
\N \in \DECIDABLE && \mbox{\qquad  by Lemma~\sref{lemma:finitedecidable}} \\
X \mbox{\  is a separable subset of $\N$} && \mbox{\qquad by definition of separable}\\
X \in \FINITE && \mbox{\qquad by Lemma~\ref{lemma:separablefinite}} \\
X \mbox{\ has a $\preceq$-least element} && \mbox{\qquad by Lemma~\ref{lemma:leastelement}, since $\m \in X$}
\end{eqnarray*}
Let $\q$ be that element.  By Lemma~\ref{lemma:orderq_helper},
$\q \s$ is the identity on $\L(\n)$.   Now let $r \prec q$,  and suppose $r \s$ is the 
identity on $\L(\n)$.  Then $r \in X$, contradiction, since $\q$ is the $\preceq$-least element of $X$.
Therefore $\q$ is the order of $\s$ on $\L(\n)$, as claimed.

\section{The Church counting axiom}

The ``Church counting axiom'' expresses the idea that iterating the Church successor function 
$j$ times starting from $\ChurchZero$ leads to the Church number $j$.  The 
formula expressing this fact is not stratified,  since $j$ as a function must get an 
index six higher than $j$ as an ``object.''    Hence if one wishes to use this principle,
it must be assumed as a new axiom.   Here is that 
axiom:

\begin{definition}\label{definition:countingaxiom}
The {\bf Church counting axiom} is 
$$\forall j \in \N\,(j\s 0 = j).$$
\end{definition}
A similar axiom was introduced by Rosser \cite{rosser1953}.
Rosser's axiom is stated using the finite Frege cardinals.
It says that $\{x\in \FregeN : x < p\}$ belongs to the 
cardinal number $p$.   In the last section of this paper,
we will prove that the two counting axioms are equivalent.
Orey proved \cite{orey1964} that the Rosser counting axiom is 
not provable in NF (unless, of course, NF is inconsistent).
Therefore our result shows that the same is true of the Church
counting axiom. 

The main result of this paper is that INF plus the Church counting axiom proves that
$\N$ is infinite and Church successor is one-to-one on $\N$.  Whether this can 
be proved with intuitionistic logic and 
without the Church counting axiom we do not know.  Our proof
appears to require that the order of (successor on) the loop be $\m$, and 
we failed to prove that without the Church counting axiom.  In this section, we 
present a proof of that fact using the Church counting axiom.

\begin{lemma} [loop counting]\label{lemma:loopcounting}
Assume the Church counting axiom.  
Suppose $\N \in \FINITE$ and $\s \k = \s \n$ and $  \k \neq \n$ and
$\k \in \Stem$ and $\n \in \N$.  Let $q \in \N$.  Then 
$$ q\s\n = \n \oplus q.$$
\end{lemma}

\noindent{\em Proof}.  We have
\begin{eqnarray*}
\n = \n \s \ChurchZero && \mbox{\qquad by the Church counting axiom}\\
q \s \n =  q \s (\n \s \ChurchZero) && \mbox{\qquad by the preceding line}\\
 q \s \n = (q \oplus  \n)\s \ChurchZero  && \mbox{\qquad by Lemma~\ref{lemma:doubleiteration}}\\
q \s \n  = (\n \oplus q)\s \ChurchZero && \mbox{\qquad by Lemma~\ref{lemma:ChurchAdditionCommutative}}\\
\n \oplus q \in \N  && \mbox{\qquad by Lemma~\ref{lemma:ChurchAdditionMaps}}\\
q \s \n  = \n \oplus q    && \mbox{\qquad by the Church counting axiom}
\end{eqnarray*}

\begin{theorem}[Order of successor on the loop is $\m$] \label{theorem:orderm}
Assume the Church counting axiom.   
 Suppose $\N \in \FINITE$ and $\s \k = \s \n$ and $  \k \neq \n$ and
$\k \in \Stem$ and $\n \in \N$.  Suppose $\m \in \N$ and $\n = \k \oplus  \m$,
and $\m$ is the $\preceq$-least solution of $\n = \k \oplus \m$.  
Suppose $q \in \N$ and $q \neq \ChurchZero$ and  $q \s$ is the identity on $\L(\n)$.   Then $\m \preceq q$.
\end{theorem}

\noindent{\em Proof}.  Suppose $q s$ is the identity on $\L(\n)$ and $q \neq \ChurchZero$. 
We must show $\m \preceq q$.  We have
\begin{eqnarray*}
q \s \n = \n \oplus q  && \mbox{\qquad by Lemma~\ref{lemma:loopcounting} and the Church counting axiom}\\
q \s \n = \n   && \mbox{\qquad since $q \s$ is the identity on $\L(\n)$}\\
\n = \n \oplus q    && \mbox{\qquad by the preceding two lines}\\
\n = \k \oplus q    && \mbox{\qquad by Lemma~\ref{lemma:kplusm}}\\
\m \preceq q   && \mbox{\qquad since $\m$ is the least solution of $\n = \k \oplus \m$}
\end{eqnarray*}
That completes the proof of the theorem.

\section{Church counting implies $\N$ is not finite}

Now we prove a series of lemmas under the hypothesis that $\N$ is finite.
With only that hypothesis,  results
proved earlier under the additional hypothesis that there is a double successor are not 
applicable;  without careful attention to the hypothesis, the reader might
get a sense of {\em deja vu}.

\begin{lemma} \label{lemma:stemseparable_helper}
Suppose $\N$ is finite and $x \in \Stem$.  Then 
$$\exists y \in \N\, (\s y = \s x \ \land \ y \neq x) \ \lor \ \neg \exists y \in \N\, (\s y = \s x \ \land \ y \neq x).$$
\end{lemma}

\noindent{\em Proof}.  Since $\N$ is finite,  it has decidable equality, 
by Lemma~\sref{lemma:finitedecidable}.  Then define
$$ R := \{ \ \langle x,y \rangle \in \N \times \N : \s y = \s x\ \ \land \ y \neq x\}.$$
The formula is stratified, giving $x$ and  $y$ index 0, with $\N \times \N$ as
a parameter.  Since $\N$ has decidable equality, $R$
 is a decidable relation $R$ on $\N$.   The 
conclusion of the lemma then follows from 
Lemma~\sref{lemma:boundedquantification}.   

\begin{lemma}\label{lemma:stemseparable}
Suppose $\N$ is finite.  Then $\Stem$ is a separable subset of $\N$.
\end{lemma}

\noindent{\em Proof}. We have to prove 
$$ \forall x \in \N\, (x \in \Stem \ \lor \ x \not \in \Stem).$$
That formula is stratified, giving $x$ index 0, with $\Stem$ a parameter.
We can therefore proceed by induction on $x$.
\smallskip

Base case. $\ChurchZero \in \Stem$, by Lemma~\ref{lemma:S1}.  Hence
$\ChurchZero \in \Stem \ \lor \ \ChurchZero \not \in \Stem$.  That completes
the base case.
\smallskip

Induction step.  
We must prove 
$$\s x \in \Stem \ \lor \ \s x \not\in \Stem.$$
The induction hypothesis is $x \in \Stem \ \lor \ x \not\in \Stem$.
We argue by cases accordingly. 
\smallskip

Case~1, $x \in \Stem$.  By Lemma~\ref{lemma:stemseparable_helper}, we 
have 
$$\exists y \in \N\, (\s y = \s x \ \land \ y \neq x) \ \lor \ \neg \exists y \in \N\, (\s y = \s x\ \land \ y \neq x).$$
We argue by cases accordingly.
\smallskip

Case~1a, $\exists y \in \N\, (\s y = \s x\ \land \ y \neq x)$.  Then $\s x \not \in \Stem$,
by Lemma~\ref{lemma:Soneone}.  That completes Case~1a.
\smallskip

Case~1b, $\neg \exists y \in \N\, (\s y = \s x\ \land \ y \neq x)$.  Then $\s x \in \Stem$,
by Lemma~\ref{lemma:S1}. That completes Case~1b.  That completes Case~1.
\smallskip

Case~2, $x \not\in \Stem$.  Then $\s x \not\in \Stem$, by Lemma~\ref{lemma:Spred}.
That completes Case~2.  That completes the induction step.
That completes the proof of the lemma.

\begin{lemma} \label{lemma:kinstem}
Suppose $\N$ is  finite.  Then there exists a double successor $\s \k = \s n$
with $\k \neq \n$ and $\k \in \Stem$.
\end{lemma}

\noindent{\em Proof}. We have
\begin{eqnarray*}
\Stem \mbox{\ is a separable subset of $\N$} &&\mbox{\qquad by Lemma~\ref{lemma:stemseparable}}\\
\N \in \FINITE  && \mbox{\qquad by hypothesis}\\
\Stem \in \FINITE && \mbox{\qquad by Lemma~\sref{lemma:separablefinite}}
\end{eqnarray*}
We would like to identify $\k$ as the maximal element of the finite set 
$\Stem$, but that is not a one-line proof, as we do not have a linear ordering
on $\Stem$ without assuming $\k \in \Stem$, which is what we are trying to prove,
so ``maximal'' makes no sense.
\smallskip

We avoid the need for a linear ordering as follows.  Define
$$ R:= \{\ \langle y,x \rangle \in \N \times \N:  x \in \Stem \ \land \ x \neq y \ \land \ \s x = \s y \ \}. $$
The formula is stratified, so the definition is legal. 
Then 
\begin{eqnarray*}
\Stem \mbox{ \ is a separable subset of $\N$} &&\mbox{\qquad by Lemma~\ref{lemma:stemseparable}}\\
\N \in \DECIDABLE && \mbox{\qquad by  Lemma~\sref{lemma:finitedecidable}}\\
R \mbox{\  is a decidable relation on $\N$}
   && \mbox{\qquad by the preceding lines}
\end{eqnarray*}
 Define
$$Z:=  \{ \ x \in \N : \exists y \in \N\, \langle y, x \rangle \in  R \ \}.$$
{\em Remark}.   $Z$ is the set of $x \in \Stem$ such that $\s x$ is a 
double successor. 
\smallskip
 
By Lemma~\ref{lemma:boundedquantification},
$Z \in \FINITE$.  
By Lemma~\sref{lemma:empty_or_inhabited}, $Z$ is empty 
or inhabited.  We argue by cases accordingly.
\smallskip

Case~1, $Z = \empty$.  Then there is no $x \in \Stem$ such that $\s x$ is  a
double successor.  Then 
\begin{eqnarray*}
\ChurchZero \in \Stem && \mbox{\qquad by Lemma~\ref{lemma:S1}}
\end{eqnarray*}
I say that $\Stem$ is closed under successor.  Suppose $x \in \Stem$;
we must show $\s x \in \Stem$.  By Lemma~\ref{lemma:S1}, it suffices
to show that $\s x$ is not a double successor; that is, it suffices
to show that 
$$ \forall v \in \N\,(\s x = \s v \imp x = v).$$
Let $v \in \N$ and $\s x = \s v$; we must show $x = v$.  Since $\N$ 
has decidable equality, we may prove that by contradiction.  Suppose $x \neq v$.
Then $x \in \Stem$ and $\s x$ is a double successor, so $x \in Z$.
But that contradicts the hypothesis $Z = \empty$ of Case~1.  That 
completes the proof that $\Stem$ is closed under successor.  Then

\begin{eqnarray*}
\N \subseteq \Stem  && \mbox{\qquad  by the definition of $\N$}\\
\Stem \subseteq \N  && \mbox{\qquad by Lemma~\ref{lemma:S1}}\\
\N = \Stem  && \mbox{\qquad by the preceding two lines}\\
\mbox{ Church successor is 
one-to-one on $\Stem$} && \mbox{\qquad by  Lemma~\ref{lemma:Soneone}}\\
\mbox{ Church successor is 
one-to-one on $\N$} &&\mbox{\qquad by the preceding two lines}\\
\mbox{ChurchSuccessor is not onto $\N$} && \mbox{\qquad by Theorem~\ref{theorem:successoromitszero}}\\
\mbox{$\N$ is infinite} && \mbox{\qquad by Definition~\sref{definition:infinite}} \\
\neg\, \N \in \FINITE  && \mbox{\qquad by Theorem~\sref{theorem:infiniteimpliesnotfinite}}
\end{eqnarray*}
But that contradicts the hypothesis that $\N$ is finite.  That completes
Case~1. 
\smallskip

Case~2, $Z$ is inhabited.  Then there exists some $x \in \Stem$ such
that $\s x$ is a double successor.  That completes the proof of the lemma.

\begin{lemma} \label{lemma:nplusm}
 Suppose $\N \in \FINITE$ and $\s \k = \s \n$ and $  \k \neq \n$ and
$\k \in \Stem$ and $\n \in \N$.  Then 
$$  \forall x \in \N\,( x \neq \ChurchZero \imp \n = \k \oplus x \imp \n = \n \oplus x).$$
\end{lemma}

\noindent{\em Proof}.  We have
\begin{eqnarray*}
\k \in \N && \mbox{\qquad by Lemma~\ref{lemma:SN}, since $\k \in \Stem$}\\
\n = \k \oplus x  && \mbox{\qquad by hypothesis}\\ 
\s \n = \s (\k \oplus x)  && \mbox{\qquad by the previous line}\\
\s \n = \k \oplus \s x  && \mbox{\qquad by Lemma~\ref{lemma:ChurchAddition_equation}}\\
\s \n= \s \k \oplus x  && \mbox{\qquad by Lemma~\ref{lemma:ChurchSuccessorShift}}\\
\s \n= \s \n \oplus x   && \mbox{\qquad since $\s \k = \s \n$}\\
\s \n= \n \oplus \s x  && \mbox{\qquad by Lemma~\ref{lemma:ChurchSuccessorShift}}\\
\s \n = \s (\n \oplus x)   && \mbox{\qquad by Lemma~\ref{lemma:ChurchAddition_equation}}\\
\n \in \L(\n) && \mbox{\qquad by Lemma~\ref{lemma:L1}}\\
\n \oplus x \in \L(\n) && \mbox{\qquad by Lemma~\ref{lemma:additionmapsloop}}\\
\n \oplus x = \n && \mbox{\qquad by Theorem~\ref{theorem:looponeone}}
\end{eqnarray*}
That completes the proof of the lemma.

\begin{lemma}\label{lemma:smloop}
 Suppose $\N \in \FINITE$ and $\s \k = \s \n$ and $  \k \neq \n$ and
$\k \in \Stem$ and $\n \in \N$.  Suppose $\m \in \N$ with $\n = \k+\m$.
Then $\m \s$ is the identity on $\L(\n)$.
\end{lemma}

\noindent{\em Remark}.    We want to say, ``by the Annihilation Theorem.''
But the domain of successor is more than just the loop, so we must consider
its restriction $f$ to $\L(\n)$, and verify that $f$ satisfies the hypotheses
of the Annihilation Theorem; and after the application, we still have to verify
that the iterates of the restriction are the restrictions of the iterates.  
\smallskip

\noindent{\em Proof}. We have
\begin{eqnarray*}
\s \n = \s (\k + \m)  && \mbox{\qquad since $\n = \k + \m$}\\
 = \k + \s\m  && \mbox{\qquad by Lemma~\ref{lemma:ChurchAddition_equation}}\\
 = \s \k + \m && \mbox{\qquad by Lemma~\ref{lemma:ChurchSuccessorShift}}\\
 = \s \n + \m && \mbox{\qquad since $\s \k = \s \n$}
\end{eqnarray*}
Define $f$ to be the restriction of Church successor  to $\L(\n)$
(which can be done by means of a stratified formula).   
One can verify that $f$ is an injection from $\L(n)$ to $\L(\n)$, in the sense of 
Definition~\ref{definition:ChurchFregenjection}.   The most important step
is that $f$ is one-to-one, by Theorem~\ref{theorem:looponeone}.
We omit the details of the verification (about 180 steps).   
\smallskip

 By Lemma~\ref{lemma:loopfinite}, $\L(\n)$ is finite.  
 By Lemma~\ref{lemma:nplusm} and the hypothesis that $\n = \k + \m$,
 we have $\n = \k + \m$.  
Since $f: \L(\n) \to \L(\n)$ is an injection, we can apply the 
Annihilation Theorem to obtain 
$$ \forall x \in \L(\n)\,(\m f x = x).$$
Then by Lemma~\ref{lemma:successoronloop}, we have 
$$ \forall x \in \L(\n)\,(\m \s x = x).$$
as desired. 
That completes the proof of the lemma.

\begin{theorem} \label{theorem:main} 
 The Church counting axiom implies that $\N$ is not finite.
\end{theorem}

\noindent{\em Proof}. Assume the Church counting axiom, and suppose   $\N$ is finite.
 By Lemma~\ref{lemma:kinstem},
there is a double successor $\s \n = \s \k$ with  
with $\k \neq \n$ and $\k \in \Stem$.   
 Then by Lemma~\ref{lemma:mexists},  there exists a $\preceq$-least
$\m \in \N$ such that $\n = \k + \m$.   Fix that $\m$. Then 
by Theorem~\ref{theorem:orderm}, $\m$ is the order  of successor
restricted to $\L(\n)$.  Explicitly, we have
\begin{eqnarray}
\forall q \in \N\,(q \neq \ChurchZero \imp q\s\n = \n \imp \m \preceq q).
\label{eq:E4523}
\end{eqnarray}
 (This is where we use the Church counting axiom, since Theorem~\ref{theorem:orderm} requires it.) 
\smallskip

We have
\begin{eqnarray*}
\n \neq \ChurchZero  && \mbox{\qquad by Lemma~\ref{lemma:nneqzero}}\\
\n \in \L(\n)        && \mbox{\qquad by Lemma~\ref{lemma:L1}}\\
\n = \s p             && \mbox{\qquad for some $p \in \L(\n)$, by Theorem~\ref{theorem:looponto}} 
\end{eqnarray*}
Now define 
$$X := \L(\n) - \{\n\}.$$
and define $f$
$$ f = ( \{ \langle x,\s x \rangle : x \in X \} - \{\langle p,\n\rangle\}) \cup \{ \ \langle p, \s \n \rangle \}.$$
  Informally,
the idea of the definition of $f$ is that 
  $f(x) = \s x$ except when $x = p$, and $f(p) = \s \n$.
  \smallskip

Our first observation about $f$ is that 
\begin{eqnarray}
\langle x, \s \n \rangle \in f \imp x = p \label{eq:E4545}
\end{eqnarray}
To prove that, suppose $\langle x, \s \n \rangle \in f$.  Then
\begin{eqnarray*}
x \in \L(\n)  && \mbox{\qquad by definition of $X$ and $f$}\\
\s \in   \L(\n) && \mbox{\qquad by Lemma~\ref{lemma:L1}}\\
\s n \neq \n  && \mbox{\qquad by Lemma~\ref{lemma:snneqn}}\\
\s x = \s n \imp x = p && \mbox{\qquad by Theorem~\ref{theorem:looponeone}}
\end{eqnarray*}
Now (\ref{eq:E4545}) follows from the definition of $f$. 
\smallskip

We have
\begin{eqnarray*}
Rel(f) && \mbox{\qquad by 17 omitted steps}\\
f \in \FUNC      && \mbox{\qquad by 60 omitted steps} \\
\s \n \neq \n  && \mbox{\qquad by Lemma~\ref{lemma:snneqn}}\\
dom(f) \subseteq X  && \mbox{\qquad by the preceding line and 30 omitted steps} \\
range(f) \subseteq X && \mbox{\qquad by Theorem~\ref{theorem:looponeone}, Lemma~\ref{lemma:snneqn}, and 46 omitted steps}\\
f: X \to X      && \mbox{\qquad by Theorem~\ref{theorem:looponeone} and 82 omitted steps}
\end{eqnarray*}
I say that $f$ is one-to-one.  Suppose $\langle x, y \rangle \in f$ 
and $\langle u, y\rangle \in f$.  We must prove $x=u$.
 Since $\N$ is finite, $\N$ has decidable equality, by Lemma~\sref{lemma:finitedecidable}.  Therefore
$$ y = \s\n \ \lor \ y \neq \s \n.$$
Case~1, $y = \s\n$.  Then 
\begin{eqnarray*}
x = p && \mbox{\qquad by (\ref{eq:E4545})}\\
u = p && \mbox{\qquad by (\ref{eq:E4545})} \\
x = u && \mbox{\qquad by the preceding two lines}
\end{eqnarray*}
That completes Case~1.
\smallskip

\noindent
Case~2, $y \neq \s \n$.  Then $y = \s x$ and $y = \s u$.  We have
\begin{eqnarray*}
x \in X  \ \land \ u \in X && \mbox{\qquad since $dom(f) \subseteq X$}\\
x \in \L(\n) \ \land \ u \in \L(\n) && \mbox{\qquad since $X \subseteq \L(\n)$}\\
\s x \in \L(\n) && \mbox{\qquad by Lemma~\ref{lemma:L1}}\\
y \in \L(\n)  && \mbox{\qquad since $y = \s x$}\\
x = u && \mbox{\qquad by Theorem~\ref{theorem:looponeone}}
\end{eqnarray*}
That completes Case~2.  That completes the proof that $f$ is one-to-one.
\smallskip

We have now proved that $f$ is an injection, since by definition that 
means $f:X \to X$, $f$ is one-to-one, $f \in \FUNC$ and $Rel(f)$, all of 
which we have verified.  
 Hence we can apply the Annihilation Theorem to $f$ and $X$
to conclude that $\m f$ is the identity on $X$. Explicitly,
\begin{eqnarray}
\forall x \in X\,(\m f x = x)\label{eq:E4598}
\end{eqnarray}
In the rest of the proof, we will show that $\m f$ is not the identity
on $X$, thus contradicting (\ref{eq:E4598}).
\smallskip
  
Let $\alpha:= \s \n$.  Then 
\begin{eqnarray*}
\s\n \neq \n && \mbox{\qquad by Lemma~\ref{lemma:snneqn}}\\
\n \in \L(\n) && \mbox{\qquad by Lemma~\ref{lemma:L1}}\\
\s \n \in L(\n) && \mbox{\qquad by Lemma~\ref{lemma:L1}}\\
\alpha \in \L(\n) && \mbox{\qquad since $\alpha = \s \n$ and $\n \in \L(\n)$}\\
\alpha \neq \n  && \mbox{\qquad since $\alpha = \s \n$  and $\s \n \neq \n$}\\
\alpha \in X && \mbox{\qquad since $X = \L(\n) - \{\n\}$ }
\end{eqnarray*}
I say that  
\begin{eqnarray}
q \neq \n \imp \s q \prec \m \imp  q f \alpha = q \s \alpha. \label{eq:E4424}
 \end{eqnarray}
We prove this by finite induction on $q$.  That is legal, since we 
can stratify that formula, giving $q$ index 6 and $\alpha$ index 0.
$\m$ occurs as a parameter, so we do not need to give it an index, but 
we could give it index 6.  $f$ gets index 3, since it contains pairs of 
objects of index 0; so $q$ contains pairs of objects of index 3; those
pairs have index 5, which is why $q$ gets index 6.  
\smallskip

Base case, $q = \ChurchZero$.  We have $\ChurchZero f \alpha = \alpha = \ChurchZero \s \alpha$, by Lemma~\ref{lemma:zeroAp}.  That completes the base case.
\smallskip

Induction step.  Since we are using finite induction, we get to 
assume 
\begin{eqnarray}
q \neq \n \label{eq:E4724}
\end{eqnarray}
We also assume 
\begin{eqnarray}
\s q \neq \n \ \land \ \s(\s q) \prec \m \label{eq:E4628} 
\end{eqnarray}
 We have to prove
\begin{eqnarray}
\s q \s \alpha = \s q f \alpha \label{eq:E4629}
\end{eqnarray}
 We have  
\begin{eqnarray*}
\s (\s q) \neq \n   && \mbox{\qquad by Lemma~\ref{lemma:precmax2},since $\s(\s q) \prec \m$}\\
\s q \prec \s (\s q)  && \mbox{\qquad by Corollary~\ref{lemma:xprecsx}, since $\s q \neq \n$}\\ 
\s q \prec \m && \mbox{\qquad by Lemma~\ref{lemma:prectrans},
since $\s (\s \q) \prec \m$}\\
q \s \alpha \in \L(\n)  && \mbox{\qquad by Lemma~\ref{lemma:xsmapsloop}, since $\alpha \in \L(\n)$} \\
\L(\n) \subseteq \N && \mbox{\qquad by Lemma~\ref{lemma:L1}}\\
q \s \alpha \in \N  && \mbox{\qquad by the preceding two lines}\\
p \in \N && \mbox{\qquad since $p \in \L(\n)$  and $\L(\n) \subseteq \N$.}
\end{eqnarray*}
Since $\N$ is finite, it has decidable equality, by Lemma~\sref{lemma:finitedecidable}.  Since $p \in \L(\n)$ and $q\s \alpha \in \N$,
we have
$$ \q \s \alpha = p \ \lor\ \q \s \alpha \neq p.$$
We argue by cases accordingly. 
\smallskip

Case~1, $q \s \alpha = p$.  Then 
\begin{eqnarray*}
\s(q \s \alpha) = \s p && \mbox{\qquad since $q \s \alpha = p$}\\
                = \n   && \mbox{\qquad since $\s p = \n$}\\
\s (\s (q \s \alpha)) = \s n && \mbox{\qquad by the previous line}\\
\s (\s (q \s \alpha)) = \alpha && \mbox{\qquad since $\s n = \alpha$}\\
\s (q \s \alpha) = \s q \s \alpha &&  \mbox{\qquad by Theorem~\ref{theorem:successorequation}}\\
\s (\s q \s \alpha) = \alpha && \mbox{\qquad by the preceding two lines}\\
\s(\s q) \s \alpha = \alpha && \mbox{\qquad by Theorem~\ref{theorem:successorequation}} \\
\s(\s q)\s x = x \mbox{ \ for all $x\in \L(\n)$} && \mbox{\qquad by Lemma~\ref{lemma:orderq_helper2}}\\
\m \preceq \s (\s q) && \mbox{\qquad by Theorem~\ref{theorem:orderm} and (\ref{eq:E4523})} \\  
\s (\s q) \prec \m  && \mbox{\qquad by (\ref{eq:E4628})}\\
\s (\s q) \preceq \m \ \land \ \s(\s q) \neq \m && \mbox{\qquad by the definition of $\prec$}\\
\s(\s q) = \m  && \mbox{\qquad by Theorem~\ref{theorem:prectrichotomy2}}
\end{eqnarray*}
But the last two lines are contradictory.  That shows that Case~1 is impossible.
\smallskip

Case~2, $q \s \alpha \neq p$.  We have $\s q \neq \n$  and $\s (\s q) \prec \m$
by hypothesis,
but in order to apply the induction hypothesis, we need $q \neq \n$ and $\s q \prec \m$.  We have $q \neq \n$ by (\ref{eq:E4724}).  Here is
a proof that $\s q \prec \m$:
\begin{eqnarray*}
\s q \prec \s (\s q) && \mbox{\qquad by Lemma~\ref{lemma:preceqsuccessor}, since $\s q \neq \n$}\\
\s (\s q) \prec \m && \mbox{\qquad by hypothesis}\\
\s q \prec \m && \mbox{\qquad by Lemma~\ref{lemma:prectrans}}
\end{eqnarray*}
Now we can use the induction hypothesis. We proceed to the 
proof of the induction step. 
\begin{eqnarray*}
\s q f \alpha = f (q f \alpha) && \mbox{\qquad by Theorem~\ref{theorem:successorequation}}\\
q f \alpha \in X && \mbox{\qquad by Lemma~\ref{lemma:xfmaps}, since $f:X \to X$}\\
\s q f \alpha = f( q \s \alpha)  && \mbox{\qquad by the induction hypothesis} \\
 = \s (q \s \alpha) && \mbox{\qquad since $q \s \alpha \neq p$, by definition of $f$}\\
 = \s q \s \alpha && \mbox{\qquad by Theorem~\ref{theorem:successorequation}}
\end{eqnarray*}
That completes Case~2.  That completes the induction step.
That completes the proof of (\ref{eq:E4629}); that is,
it completes the induction step.  That completes the proof of (\ref{eq:E4424}).
\smallskip

We have $\m \neq \ChurchZero$, since if $\m = \ChurchZero$ then $\n = \k + \m = \k + \ChurchZero = \k$, contradiction.  Then by
 Lemma~\ref{lemma:predecessornotn}, there exists $m_1$  such that
\begin{eqnarray*}
\s m_1 = \m && \\
m_1 \in \N \ \land \ m_1 \neq \n &&
\end{eqnarray*}
(The variable names $m_1$ and $m_2$ in this proof are meant to suggest $m-1$ and $m-2$,
although subtraction has not been defined.) 
I say  $ m_1 \neq \ChurchZero$.  Here is the proof:
\begin{eqnarray*}
 m_1 = \ChurchZero && \mbox{\qquad assumption} \\
  \s m_1 =   \s \ChurchZero && \mbox{\qquad by the previous line}\\
 \s m_1 = \m  && \mbox{\qquad by construction of $m_1$}\\
 \m = \s \ChurchZero && \mbox{\qquad by the preceding lines} \\
 \k \oplus \s \ChurchZero = \n && \mbox{\qquad since $\k + \m = \n$}\\
 \s (\k \oplus \ChurchZero) = \n && \mbox{\qquad by Lemma~\ref{lemma:ChurchAddition_equation}}\\
 \s \k = \n && \mbox{\qquad by Lemma~\ref{lemma:ChurchZero_equation}}\\
 \s (\s \k) = \s \n && \mbox{\qquad by the previous line}\\
 \s (\s \k) = \s \k && \mbox{\qquad since $\s \n = \s \k$}\\
 \s (\s \k) \neq \s \k && \mbox{\qquad by Lemma~\ref{lemma:snneqn}}
 \end{eqnarray*}
That contradiction completes the proof that $m_1 \neq \ChurchZero$.
\smallskip

Then by Lemma~\ref{lemma:predecessornotn},
there exists $m_2$  such that
\begin{eqnarray*}
\s m_2 = m_1 &&\\
m_2 \in \N \ \land \ m_2 \neq \n &&
\end{eqnarray*}
Then 
\begin{eqnarray*}
m_1 \prec \s m_1   && \mbox{\qquad by Lemma~\ref{lemma:xprecsx}, since $m_1 \neq \n$}\\
\s m_2  \prec \m  && \mbox{\qquad since $\s m_2 = \m_1$ and $\s m_1 = \m$}
\end{eqnarray*}  
 Since $\s m_2 \prec \m$ and $m_2 \neq \n$, we have 
\begin{eqnarray}
 m_2  f \alpha = m_2\s \alpha &&\mbox{\qquad by (\ref{eq:E4424})}.\label{eq:4455} 
 \end{eqnarray}
We also have
\begin{eqnarray*}
 \s (\s(p)) = \s \n = \alpha  && \mbox{\qquad by the definitions of $p$ and $\alpha$}\\
 \m \s \alpha = \alpha   && \mbox{\qquad by Lemma~\ref{lemma:smloop}}\\
\s (\s  m_2) \s \alpha = \alpha && \mbox{\qquad since $\s(\s m_2) = \m$} \\
\s (\s  m_2) \s \alpha = \s \n && \mbox{\qquad since $\alpha = \s \n$}\\
\s (\s m_2 \s \alpha) = \s \n && \mbox{\qquad by Theorem~\ref{theorem:successorequation}}\\
(\s  m_2) \s \alpha \in \L(\n) && \mbox{\qquad by Lemma~\ref{lemma:xsmapsloop}}\\
(\s  m_2) \s \alpha =  \n && \mbox{\qquad by Theorem~\ref{theorem:looponeone}}\\
(\s  m_2) \s \alpha = \s p && \mbox{\qquad since $\n = \s p$}\\
\s ( m_2 \s \alpha) = \s p && \mbox{\qquad by Theorem~\ref{theorem:successorequation}}\\
m_s \s \alpha \in \L(\n) && \mbox{\qquad by Lemma~\ref{lemma:xsmapsloop}}
\end{eqnarray*}
\begin{eqnarray}
\hspace{-20pt} m_2 \s \alpha = p &&\mbox{\qquad by Theorem~\ref{theorem:looponeone}} \label{eq:E4459}
\end{eqnarray}
\begin{eqnarray*}
 m_2 f \alpha = p  && \mbox{\qquad by (\ref{eq:4455}) and (\ref{eq:E4459})}\\
 f(m_2 f \alpha) = f(p) && \mbox{\qquad applying $f$ to both sides} \\
 (\s m_2) f(\alpha) = f(p) && \mbox{\qquad by Theorem~\ref{theorem:successorequation}}\\
 m_1  f \alpha = f(p) && \mbox{\qquad since $\s m_2 = m_1$} \\
 f(p) = \alpha && \mbox{\qquad by the definition of $f$}\\
  m_1  f \alpha = \alpha && \mbox{\qquad by the preceding two lines}\\
f( m_1 f \alpha) = f(\alpha) && \mbox{\qquad applying $f$ to both sides} \\
(\s m_1) f \alpha = f(\alpha) && \mbox{\qquad by Theorem~\ref{theorem:successorequation}}\\
  \m f \alpha = f(\alpha) && \mbox{\qquad since $\s m_1 = \m$}\\
  \s (\s \n) \neq \n   && \mbox{\qquad by Lemma~\ref{lemma:ssnneqn}}\\
\alpha = p \imp \s (\s \n) = \n &&\mbox{\qquad since $  \s \n = \alpha$ and $\s p = \n $}\\ 
\alpha \neq p  &&  \mbox{\qquad by the previous two lines} \\
f(\alpha) = \s \alpha &&\mbox{\qquad by definition of $f$, since $\alpha \neq p$} \\
\m f \alpha  = \s \alpha  && \mbox{\qquad since $\m f \alpha = f(\alpha)$} \\
\s \alpha  \neq \alpha  && \mbox{\qquad by Lemma~\ref{lemma:snneqn}} \\
\m f(\alpha) \neq \alpha  && \mbox{\qquad by the previous two lines} \\
\m f(\alpha) = \alpha  && \mbox{\qquad by the Annihilation Theorem}
  \end{eqnarray*}
That contradiction completes the proof of the theorem.

\section{$\N$ not finite implies $\N$ is infinite}

In this section we will show that if $\N$ is not finite, then
$\N$ is infinite, and indeed (what is more)
 Church successor is one-to-one on $\N$.
  Since we proved that the Church counting axiom implies
$\N$ is not finite, it will follow that the Church counting axiom implies
$\N$ is infinite and $\s$ is one-to-one.  
\smallskip

That Church successor is one-to-one means $\s x = \s y \imp x = y$.  
That it is weakly one-to-one means $x \neq y \imp \s x \neq \s y$.
One can check that if successor is weakly one-to-one, then $\N$ has
decidable equality (by induction, with Lemma~\ref{lemma:decidable0} as 
the base case).  With decidable equality, weakly one-to-one implies one-to-one.
\smallskip

The idea of the proof can be explained simply. 
We start at $\ChurchZero$ and make dots on our paper for $0,1,2,\ldots$.
At any moment the set of dots so far written is finite.  If we come to 
a double successor (as shown in Fig.~\ref{figure:rho}),
 then we have a set that contains $\ChurchZero$ and 
is closed under successor, so it is all of $\N$; but then $\N$ is finite,
so that cannot happen.  Instead we continue on indefinitely, i.e., successor
is one-to-one.

 To make that idea rigorous, we will define
 a relation $\JLift$, whose intended interpretation
 is that if $\langle \{x\},y \rangle \in \JLift$,
 then $y$ is the set of dots written down after $x$ steps
of the drawing process described above.  We use $\{x\}$ instead
of $x$ to achieve stratification.   The idea is to 
define $\JLift$ in such a way that $\JLift$ is the least relation such that 
\begin{eqnarray*}
&& \langle \{ \ChurchZero \}, \{\ChurchZero\} \rangle \in \JLift     \\
&& \forall x,y \in \N\,( \langle \{x\}, y \rangle \in \JLift \imp \s x \not\in y \imp \langle \s x, y \cup \{ \s x\} \rangle \in \JLift)   
\end{eqnarray*}
Of course, a proper definition cannot mention $\JLift$ on the right.  Here is 
a proper definition:  

\begin{definition} \label{definition:Jlift}
  $\JLift$ is the set
of all ordered pairs $\langle \{ p \}, q\rangle$ with $p,q \in \N$ such that $\langle \{ p\}, q\rangle$ belongs
  to every set $w$
satisfying the following conditions:
\begin{eqnarray*}
&& \langle \{ \ChurchZero \}, \{\ChurchZero\} \rangle \in w     \\
&& \forall x,y \in \N\,( \langle \{x\}, y \rangle \in w \imp \s x \not\in y \imp \langle \{ \s x\}, y \cup \{ \s x\} \rangle \in w)   
\end{eqnarray*}
\end{definition}
The formula is stratified, giving $x$ index 0, $y$ index 1, so $\langle \{x\},y\rangle$
gets index 3; then $2$ gets index 4. $\ChurchZero$ is a parameter, so does not 
need an  index, but we could give it index 0.  Either way, the formula is 
stratified, so the definition can be given in INF.

\begin{lemma} \label{lemma:JLift0} 
 $\langle \{\ChurchZero\}, \{ \ChurchZero \}\rangle \in \JLift$.
\end{lemma}

\noindent{\em Proof}. Immediate from the definition of $\JLift$.

\begin{lemma} \label{lemma:JLiftRec}
 
$$ \langle \{x\}, y \rangle \in \JLift \imp \s x \neq y \imp 
\langle \{\s x\}, y \cup \{ \s x\} \in \JLift.$$
\end{lemma}

\noindent{\em Proof}. Follows  from Definition~\ref{definition:Jlift}
in about 25 steps (omitted here).

\begin{lemma} \label{lemma:JLiftfinite}
 Suppose $\langle \{ x\}, y\rangle \in \JLift$.
Then $y \in \FINITE$ and $y \subseteq \N$ and $x \in \N$.
\end{lemma}

\noindent{\em Proof}.  Let $W$ be the set of members of $\JLift$ satisfying
the conditions in the lemma; explicitly,
$$W = \{\ \langle \{ x\}, y\rangle \in \JLift: y \in \FINITE \ \land \ y \subseteq \N \ \land \ x \in \N\ \}.$$  
Then $W$ satisfies the closure 
conditions in  Definition~\ref{definition:Jlift}:
\begin{eqnarray*}
\langle \{\ChurchZero\}, \{\ChurchZero\}\rangle \in W && \mbox{\qquad by  
Lemma~\ref{lemma:JLift0}}\\
\{\ChurchZero\} \in \FINITE && \mbox{\qquad by Lemma~\sref{lemma:singletons_finite}}\\
y \cup \{ \s x \} \in \FINITE && \mbox{\qquad if $\s x \not\in y$, by 
Lemma~\sref{lemma:finite_adjoin}}
\end{eqnarray*}
The details, omitted here, take about 90 steps. 
Therefore $Z \subseteq W$.
That completes the proof of the lemma.

\begin{lemma} \label{lemma:JLiftMaps_helper}
Suppose $\langle \{x\}, y \rangle \in \JLift$. Then 
\begin{eqnarray*}
&&\ChurchZero \in y \\
&& x \in y \\
&&\forall u\,(u \in y \imp u \neq x \imp \s x \in y)
\end{eqnarray*}
\end{lemma}

\noindent{\em Remark}.  The last condition, expressed 
in words, is ``$y$ is closed under successor except $x$.''
\smallskip

\noindent{\em Proof}.  Let $W$ be defined as the 
set of all $\langle \{x\},y \rangle \in \JLift $ such that   
conditions of the lemma are satisfied.  Since the formulas in the lemma
are stratified, $W$ can be defined in INF.
 We will prove $W$ satisfies
the closure conditions in the definition of $\JLift$.

First, $\langle \{\ChurchZero\}, \{\ChurchZero\}\rangle \in W$; it
belongs to $Z$ by Lemma~\ref{lemma:JLift0}, and the other conditions 
are straightforward.

Second, assume $\langle \{x\}, y\rangle \in W$ and $\s x \not \in y$.
We must show $\langle \{\s x\}, y \cup \{ \s x\} \rangle \in W$.
By Lemma~\ref{lemma:JLiftRec}, it belongs to $\JLift$.
 
Since $\langle \{x\}, y\rangle \in W$, we have $\ChurchZero \in y$.
Hence $\ChurchZero \in y \cup \{ \s x\}$.

By Lemma~\ref{lemma:JLiftfinite}, we have $y \in \FINITE$. 
Then  $y \in \DECIDABLE$, by Lemma~\sref{lemma:finitedecidable}.  
We have to show $y \cup \{\s x\}$
is closed under successor except $\s x$. Let $u \in y \cup \{ \s x\}$
with $u \neq \s x$.  Then $u \in y$.  Since $y \in \DECIDABLE$,
we have $u = x \ \lor\  u \neq x$.  If $u = x$ then $\s u = \s x \in y \cup \{ \s x\}$.
If $ u \neq x$ then $\s x \in y$ since $\langle \{x\}, y\rangle \in W$;
therefore $\s x \in y \cup \{\s x\}$ as well.
\smallskip

That completes the proof that $W$ satisfies the closure conditions.
That completes the proof of the lemma.

\begin{lemma}[No loops]\label{lemma:noloops}
Assume $\N$ is not finite.  Suppose $\langle \{x\}, y\rangle \in \JLift$.
Then $\s x \not \in y$.
\end{lemma}

\noindent{\em Proof}.  Suppose $\s x \in y$.  Then 
\begin{eqnarray*}
\ChurchZero \in y  && \mbox{\qquad by Lemma~\ref{lemma:JLiftMaps_helper}}\\
u \in y \imp u \neq x \imp \s u \in y &&\mbox{\qquad by Lemma~\ref{lemma:JLiftMaps_helper}}\\
\s x \in y  && \mbox{\qquad by hypothesis}\\
y \in \FINITE && \mbox{\qquad by Lemma~\ref{lemma:JLiftfinite}}\\
y \in \DECIDABLE && \mbox{\qquad by Lemma~\sref{lemma:finitedecidable}}\\
x \in y          && \mbox{\qquad by Lemma~\ref{lemma:JLiftMaps_helper}} \\
u \in y \imp u = x \ \lor \ u \neq x  &&\mbox{\qquad by the preceding lines}\\
u \in y \imp \s u \in y  && \mbox{\qquad by the preceding lines}\\
\N  \subseteq y          && \mbox{\qquad by the definition of $\N$}\\
y \subseteq \N           && \mbox{\qquad by Lemma~\ref{lemma:JLiftfinite}}\\
y = \N                    && \mbox{\qquad by the preceding two lines}\\
\N \in \FINITE           && \mbox{\qquad since $y \in \FINITE$}
\end{eqnarray*}
But that contradicts the hypothesis.  That completes the  proof of the lemma.

\begin{lemma} \label{lemma:JLiftMaps}
Assume $\N$ is not finite.
Then 
$$\forall x \in \N\, \exists y\,(\langle \{x\},y\rangle \in \JLift).$$
\end{lemma}

\noindent{\em Proof}.
The formula in the lemma is stratified, giving $x$ and $y$ index 0, since $Z$ is 
a definable relation (occurring here as a parameter).  Therefore we
may prove it by induction on $x$. 
\smallskip

Base case, $x = \ChurchZero$, holds by Lemma~\ref{lemma:JLift0}.
\smallskip

Induction step.  Suppose $\langle \{x\},y \rangle \in \JLift$.
By Lemma~\ref{lemma:noloops}, $\s x \not \in y$.
Then by Lemma~\ref{lemma:JLiftRec}, $\{\s x\}, y \cup \{\s x\} \in \JLift$.
That completes the induction step.  That completes the proof of the lemma.

\begin{lemma}\label{lemma:JLiftdom} Suppose $\langle \{   x \}, p \rangle \in \JLift$.   Then $x = \ChurchZero$,  or $x = \s u$ for some $u \in \N$.
\end{lemma}

\noindent{\em Proof}.  Define 
$$Z := \{ \ \langle \{ x\}, p\rangle \in \JLift : x = \ChurchZero \ \lor \ 
\exists u \in \N\, (\s u = x) \}.$$
Then $Z$ satisfies the conditions in the definition of $\JLift$,
as one verifies in about 70 steps (here omitted).  
Therefore $\JLift \subseteq Z$.  To finish the proof:
\begin{eqnarray*}
\langle \{   x \}, p \rangle \in \JLift  && \mbox{\qquad assumption}\\
 \langle \{   x \}, p \rangle \in Z && \mbox{\qquad since $\JLift \subseteq Z$}\\
 x = \ChurchZero \ \lor \ 
\exists u \in \N\, (\s u = x) && \mbox{\qquad by the definition of $Z$} 
\end{eqnarray*}
The two resulting cases are just the conditions that $Z$ has been proved
to satisfy.  That completes the proof of the lemma.

\begin{lemma} \label{lemma:JLift1} Suppose $\langle \{ \ChurchZero \}, p \rangle \in \JLift$.   Then $p = \{ \ChurchZero\}$.
\end{lemma}

\noindent{\em Proof}.  Suppose $p \neq \{ \ChurchZero\}$.  Then define
$$Z:= \JLift - \{ \langle \{ \ChurchZero \}, p \rangle \}.$$
One can verify that $Z$ satisfies the closure
conditions in the definition of $\JLift$.  (It takes about 60 steps,
omitted here, using several of the lemmas above, including Lemmas~\ref{lemma:JLiftRec}
and \ref{lemma:JLiftdom}.) 
Therefore $\JLift \subseteq  Z$.
But that is a contradiction.  Therefore $ \neg\neg\, p = \{\ChurchZero\}$.
Now
\begin{eqnarray*}
p \in \FINITE && \mbox{\qquad by Lemma~\ref{lemma:JLiftfinite}}\\
p \in \DECIDABLE && \mbox{\qquad by Lemma~\sref{lemma:finitedecidable}}\\
\ChurchZero \in p && \mbox{\qquad by Lemma~\ref{lemma:JLiftMaps_helper}}\\
 \neg\neg\, \forall u \in p\,  
( u = \ChurchZero)   && \mbox{\qquad since $ \neg\neg\, p = \{\ChurchZero\}$}\\  
 \forall u \in p\, (\neg\neg\,
( u = \ChurchZero))  && \mbox{\qquad by intuitionistic logic}\\
\forall u \in p\,  
( u = \ChurchZero)  && \mbox{\qquad since $p \in \DECIDABLE$}\\
p = \{ \ChurchZero\} && \mbox{\qquad by the preceding line and $\ChurchZero \in p$}
\end{eqnarray*}
That completes the proof of the lemma.

\begin{lemma} \label{lemma:JLift2} Suppose $x \in \N$ and 
$\langle \{ \s x\}, y \rangle \in \JLift$.
Then there exist $u$ and $p$ such that 
\begin{eqnarray*}
&& \s x = \s u  \\
&& \langle \{ u\}, p \rangle \in \JLift \\
&& \s u \not \in p     \\
&& y = p \cup \{ \s u\}
\end{eqnarray*}
\end{lemma}

\noindent{\em Remark}.  The point of the lemma (and the preceding one)
is that everything
in $\JLift$ is in $\JLift$ because it has been constructed according to 
the two construction rules in the definition.
\medskip

\noindent{\em Proof}.  Define 
\begin{eqnarray*}
Z:= \{\ z \in \JLift : && z = \langle \{\ChurchZero\},\{\ChurchZero\}\rangle \ \lor \\
 && (\exists x,y\,(x \in \N \ \land \ z = \langle \{\s x\},y\rangle) \ \land \\
 && (\forall x,y\,(x \in \N \imp  z = \langle \{\s x\},y\rangle) \imp \\
&& \exists u,p\,( \langle \{ u\}, p\rangle \in \JLift \ \land \ x \in \N \ \land \ u \in \N \ \land \\
&& \s u = \s x \ \land \ \s u \not \in p \ \land \ y = p \cup \{\s u\}) \ \}
\end{eqnarray*}
The formula is stratified, giving $x$ and $u$ index 0, and $y$ and $p$ index 1.
Then the ordered pairs are pairs of type 1 objects, so they get type 3.
So $z$ gets index 3, and $\JLift$ is a parameter.  Therefore the definition
is legal in INF. 
 
Then one can verify that $Z$ satisfies the closure conditions in 
the definition of $\JLift$.  (It takes about 110 steps, omitted here.)
There are several variations of the definition of $Z$ that 
look equally convincing but are in fact not correct.  Once the definition
is correct, the 110 steps mentioned are fairly straightforward.  
\smallskip 
 
Having derived that $Z$ satisfies the closure conditions, we have 
$\JLift \subseteq Z$, by definition of $\JLift$.  Now suppose $x \in \N$
and $\langle \{\s x\}, y \rangle \in \JLift$.  Then since $\JLift \subseteq Z$
we have $\langle \{\s x\}, y \rangle \in Z$.  Substituting 
$\langle \{\s x\}, y \rangle$ for $z$ in the definition of $Z$, 
the disjunction on the right of the definition gives rise to two cases.
\smallskip

Case~1, $\langle \{\s x\}, y \rangle   = 
\langle \{\ChurchZero\},\{\ChurchZero\}\rangle$.
Then
\begin{eqnarray*}
\{ \s x\} = \{ \ChurchZero\}  && \mbox{\qquad by Lemma~\sref{lemma:ordered_pair_equality}}\\
\s x = \ChurchZero  && \mbox{\qquad by Lemma~\sref{lemma:single_oneone}}\\
x \in \N           && \mbox{\qquad by hypothesis}\\
\s x \neq \ChurchZero && \mbox{\qquad by Theorem~\ref{theorem:successoromitszero}}
\end{eqnarray*}
That disposes of Case~1.   (Note the necessity of including $x \in \N$ as a 
hypothesis of the lemma; we cannot rule out the strange possibility that 
$\s x$ might be $\ChurchZero$ for some $x$ that is not a Church number.)
\smallskip

Case~2, the other disjunction of the definition of $Z$ holds with 
$z:= \langle \{\s x\}, y \rangle$.  Then it is a straightforward 
ten steps (which we omit here) to deduce the conclusion of the lemma.
That completes the proof of the lemma.

\begin{definition} \label{definition:comparable}
Two subsets $x$ and $y$ of $\N$ are {\bf comparable}:
$$  \comp(x,y) :=  x \subseteq \N \ \land \ y \subseteq \N  \ \land \ \neg\neg\, (x \subseteq y  \ \lor \ y \subseteq x).$$
\end{definition}

\noindent{\em Remark}.  Comparability is symmetric and reflexive, but not 
 transitive.  Perhaps this would have worked without the double negation,
 but it certainly does work with the double negation. 
\smallskip

\begin{lemma}\label{lemma:JFUNC_helper2}
Suppose $\langle \{x\}, y \rangle \in \JLift$ and 
$y$ is comparable to every element in the rangle of $\JLift$ and 
$$\forall z\,\langle \{x\}, z \rangle \in \JLift \imp y = z.$$
Then $y \cup \{\s x\}$ is comparable to every element in the 
range of $\JLift$.
\end{lemma}

\noindent{\em Remark}. Functionality at $\{x\}$ and comparability at $y$
imply comparability at $y \cup \{\s x\}$. 
\medskip

\noindent{\em Proof}.  Suppose
 $\langle \{x\}, y\rangle \in \JLift$
 and $\langle \{ u\}, p\rangle \in \JLift$.
Then $ \comp (y,p)$, by the first hypothesis.  We must show
$ \comp (y \cup \{ \s x\}, p)$.  It suffices to prove it from 
$$ x = u \ \lor  x \neq u \imp y \subseteq p \ \lor \ p \subseteq y \imp  \comp (y \cup \{ \s x\}, p),$$
since double-negating that statement yields the desired
$$ \comp(y,p) \imp \comp (y \cup \{ \s x\}, p).$$
Therefore we may assume 
\begin{eqnarray}
x = u \ \lor  x \neq u  && \mbox{\qquad assumption} \label{eq:E5477}\\
 y \subseteq p \ \lor \ p \subseteq y && \mbox{\qquad assumption}\nonumber
\end{eqnarray}
We argue by cases accordingly to prove $\comp (y \cup \{ \s x\}, p)$.
\smallskip

Case~1, $p \subseteq y$.  Then $p \subseteq y \cup \{ \s x\}$, done.
\smallskip

Case~2. $y \subseteq p$. Then 
\begin{eqnarray}
x \in y  && \mbox{\qquad by Lemma~\ref{lemma:JLiftMaps_helper}}\nonumber\\
x \in p  && \mbox{\qquad since $y \subseteq p$\nonumber}\\
x \neq u \imp \s x \in p && \mbox{\qquad by Lemma~\ref{lemma:JLiftMaps_helper}}\nonumber\\
x \neq u \imp y \cup \{\s x\} \subseteq p && \mbox{\qquad since $y \subseteq p$} \label{eq:E5491}\\
x = u \ \lor  x \neq u && \mbox{\qquad by (\ref{eq:E5477})} \nonumber
\end{eqnarray}
We argue by cases accordingly. 
\smallskip

Case~1, $x = u$.  
  Then $\langle \{x\}, y \rangle \in \JLift$ and 
$\langle \{x\}, p\rangle \in \JLift$.  So by hypothesis $y = p$.  Then 
$p \subseteq y \cup \{ \s x\}$, so $\comp(y,p)$.
\smallskip

Case~2, $x \neq u$.  Then 
\begin{eqnarray*}
y \cup \{\s x\} \subseteq p && \mbox{\qquad by (\ref{eq:E5491})}\\
\comp(y \cup \{\s x\}, p    && \mbox{\qquad by the definition of $\comp$ and logic}
\end{eqnarray*}
That completes the proof of the lemma.

\begin{lemma}\label{lemma:JFUNC_helper3}
Assume $\N \not\in \FINITE$.
Suppose $x \in \N$ and    
\begin{eqnarray}
&&\langle \{x\}, y \rangle \in \JLift  \nonumber\\
&&\forall t,u\,(\langle\{t\},u\rangle \in \JLift \imp \comp(y,u) \label{eq:comparability} \\
&& \forall z\,\langle \{x\}, z \rangle \in \JLift \imp y = z  \label{eq:functionality}\\
&& \forall t,u\, (\langle \{t\},u\rangle \in \JLift \imp t = x \ \lor \ t \neq x) 
                                                   \label{eq:domaindecidability} \\
&& \langle \{\s x\},u\rangle \in \JLift   \nonumber \\
&&  \langle \{\s x\}, v \rangle \in \JLift. \label{eq:E5521}
  \end{eqnarray}
 Then $u=v$.
\end{lemma}

\noindent{\em Remark}. It may help to attach names to the formulas.

(\ref{eq:comparability}) is ``comparability''. 

(\ref{eq:functionality}) is ``functionality''.

(\ref{eq:domaindecidability}) is ``domain decidability''.

\noindent
Then the lemma says:
Functionality at $\{x\}$ and domain decidability at $x$ 
 and comparability at $y$
imply functionality at $\{ \s x \}$. 
\medskip

\noindent{\em Proof}.  
\begin{eqnarray*}
u = p \cup  \{\s t \}  \ \land \ \s t = \s x \ \land \ \langle \{t\}, p \rangle \in \JLift  && \mbox{\qquad by Lemma~\ref{lemma:JLift2}}\\
v = q \cup \{\s r \}  \ \land \ \s r = \s x \ \land \ \langle \{r\}, q \rangle \in \JLift  && \mbox{\qquad by Lemma~\ref{lemma:JLift2}}\\
\langle \{\s x\}, y \cup \{\s x\} \rangle \in \JLift && \mbox{\qquad by Lemma~\ref{lemma:JLiftRec}}\\
\langle \{x\}, y \rangle \in \JLift && \mbox{\qquad by hypothesis}\\
\s x \not\in y && \mbox{\qquad by Lemma~\ref{lemma:noloops}, since $\N \not\in \FINITE$}\\
\comp(y,p)  && \mbox{\qquad by hypothesis, since $\langle \{t\}, p \rangle \in \JLift $ }
\end{eqnarray*}
By Lemma~\ref{lemma:JLiftMaps_helper}, $y$ is closed under successor except $x$,
and $p$ is closed under successor except $t$, and 
$q$ is closed under successor except $r$.  Explicitly,
\begin{eqnarray}
\forall z \in y\, (z \neq x \imp \s z \in y) \label{eq:E5339}  \\
\forall z \in p\, (z \neq t \imp \s z \in p) \label{eq:E5340} \\
\forall z \in q\, (z \neq r \imp \s z \in q) \label{eq:E5341} 
\end{eqnarray}
\smallskip

Now I say that 
\begin{eqnarray}
 t = x  \label{eq:E5538}
\end{eqnarray}
By (\ref{eq:domaindecidability}), we have $t \in x \ \lor \ t \neq x$.
In case $t = x$, we have (\ref{eq:E5538}) immediately; so we may 
assume 
\begin{eqnarray}
  t \neq x  \label{eq:E5547}
\end{eqnarray}
We must derive a contradiction.
I say that $t \not \in y$:
\begin{eqnarray*}
t \in y  && \mbox{\qquad assumption, for contradiction}\\
\s t \in y && \mbox{\qquad by (\ref{eq:E5339}) and (\ref{eq:E5547})} \\
\s x \in y &&\mbox{\qquad since $\s x = \s t$}\\
\s x \not \in y && \mbox{\qquad as shown above} 
\end{eqnarray*}
Therefore $t \not\in y$, as claimed.
Similarly $x \not \in p$. 
We have
\begin{eqnarray*}
t \in p  && \mbox{\qquad by Lemma~\ref{lemma:JLiftMaps_helper}} \\
\neg\,(p \subseteq y) && \mbox{\qquad since $t \not \in y$}\\
x \in y  && \mbox{\qquad by Lemma~\ref{lemma:JLiftMaps_helper}} \\
\neg\,(y \subseteq p) && \mbox{\qquad since $x \not \in p$} \\
\neg\, (p \subseteq y \ \lor \ y \subseteq p)  && \mbox{\qquad by logic}\\
\neg\, (y \subseteq p \ \lor \ p \subseteq y)  && \mbox{\qquad by logic}\\
\neg\, \comp(y,p) && \mbox{\qquad by the definition of $\comp$} 
\end{eqnarray*}
But we have derived $\comp(y,p)$ above. That
contradiction completes the proof of (\ref{eq:E5538}), namely $ t = x$.
\smallskip
 
Proceeding, we have 
\begin{eqnarray*}
\langle \{x\},y \rangle \in \JLift && \mbox{\qquad by hypothesis}\\
\langle \{t\},p \rangle \in \JLift && \mbox{\qquad derived above} \\
\langle \{x\},p \rangle \in \JLift && \mbox{\qquad since $t = x$}\\
 y = p  && \mbox{\qquad by the functionality hypothesis}
\end{eqnarray*}
Interchanging $r$ for $t$ and $q$ for $p$,  and using (\ref{eq:E5341})
instead of (\ref{eq:E5340}), we similarly derive
$r=x$ and $y = q$.  Then $t = r$, since both are equal to $x$,
and $p = q$, since both are equal to $y$. 
Then $u = v$, since
$$ u = p \cup \{ \s t\} = q \cup \{ \s r \} = v.$$
That completes the proof of the lemma.

\begin{lemma}\label{lemma:JFUNC_helper4}
Assume $\N \not\in \FINITE$.
Suppose $x \in \N$ and  $\langle \{x\}, y \rangle \in \JLift$ and
\begin{eqnarray*}
&&\forall t,u\,(\langle\{t\},u\rangle \in \JLift \imp \comp(y,u))  \\
&& \forall z\,\langle \{x\}, z \rangle \in \JLift \imp y = z  \\
&& \forall t,u\, (\langle \{t\},u\rangle \in \JLift \imp t = x \ \lor \ t \neq x) 
                                                   \label{eq:domaindecidability2} \\
&& \langle \{\s x\},u\rangle \in \JLift  
  \end{eqnarray*}
 Then 
$$ \forall t,q\,( \langle\{t\},q \rangle \in \JLift \imp \s x = t \lor \s x \neq  t).$$
\end{lemma}

\noindent{\em Remark}.  This lemma adds to Lemma~\ref{lemma:JFUNC_helper3}
by extending ``domain decidability'' from $x$ to $\s x$.
\medskip

\noindent{\em Proof}. 
\begin{eqnarray*}
\langle\{t\},q \rangle \in \JLift  && \mbox{\qquad assumption} \\ 
\langle \{x\}, y \rangle \in \JLift &&\mbox{\qquad hypothesis} \\
\s x \not \in y && \mbox{\qquad by Lemma~\ref{lemma:noloops}}\\
\langle \{\s x\}, y \cup \{ \s x\} \rangle \in \JLift &&\mbox{\qquad by Lemma~\ref{lemma:JLiftRec}}\\
t \in \N  && \mbox{\qquad by Lemma~\ref{lemma:JLiftfinite}}\\
x \in \N  && \mbox{\qquad by Lemma~\ref{lemma:JLiftfinite}} \\
\langle \{\s x\}, z  \rangle \in \JLift \imp z = y \cup\{\s x\} &&
    \mbox{\qquad by Lemma~\ref{lemma:JFUNC_helper3}}
\end{eqnarray*}

By Lemma~\ref{lemma:decidable0}, $t = \ChurchZero \ \lor\ 
t \neq \ChurchZero$.  We argue by cases accordingly.
\smallskip

Case~1, $t= \ChurchZero$.  Then  $\s x \neq t$, by 
Theorem~\ref{theorem:successoromitszero}.  That completes Case~1.
\smallskip

Case~2, $t \neq \ChurchZero$.  Then 
\begin{eqnarray*}
t = \s m && \mbox{\qquad for some $m \in \N$, by Lemma~\ref{lemma:predecessor}}\\
\langle\{\s m\},q \rangle \in \JLift&& \mbox{\qquad since 
$\langle\{t\},q \rangle \in \JLift$}\\
m = x \ \lor \ m \neq x && \mbox{\qquad by (\ref{eq:domaindecidability2})}
\end{eqnarray*}

Case~2a, $m = x$. Then  $\s x = \s m$, so we are done.
\smallskip

Case~2b, $m \neq x$.  I say that $\s x \neq \s m$.  Suppose
$\s x = \s m$. Then
\begin{eqnarray*}
\langle\{\s m\},q \rangle \in \JLift && \\
\langle\{\s x\},q \rangle \in \JLift && \mbox{\qquad since $\s x = \s m$}\\
\langle \{\s x\}, y \cup \{ \s x\}\rangle \in \JLift &&\\
q = y \cup \{\s x\} && \mbox{\qquad by Lemma~\ref{lemma:JFUNC_helper3}}\\
\langle\{x\},y \rangle\in \JLift && \mbox{\qquad by hypothesis}\\
\langle\{m\}, u \rangle \in \JLift && \mbox{\qquad for some $u \in\N$,
                 by Lemma~\ref{lemma:JLiftMaps}} \\
m \in u    && \mbox{\qquad by Lemma~\ref{lemma:JLiftMaps_helper}}\\
\s m \not \in u && \mbox{\qquad by Lemma~\ref{lemma:noloops}}\\
\langle\{\s m\}, u \cup \{\s m\} \rangle\in \JLift && \mbox{\qquad by Lemma~\ref{lemma:JLiftRec}}\\
\langle\{\s x\}, u \cup \{\s x\} \rangle\in \JLift && \mbox{\qquad  since $\s x = \s m$}\\
u \cup \{\s x \} = y \cup \{\s x\} && \mbox{\qquad by Lemma~\ref{lemma:JFUNC_helper3}}\\
\s x \not \in y && \mbox{\qquad by Lemma~\ref{lemma:noloops}, since $\langle\{x\},y \rangle\in \JLift$}\\
\s m \not \in u && \mbox{\qquad by Lemma~\ref{lemma:noloops}, since $\langle\{m\}, u \rangle \in \JLift$}\\
\s x \not \in u  && \mbox{\qquad since $\s x = \s m$}\\
u = y  && \mbox{\qquad since $u \cup \{\s x \} = y \cup \{\s x\} $ }\\
m \in y && \mbox{\qquad since $m \in u$ and $u = y$}\\
\forall q \in y\,(q \neq x \imp \s q \in y) && \mbox{\qquad by Lemma~\ref{lemma:JLiftMaps_helper}, since $\langle\{x\},y\rangle \in \JLift$}\\
\s m \in y  && \mbox{\qquad since $m \neq x$ and $m \in y$}\\
\s x \in y  && \mbox{\qquad since $\s x = \s m$} \\
\s x \not \in y && \mbox{\qquad as proved above, by  Lemma~\ref{lemma:noloops}}
\end{eqnarray*}
That contradiction completes the proof that $\s x \neq \s m$.
That completes Case~2b.  That completes Case~2.  That completes
the proof of the lemma.
\medskip

Now we are in a position to prove $\JLift\in \FUNC$; that is,
the value $y$ such that $\langle \{\s x\},y \rangle \in \JLift$
is uniquely determined by $\s x$.  We prove this property 
simultaneously with the property that equality is decidable between
$x$ and any element of the domain of $\JLift$, and $y$ is comparable
to any element of the range of $\JLift$.  The last three lemmas
together have the information needed to carry out the induction step.%
\footnote{If any readers think  this proof is too complicated, I can 
assure them there are several simpler ``proofs'' that are not correct.
This one may be complicated, but it {\em is} correct. It is, however,
annoying that the picture is so much simpler than the proof.  At least
a part of the problem is that we do not have decidability of equality on 
$\N$ at this point.
}

 \begin{lemma}\label{lemma:JLiftFUNC}  
Assume $\N \not\in \FINITE$.  Then $\JLift$ is a functional relation,
in the sense that if $\langle \{ x\}, y \rangle \in \JLift$
and $\langle\{ x\}, z\rangle \in \JLift$, then $y = z$.
\end{lemma}

\noindent{\em Proof}.  As described above,  we actually prove a more
complicated proposition.  Namely, the conjunction of these three:
\begin{eqnarray*}
\mbox{functionality}&& \forall y,z\, \langle \{ x\}, y \rangle \in \JLift \imp 
 \langle\{ x\}, z\rangle \in \JLift \imp y = z  \\
\mbox{domain decidability}&& \forall y,p,t\, \langle \{x\},y \rangle \in \JLift \imp 
 \langle\{ t\}, p\rangle \in \JLift \imp x = t \ \lor \ x \neq t  \\
\mbox{comparability} &&  \forall y,p,t\, \langle \{x\},y \rangle \in \JLift \imp 
 \langle\{ t\}, p\rangle \in \JLift \imp \comp(y,p) 
\end{eqnarray*}
These formulas are all stratified, giving $x$ and $t$ index 0 and $y$, $z$, and $p$
index 1.  Therefore we may proceed by induction.
\smallskip

Base case, $x = \ChurchZero$.  By Lemma~\ref{lemma:JLift1}, 
we have $\langle \{ \ChurchZero \}, y \rangle \in \JLift$ if 
and only if $y = \{ \ChurchZero \}$.  That takes care of 
functionality when $x = \ChurchZero$. By Lemma~\ref{lemma:decidable0},
we have $t \in \N \imp \ChurchZero = t \ \lor \ \ChurchZero \neq t$.
That takes care of domain decidability when $x = \ChurchZero$.
To prove comparability, it suffices by Lemma~\ref{lemma:JLift1}
to prove that $\{ \ChurchZero\}$ is comparable to any $y$ such 
that $\langle\{ x\}, y\rangle \in \JLift$ for some $x,y$.  But 
any such $y$ contains $\ChurchZero$, by Lemma~\ref{lemma:JLiftMaps_helper}.
Therefore $\{ \ChurchZero\} \subseteq y$. Therefore $\comp(\{\ChurchZero\},y)$.
That completes the base case. 
\smallskip

Induction step.  We will use the three
lemmas \ref{lemma:JFUNC_helper2}, \ref{lemma:JFUNC_helper3},
and \ref{lemma:JFUNC_helper4} to carry out the induction step.
We will spell out the logic explictly here.  To that end,
let $A$, $B$, and $C$ be the three sets, respectively, of $x$
satisfying functionality, domain decidability, and comparability,
as explicitly written about above.  Let $Z:= A \cap B \cap C$.
 Then the induction hypotheses
is $x \in X$.  We have
\begin{eqnarray*}
x \in Z \imp \s x \in A &&\mbox{\qquad by Lemma~\ref{lemma:JFUNC_helper3}}\\
x \in Z \imp \s x \in B &&\mbox{\qquad by Lemma~\ref{lemma:JFUNC_helper4}}\\
x \in Z \imp \s x \in A  \imp \s x \in C &&\mbox{\qquad by Lemma~\ref{lemma:JFUNC_helper2}}
\end{eqnarray*}
Combining these three implications, we have $x \in Z \imp \s x \in Z$.
Note that the induction step for $A$ is used again in proving
the induction step for $C$.  Let us look at that part of the argument,
i.e., at the third implication listed above. 
\begin{eqnarray*}
\langle \{\s x\},y \rangle \in \JLift && \mbox{\qquad assumption}\\
\langle \{ t\}, p \rangle \in \JLift  && \mbox{\qquad assumption} \\
\langle \{ x\}, u \rangle \in \JLift  && \mbox{\qquad for some $u$, by Lemma~\ref{lemma:JLiftMaps}}\\
\langle \{\s x\}, u \cup \{ \s x\}\rangle && \mbox{\qquad by Lemma~\ref{lemma:JLiftRec}}
\end{eqnarray*}
and we have to prove $\comp(y,p)$. We have
\begin{eqnarray*}
\s x \not \in u && \mbox{\qquad by Lemma~\ref{lemma:noloops}} 
\end{eqnarray*}
Then the crucial step is
\begin{eqnarray*}
y = u \cup \{\s x\} 
\end{eqnarray*}
which comes from the induction step for functionality, i.e. $\s x \in A$.
Then we have
\begin{eqnarray*}
\forall z\,(\langle \{x\},z\rangle \in \JLift \imp u = z) && \mbox{\qquad since $x \in A$}\\
\forall r,p\,(\langle \{r\},p \rangle \in \JLift \imp \comp(u,p) && \mbox{\qquad since $x \in C$} \\
\comp ( u \cup \{\s x\},p) && \mbox{\qquad by Lemma~\ref{lemma:JFUNC_helper2}} \\
\comp (y,p)  && \mbox{\qquad since $y = u \cup \{\s x\}$}
\end{eqnarray*}
That is the desired goal of the third implication.
\smallskip

 The first two implications are 
straightforward applications of the lemmas; they take about 
100 ``bookkeeping'' steps, which we omit here.
 That completes the proof of the lemma.
\smallskip

\noindent{\em Remark}.  The original plan was to prove that 
  $\JLift$ is the graph of a function, and then introduce
 the function itself by a comprehension term, so that 
 \begin{eqnarray*}
 \JC(\{\ChurchZero\}) = \{\ChurchZero\} \\
  \JC(\{\s x\}) = \JC(x) \cup \{ \s x\}  
\end{eqnarray*}
That is now easily done; but we do not do it, because from this 
point we can reach the main theorem directly from $\JLift$.

\begin{theorem}\label{theorem:Churchsuccessorweaklyoneone}
Suppose $\N$ is not finite.  Then Church successor is weakly one-to-one on $\N$,
in the sense that $$ x \neq y \imp \s x \neq \s y,$$
and $\N$ is (therefore) infinite.  
\end{theorem}

\noindent{\em Proof}.  Suppose $x \in \N$ and $t \in \N$ and $\s x = \s t$.
and $x \neq t$.  We must derive a contradiction.  We have
\begin{eqnarray*}
\langle \{ x\}, y \rangle \in \JLift && \mbox{\qquad for some $y$, by Lemma~\ref{lemma:JLiftMaps}}\\
\s x \not \in y && \mbox{\qquad by Lemma~\ref{lemma:noloops}} \\
\langle \{ \s x\}, y \cup \{\s x\} \rangle \in \JLift
                     && \mbox{\qquad by Lemma~\ref{lemma:JLiftRec}}\\
\langle \{ t\}, p \rangle \in \JLift && \mbox{\qquad for some $p$, by Lemma~\ref{lemma:JLiftMaps}}\\
\s t \not \in p && \mbox{\qquad by Lemma~\ref{lemma:noloops}} \\
\langle \{ \s t\}, p \cup \{\s t\} \rangle \in \JLift
                     && \mbox{\qquad by Lemma~\ref{lemma:JLiftRec}}\\
\langle \{ \s x\}, p \cup \{\s x\} \rangle \in \JLift
                     && \mbox{\qquad since $\s t = \s x$}\\
\s x \in \N          && \mbox{\qquad by Lemma~\ref{lemma:successorN}}\\
\s t \in \N           && \mbox{\qquad by Lemma~\ref{lemma:successorN}}\\
p \cup \{\s x\} = y \cup \{\s x\}  && \mbox{\qquad by Lemma~\ref{lemma:JLiftFUNC}}\\
\s x \not\in p        && \mbox{\qquad since $\s t \not \in p$ and $\s x = \s t$}\\
\s x \not \in y       && \mbox{\qquad  proved above }\\
p = y                && \mbox{\qquad by the preceding three lines}
\end{eqnarray*}
Now by Lemma~\ref{lemma:JLiftMaps_helper},
$y$ is closed under successor except $x$, and 
$p$ is closed under successor except $t$.  Explicitly,
\begin{eqnarray*}
\forall u \in y\,(u \neq x \imp \s u \in y) && \mbox{\qquad by Lemma~\ref{lemma:JLiftMaps_helper}}\\
\forall u \in p\,(u \neq t \imp \s u \in p) && \mbox{\qquad by Lemma~\ref{lemma:JLiftMaps_helper}}\\
y \in \FINITE      && \mbox{\qquad by Lemma~\ref{lemma:JLiftfinite}}\\
y \in   \DECIDABLE   &&\mbox{\qquad by Lemma~\sref{lemma:finitedecidable}} \\
x \in y   && \mbox{\qquad by Lemma~\ref{lemma:JLiftMaps_helper}}\\
x \in p   && \mbox{\qquad since $y=p$}\\
t \in y   && \mbox{\qquad since $t \in p$ and $y = p$}\\
\ChurchZero \in y && \mbox{\qquad by Lemma~\ref{lemma:JLiftMaps_helper}}\\
t = x \ \lor \ t \neq x && \mbox{\qquad since $y \in \DECIDABLE$} \\
\end{eqnarray*}
Now I say that $y$ is closed under successor.  To prove that,
suppose $u \in y$.   We have to prove $\s u \in y$.
Since $y \in \DECIDABLE$, we have 
$$ u = x \ \lor u \neq x.$$
We argue by cases accordingly.
\smallskip

Case~1, $u = x$.   Then
\begin{eqnarray*}
x \neq t && \mbox{\qquad by hypothesis}\\
p = y    && \mbox{\qquad proved above}\\
\s x \in y  && \mbox{\qquad since $p$ is 
closed under successor except $t$ and $y=p$} 
\end{eqnarray*}
That completes Case~1.
\smallskip

Case~2, $u \neq x$.  Since $y$ is closed under successor except $x$,
we have $\s u \in y$.  That completes Case~2.
That completes the proof that $y$ is closed under successor.
\smallskip

Now we have proved that $y$ contains $\ChurchZero$ and is closed 
under successor.  By definition of $\N$, we have $\N \subseteq y$.
By Lemma~\ref{lemma:JLiftfinite}, we have $ y \subseteq \N$.  Hence $y = \N$.
But $y \in \FINITE$ and $\N \not \in \FINITE$.  That contradiction
completes the proof of the theorem.

 \begin{corollary} \label{lemma:oneoneimpliesdecidableequality}
If $\N$ is not finite, then $\N$ has decidable equality.
\end{corollary}

\noindent{\em Proof}.
 We prove by induction on $x$ 
that 
\begin{eqnarray}
 \forall y \in \N\,( x = y \ \lor \ x \neq y).\label{eq:E3930}
 \end{eqnarray}
The base case is Lemma~\ref{lemma:decidable0}.
For the induction step, the induction hypothesis is (\ref{eq:E3930}).
We have to prove $\s x = y \ \lor \ \s x \neq y$.
Let $y \in \N $ be given.
By Lemma~\ref{lemma:decidable0}, $y = 0 \ \lor \ y \neq 0$.
If $y = 0$ then we are done by Lemma~\ref{lemma:decidable0}
(or by Theorem~\ref{theorem:successoromitszero}), 
so we may assume $y \neq 0$.  Then by Lemma~\ref{lemma:predecessor},
$y = \s z$ for some $z \in \N$.
Then we have to prove $\s x = \s z \ \lor \ \s x \neq \s z$.
By (\ref{eq:E3930}), we have
$x = z \ \lor x \neq z$.
We argue by cases.
\smallskip

Case~1, $x = z$.  Then $\s x = \s z$.
That completes Case~1.
\smallskip

Case~2, $x \neq z$.  Then by Theorem~\ref{theorem:Churchsuccessorweaklyoneone},
$\s x \neq \s z$.  That completes Case~2.
 That completes the 
proof of the lemma.

\begin{theorem} \label{theorem:Churchsuccessoroneone}
If $\N$ is not finite, then Church successor is one-to-one.
\end{theorem}

\noindent{\em Proof}.  Suppose $\N$ is not finite.  Suppose 
$x\in \N$ and $y \in \N$, and $\s x = \s y$.  We must show $x = y$.
By Theorem~\ref{theorem:Churchsuccessorweaklyoneone}, 
we have $\neg\neg\,(x = y)$.
By Corollary~\ref{lemma:oneoneimpliesdecidableequality},
$\N$ has decidable equality.  Therefore $x = y$.
That completes the proof. 

\begin{theorem} \label{theorem:notfiniteimpliesinfinite}
If $\N$ is not finite, then $\N$ is infinite.
\end{theorem}

\noindent{\em Proof}. Suppose $\N \not\in \FINITE$. 
  Define $Z:= \N - \{ \ChurchZero\}$.  By
Theorem~\ref{theorem:successoromitszero}, $Z$ is a proper subset of $\N$.
By Lemma~\ref{lemma:predecessor}, Church successor maps $\N$ onto $Z$.
By Theorem~\ref{theorem:Churchsuccessoroneone}, Church successor
is one-to-one; therefore it is a similarity from $\N$ to $Z$.
One can verify straightforwardly that Church successor is a similarlity
between $\N$ and $Z$ (120 steps omitted here).    By the 
definition of infinite, $\N$ is infinite.  That completes
the proof of the theorem.

\begin{theorem} \label{theorem:infinity1} The Church counting axiom 
implies $\N$ is infinite and Church successor is one-to-one.
\end{theorem}

\noindent{\em Proof}. By Theorem~\ref{theorem:main}, the 
Church counting axiom implies $\N$ is not finite. Then by 
Theorem~\ref{theorem:notfiniteimpliesinfinite}, $\N$ is infinite.
By Theorem~\ref{theorem:Churchsuccessoroneone}, Church successor
is one-to-one.  That completes the proof of the theorem.

\begin{theorem} \label{theorem:interpretability} 
Heyting's arithmetic HA  can be interpreted in INF plus the
Church counting axiom.
\end{theorem}

\noindent{\em Remark}. Since this is a meta-theorem,  not a theorem of INF,
we have not checked it in Lean as we did all the other proofs in this paper.
\medskip

\noindent{\em Proof}.  We have already shown in \cite{inf-basics} that
one may conservatively add comprehension terms to INF.  We have defined
Church successor, Church addition, and Church multiplication by such 
comprehension terms.  The interpretation of a formula $A$ of HA is defined by 
replacing the function terms of $A$ by comprehension terms involving
the symbols for Church successor, addition, and multiplication, and the 
constant $0$ of HA by the comprehension term defining $\ChurchZero$.  
The quantifiers of $A$ are replaced by bounded quantifiers restricted to $\N$,
which is defined by a comprehension term.  The interpretations of the 
axioms for addition and multiplication hold, by Lemmas~\ref{lemma:ChurchAdditionAssociative},
 \ref{lemma:ChurchAdditionCommutative},
\ref{lemma:ChurchMultiplicationAssociative},  
\ref{lemma:ChurchMultiplicationCommutative},
\ref{lemma:ChurchAddition_equation}, and 
\ref{lemma:ChurchZero_equation}.
The interpretation of any formula of HA is a stratified formula, giving all 
the variables index 0,  so the interpretation of the induction axiom schema
follows from the definition of $\N$.  The main point of 
interest is that we need Theorem~\ref{theorem:Churchsuccessoroneone}  
to verify that
successor is one-to-one.   That completes the proof of the theorem.

\section{If $\F$ is infinite, so is $\N$}
The plan of this section is to prove that (assuming $\F$ is infinite),
every Church number is the order of a cyclic permutation on some finite set.
Then if $\N$ is finite, we have $\s \n = \s k$ with $\n = \k + \m$ and
$\m \neq \ChurchZero$,  so $\m$ has a predecessor $r$.  Then there is 
a finite set $X$ and cyclic permutation $f$ of $X$ whose order is $r$.  Then 
$\m f x = \s r f x = f(r f x) = f x$.  But by the Annihilation Theorem,
$\m f x = x$, contradiction, since a cyclic permutation is not the identity.
While this proof is conceptually simple, there are many details to supply.

\begin{lemma} \label{lemma:enlarge}
Suppose Frege successor maps $\F$ to $\F$.  Then
for every finite set $X$, not-not there exists $c$ with $c \not\in X$.
\end{lemma}

\noindent{\em Remark}.  That is, every finite set is not-not enlargeable,
as $X \cup \{c\}$ is a finite set properly containing $X$.  We were not
able to eliminate the double negation (which would have simplified the subsequent arguments).
\medskip

\noindent{\em Proof}.  Let $p = Nc(X)$.  Then
\begin{eqnarray*}
p \in \F  && \mbox{\qquad by Lemma~\sref{lemma:finitecardinals3}}\\
p^+ \in \F && \mbox{\qquad since $\F$ is closed under successor, by hypothesis}\\
u \in p^+  && \mbox{\qquad for some $u$, by Lemma~\sref{lemma:cardinalsinhabited}} \\
u \in \FINITE && \mbox{\qquad by Lemma~\sref{lemma:finitecardinals2}}
\end{eqnarray*}
 By the definition of Frege successor, 
there exist $v$ and $c$ such that
$u = v \cup \{c\}$, $v \in p$, and $c \not \in v$.  Then
\begin{eqnarray*}
X \in p && \mbox{\qquad by Lemma~\sref{lemma:xinNcx}}\\
X \sim v &&\mbox{\qquad  by Lemma~\sref{lemma:finitecardinals2}}
\end{eqnarray*}

I say $u \not \subseteq X$.  To prove that, assume $u \subseteq X$.
Then 
\begin{eqnarray*}
c \in u && \mbox{\qquad since $u = v \cup \{c\}$}\\
c \in X && \mbox{\qquad since $u \subseteq X$}\\
v \neq X && \mbox{\qquad since $c \not \in v$ and $c \in X$} \\
v \subseteq X && \mbox{\qquad since $v \subseteq u$ and $u \subseteq X$}\\
X \mbox{\ is infinite} && \mbox{\qquad by Definition~\sref{definition:infinite}} \\
X \in \FINITE && \mbox{\qquad by hypothesis}  
\end{eqnarray*}
But then $X$ is both finite and infinite, contradicting Theorem~\sref{theorem:infiniteimpliesnotfinite}.  That
completes the proof that $u \not\subseteq X$. 
\smallskip

By Lemma~\sref{lemma:finitecardinals3}, since $u \in p^+$ and $p^+ \in \F$,
we have $u \in \FINITE$.  
Recall that we can move double negation both ways across a finite universal quantifier:
\begin{eqnarray*}
\forall z \in u\, \neg\neg\, (z \in X) \imp \neg\neg\, \forall z \in u (z \in X)
          && \mbox{\qquad by Lemma~\sref{lemma:finiteDNS}} 
\end{eqnarray*}
and double negation moves in the other direction by pure logic.
Using these facts we have 
\begin{eqnarray*}
\neg\, \forall z \in u\,(z \in X) && \mbox{\qquad by the definition of $\subseteq$}\\
\neg\, \forall z \in u\, \neg\neg\,(z \in X) && \mbox{\qquad by Lemma~\sref{lemma:finiteDNS}, since $u$ and $X$ are finite}\\
\neg\neg\, \exists z \in u \, \neg\, (z \in X) && \mbox{\qquad by logic}
\end{eqnarray*}
That completes the proof of the lemma.

For convenience we repeat Definition~\ref{definition:permutation2}.
\begin{definition}  
$f$ is a {\bf permutation} of a finite set $X$ if and only if $f:X \to X$,
and $Rel(f)$ and $f \in \FUNC$, and $dom(f) \subseteq X$,
and $f$ is both one-to-one and onto from $X$ to $X$.
\end{definition} 

\begin{definition}\label{definition:cyclic}
$f$ is a {\bf cyclic permutation} of a finite set $X$ with {\bf generator} $a$ if 
$f$ is a permutation of $X$ and  $a \in X$ and 
$$\forall z \in X\, \exists r \in \N\, (z = rfa).$$
\end{definition} 

\begin{definition} \label{lemma:permorder} Let $X$ be a finite set and 
let $f:X \to X$ be a cyclic permutation of $X$ with generator $a$.
Suppose $q \in \N$ and $qfa = a$   and $q$ is the $\preceq$-least such 
Church number, i.e.,
$$ \forall r \in \N\,(r \preceq q \imp (\forall x\, (rfx = x) ) \imp r = q).$$
Suppose also that 
$$ \forall x \in X\, \exists r \in \N\,(r \preceq q \ \land \ r f a = x).$$
Then $q$ is the {\bf order} of $f$.
\end{definition}

\noindent{\em Remark}. The second condition is usually omitted, and could also
be omitted in this context, but then we would have to prove that when we divide $x$
by $y$, the remainder is $\prec y$.  That can be done, but it is easier to 
just use this stronger definition of ``order.''
\medskip

\begin{lemma} \label{lemma:permorder2} Let $X$ be a finite set and 
let $f:X \to X$ be a cyclic permutation of $X$  with generator $a$.
Let $q \in \N$ and suppose $qfa = a$.  
Then $\forall x \in X\, (qfx = x)$.
\end{lemma}

\noindent{\em Proof}. Let $x \in X$.
Since $f$ is a cyclic permutation with generator $a$,
we have $x = jfa$ for some $j \in \N$.  Suppose $qfa = a$.  
Then 
\begin{eqnarray*}
qfx = qf(jfa)  && \mbox{\qquad since $x = jfa$}\\
 = (q\oplus j)fa   && \mbox{\qquad by Lemma~\ref{lemma:doubleiteration}}\\
 = (j\oplus q)fa  && \mbox{\qquad since  $\oplus$ is commutative}\\
 =jf(qfa)  &&  \mbox{\qquad by Lemma~\ref{lemma:doubleiteration}}\\
 = jfa  && \mbox{\qquad since $qfa = a$}\\
 = x    &&  \mbox{\qquad since $x = jfa$}
\end{eqnarray*}

\begin{lemma} \label{lemma:orderstep}  Suppose $X$ is a finite set,
$f:X\to X$ is a cyclic permutation of $X$ with generator $a$, and
$c \not\in X$.  Suppose $\N$ is finite with $\s \k = \s \n$ and $k \in \Stem$.
Let $q$ be the order of $f$.  Suppose $q \neq \n$ and $q \neq \ChurchZero$.  
Then there is a cyclic permutation $g$ of $X \cup \{c\}$ with generator $a$ and
order $\s q$.
\end{lemma}

\noindent{\em Proof}.  Let $a$ be a generator of $f$.   
Let $b = f^{-1}(a)$, so $a = f(b)$.   Define
$$ g = f - \{\langle b,a \rangle\} \cup \{\langle b,c \rangle \} \cup \{ \langle c,a \rangle \}.$$
Using functional notation may clarify the idea of the definition of $g$,
although we are not legally entitled to do so until we prove that $g$ is a function.
$$ g(x) = \threecases {f(x)  \mbox{\qquad if $x \neq b \ \land \ x \neq c$}}
                    { c     \mbox{\ \quad\qquad if $x = b$} }  
                    { a     \mbox{\ \quad\qquad if $x = c$}}
$$
 Then $g$ is a permutation of $X \cup \{c\}$.  The formal proof
of that fact requires about 850 steps,  since there are several cases in 
the definition of ``permutation'', and three cases in the definition of $g$,
making nine cases for each case in the definition of ``permutation''.  We
omit those 850 steps.
\smallskip

  Since $q \neq \ChurchZero$,
 by Lemma~\ref{lemma:predecessornotn} 
there exists $t \in \N$ with 
\begin{eqnarray}
\s t = q \ \land \ t \neq \n. \label{eq:E5781}
\end{eqnarray}
  I say that 
 \begin{eqnarray}
 r \prec t \imp rga = rfa \ \land \ rfa \neq b  \label{eq:E5752} \\
 tfa = tga = b  \label{eq:E5846} 
 \end{eqnarray}
In spite of the intuitive conviction that a simple picture of the 
situation provides, the proofs of these assertions are not short;
and unlike the proof that $g$ is a permutation, they are not particularly
straightforward or obvious.  
To keep the 
 length and logical complexity of proofs manageable, we prove  
 these formulas in three separate lemmas.  Logically those lemmas should come before this one, but for 
 readability we postpone them.  We will first finish this proof
under the additional assumptions  (\ref{eq:E5752}) and (\ref{eq:E5846}).
\smallskip

We then have 
\begin{eqnarray}
tfa \in X && \mbox{\qquad by Lemma~\ref{lemma:xfmaps}}\nonumber\\
f (tfa)= f (b)  && \mbox{\qquad by (\ref{eq:E5846})}\nonumber \\
\s tfa =  f(b)  && \mbox{\qquad by Theorem~\ref{theorem:successorequation}}\nonumber\\
qfa = f(b)     && \mbox{\qquad since $\s t = q$}\nonumber\\
qfa = a         && \mbox{\qquad since $f(b) = a$}  \label{eq:qfa}
\end{eqnarray}
Similarly,
\begin{eqnarray}
tga \in X \cup \{c\} && \mbox{\qquad by Lemma~\ref{lemma:xfmaps}}\nonumber\\
g(tga) = g(b) && \mbox{\qquad by (\ref{eq:E5846})} \nonumber\\
\s tga =  g(b)  && \mbox{\qquad by Theorem~\ref{theorem:successorequation}}\nonumber\\
qga = g(b)     && \mbox{\qquad since $\s t = q$}\nonumber\\
qga = c         && \mbox{\qquad since $g(b) = c$} \label{eq:qga}\\
g(qga) = g(c)   && \mbox{\qquad by the preceding line} \nonumber\\
\s q ga = g(c)    && \mbox{\qquad by Theorem~\ref{theorem:successorequation}}\nonumber\\
\s q ga = a       && \mbox{\qquad since $g(c) = a$} \label{eq:sqga}
 \end{eqnarray}
 
We must prove that $g$ is a cyclic permutation of $X \cup \{c\}$
 with generator $a$ and order $\s q$.  Let $z \in X \cup \{c\}$.  We must 
 show that $z = rga$ for some $r \in \N$ with $r \preceq \s q$.
 We have
 
 \begin{eqnarray*}
 X \in \FINITE  && \mbox{\qquad by hypothesis} \\
 X \cup \{c\} \in \FINITE && \mbox{\qquad by Lemma~\ref{lemma:finite_adjoin}}\\
 X \cup \{c\} \in \DECIDABLE && \mbox{\qquad by Lemma~\ref{lemma:finitedecidable}}\\
 z = a \ \lor \ z \neq a && \mbox{\qquad by definition of $\DECIDABLE$} \\
 z = c \  \lor \ z \neq c && \mbox{\qquad by definition of $\DECIDABLE$} 
 \end{eqnarray*}
 If $z=a$ we take $r = \ChurchZero$; then $rga=a=z$.  
 If $z=c$ we take $r = q$; then by (\ref{eq:qga}) we have $rga = qga = c = z$.
  Therefore 
 we may assume 
 \begin{eqnarray}
 z \neq a \ \land \ z \neq c \label{eq:E6855}
 \end{eqnarray}
 Since $z \in X \cup \{c\}$ and $z \neq c$, we have $z \in X$.
 Since $f$ is a cyclic permutation of $X$ of order $q$, there 
 exists $r \preceq q$ such that $z =  rfa$.  By Theorem~\ref{theorem:prectrichotomy1}, 
 we have $r \prec t \ \lor\ r = t \ \lor \ t \prec r$.  I say first that 
\begin{eqnarray}
r \neq t  \label{eq:E5865}
\end{eqnarray}

 Case~1, $r \prec t$.  Then  
\begin{eqnarray*}
rga = rfa \ \land \ rfa \neq b &&\mbox{\qquad by (\ref{eq:E5752})}\\
z = rga  && \mbox{\qquad since $z = rfa$ and $rga = rfa$} \\
q \prec \s q  && \mbox{\qquad by Lemma~\ref{lemma:xpreceqsx}}\\
r \preceq \s q  && \mbox{\qquad by Lemma~\ref{lemma:preceqtrans}, since $r \preceq q$}
\end{eqnarray*}
That completes Case~1.
\smallskip

Case~2, $r=t$.  
Then
 \begin{eqnarray*}
 z = tfa  && \mbox{\qquad since $z = rfa$ and $r=t$}\\
  = tga  && \mbox{\qquad by (\ref{eq:E5846}) }\\
  = rga  && \mbox{\qquad since $r=t$}
\end{eqnarray*}
It remains to show $r \preceq \s q$. We have
\begin{eqnarray*}
q \neq \n && \mbox{\qquad by hypothesis}\\
\s r \neq \n && \mbox{\qquad since $q = \s t$ and $r = t$}\\
r \neq \n && \mbox{\qquad by (\ref{eq:E5781}), since $r=t$}\\
\s r \preceq \s \s r && \mbox{\qquad by Lemma~\ref{lemma:xpreceqsx}, since $\s r \neq \n$ }\\
r \preceq \s r && \mbox{\qquad by Lemma~\ref{lemma:xpreceqsx}, since $ r \neq \n$ }\\
r \preceq \s \s r && \mbox{\qquad by Lemma~\ref{lemma:preceqtrans}}\\
r \preceq \s q  && \mbox{\qquad since $q = \s t$ and $r = t$}
\end{eqnarray*}
That completes Case~2.
\smallskip

 Case~3, $t \prec r$. 
Then 
\begin{eqnarray*}
t \neq r \ \land \ t \preceq r && \mbox{\qquad by definition of $\prec$}\\
r \preceq q = \s t && \mbox{\qquad proved above}\\
r \preceq t \ \lor \ r = \s t && \mbox{\qquad by Lemma~\ref{lemma:preceqsuccessor}, since $t \neq \n$}\\
\neg\, (r\preceq t) && \mbox{\qquad by Theorem~\ref{theorem:prectrichotomy3}, since 
$t \preceq r$ and $t \neq r$} \\
r = \s t  && \mbox{\qquad by the preceding two lines}\\
tfa = tga && \mbox{\qquad by (\ref{eq:E5846})}\\
r = q  && \mbox{\qquad since $q = \s t$} \\
z = qfa  && \mbox{\qquad since $z = rfa$ and $r = q$}\\
\exists r \preceq q\, (z = rfa) && \mbox{\qquad namely, $r = q$}
\end{eqnarray*}
That completes Case~3.
\smallskip

I say that
$\s q$ is the order of $g$. 
By (\ref{eq:sqga}),   $\s q g a= a$.   Suppose $r g a = a$.  We must 
prove $\s q \preceq r$.  By Theorem~\ref{theorem:prectrichotomy3},
it suffices to derive a contradiction from $r \prec \s q$.  
We now assume $r \prec \s q$ 

Then by Lemma~\ref{lemma:preceqsuccessor},
$ r \preceq q$.  We have 
\begin{eqnarray*}
qga = c && \mbox{\qquad by (\ref{eq:qga}) }\\
rga = a && \mbox{\qquad by the definition of order}\\
a \neq c  && \mbox{\qquad since $a \in X$ but $c \not\in X$}\\
r \neq q  && \mbox{\qquad since if $r=q$ then $a =c$}\\
\end{eqnarray*}
\begin{eqnarray*}
r \preceq \s t && \mbox{\qquad since $q = \s t$}\\
r \preceq t \ \lor r = \s t && \mbox{\qquad by Lemma~\ref{lemma:preceqsuccessor}}\\
r \neq \s t && \mbox{\qquad since $q = \s t$ and $r \neq q$}\\
r \preceq t && \mbox{\qquad by the preceding two lines}\\
rga= a && \mbox{\qquad since $r$ is the order of $g$}\\
tga = b  && \mbox{\qquad by (\ref{eq:E5846}) }\\
r \neq t  && \mbox{\qquad by the preceding lines, since $a \neq b$}\\
r \prec t && \mbox{\qquad by the definition of $\prec$}\\
rga = rfa && \mbox{\qquad by (\ref{eq:E5752})}\\
rfa = a && \mbox{\qquad since $rga = a$}\\
q \preceq r && \mbox{\qquad since $q$ is the order of $f$}\\
q = r && \mbox{\qquad by Theorem~\ref{theorem:prectrichotomy3}, since $r \preceq q$}
\end{eqnarray*}
That contradicts $r \neq q$, which was proved above. 
We have now proved that if $rga = a$, then $\s q \preceq r$.
That is, however, only have the definition of ``$\s q$ is the order of $g$.''
\smallskip

We still have to prove that every $x \in X \cup \{c\}$ has the 
form $rga$ for some $r \preceq \s q$, which is the second part of the definition.
Let $x \in X \cup \{c\}$.  Since $X \cup \{c\}$ is finite, it has 
decidable equality, so we have 
$$ x = c \ \lor \ x = b \ \lor x = a \ \lor \ (x \neq c \ \land \ x \neq b \ \land \ x \neq a).$$
We argue by cases.
\smallskip

Case~1, $x = c$.  Then take $r = q$. We have $qga = c$
by (\ref{eq:qga}), and $q \preceq \s q$ since $q \neq \n$.
\smallskip

Case~2, $x = b$.  Then take $r = t$.  We have $tga = b$ by (\ref{eq:E5846}),
and $t \preceq \s t \preceq \s (\s t) = \s q$, since $t \neq \n$ and $q \neq \n$.
\smallskip

Case~3, $x = a$.  Then take $r = \s q$.  We have $\s q ga = a$ by (\ref{eq:sqga}).
\smallskip

Case~4, $x \neq c \ \land \ x \neq b \ \land \ x \neq a$.  Then $x \in X$.
By the induction hypothesis, $q$ is the order of $f$, so there exists $\rho \preceq q$
such that $rfa = x$.  We have
\begin{eqnarray*}
\rho \preceq \s t && \mbox{\qquad since $q = \s t$}\\
\rho \preceq t \ \lor \ r = \s t && \mbox{\qquad by Lemma~\ref{lemma:preceqsuccessor}} \\
\rho= t \ \lor \ \rho \neq t && \mbox{\qquad since $\N \in \DECIDABLE$} \\
\rho \prec t \ \lor \ \rho = t \ \lor \ \rho = \s t && \mbox{\qquad by the definition of $\prec$}
\end{eqnarray*}
If $\rho \prec t$ then by  (\ref{eq:E5752}),  $\rho ga = \rho fa = x$, so we can take $r=\rho$.
If $\rho = t$ then by  (\ref{eq:E5846}), then $tga = b$, so $x=b$, 
contradicting the hypothesis of Case~4. If $\rho= \s t = q$, then 
$x = \rho ga = qga = c$ by (\ref{eq:qga}), contrary to the hypothesis $x \neq c$
of Case~4. 
That completes Case~4.  
That completes the proof that $\s q$ is the order of $g$. 
That completes the proof of the lemma,
under the assumptions (\ref{eq:E5752}) and (\ref{eq:E5846}).

 \begin{lemma} \label{lemma:formula62}
  Suppose $X$ is a finite set,
$f:X\to X$ is a cyclic permutation of $X$ with generator $a$, and
$c \not\in X$.  Suppose $\N$ is finite with $\s \k = \s \n$ and $k \in \Stem$.
Let $q$ be the order of $f$.  Suppose $q \neq \n$ and $q = \s t$ with $t \neq \n$,
and 
$$ g = f - \{\langle b,a \rangle\} \cup \{\langle b,c \rangle \} \cup \{ \langle c,a \rangle \}.$$
Then 
$$ r \prec t \imp rga = rfa \ \land \ rfa \neq b.$$
 \end{lemma}
 
\noindent{\em Proof}. 
 The conclusion of the lemma is a stratified formula, giving  $a$ and $b$ index 0, $f$ index 3,
 and $t$ and $r$ index 6.  We prove it by induction on $r$.
 \smallskip
 
 Base case, $r = \ChurchZero$.  Then $rga = rfa = a$.  If $a=b$ then 
 $f(a) = a$, so $f$ is the identity, contradiction.  That completes the 
 base case.
 \smallskip
 
Induction step.  We use ``finite induction'', so we may assume $r \neq \n$. 
\begin{eqnarray*}
\s r \prec t  && \mbox{\qquad assumption} \\
r \prec \s r  && \mbox{\qquad by Lemma~\ref{lemma:xprecsx}, since $r \neq \n$}\\
r \prec t      && \mbox{\qquad by Lemma~\ref{lemma:prectrans}}\\
 rga  = rfa       && \mbox{\qquad by the induction hypothesis} \\
g(rga) = g(rfa)   && \mbox{\qquad by the preceding line} \\
f:X \to X         && \mbox{\qquad by definition of ``permutation''}\\
rf:X \to X         && \mbox{\qquad by Lemma~\ref{lemma:xfmaps}}\\
rfa \in X          && \mbox{\qquad since $a \in X$ and $rf: X \to X$}\\
rfa \neq b        && \mbox{\qquad by the induction hypothesis}\\
g( rfa) = f (rfa)  && \mbox{\qquad by the definition of $g$, since $rfa \neq b$}\\
= \s r fa   && \mbox{\qquad by Theorem~\ref{theorem:successorequation}} \\
\s rga = g (rga)  && \mbox{\qquad by Theorem~\ref{theorem:successorequation}}\\
 = \s rfa && \mbox{\qquad  by the preceding lines} 
\end{eqnarray*}
It remains to show that $\s rfa \neq b$.  We have
\begin{eqnarray*}
\s rfa = b  && \mbox{\qquad by assumption, for proof by contradiction}\\
 f(\s rfa)) = f(b) && \mbox{\qquad by the previous line}\\
 f(\s rfa) = a    && \mbox{\qquad since $f(b) = a$} \\
 (\s \s r) f a = a &&\mbox{\qquad by Theorem~\ref{theorem:successorequation}} \\
q \preceq \s \s r  && \mbox{\qquad since $q$ is the order of $f$}  \\
\s t \preceq \s \s r && \mbox{\qquad since $\s t = q$ }\\
\neg \n \prec t   && \mbox{\qquad by Lemma~\ref{lemma:precmax2}}\\
\s r \prec t      && \mbox{\qquad as assumed above for the induction step}\\
\s r \neq \n       && \mbox{\qquad by the preceding two lines }\\
t \neq \n          && \mbox{\qquad by (\ref{eq:E5781})}\\
t \preceq \s r  && \mbox{\qquad by Lemma~\ref{lemma:precpred}, since $\s t \preceq \s \s r$ and $t\neq \n$ }\\
\s r \preceq t && \mbox{\qquad by definition of $\prec$, since $\s r \prec t$} \\
\s r = t        && \mbox{\qquad by Theorem~\ref{theorem:prectrichotomy3}} \\
\s r \neq t     && \mbox{\qquad by definition of $\prec$, since $\s r \prec t$} 
\end{eqnarray*}
That contradiction completes the proof that $\s rga \neq b$. 
That completes the induction step.  That completes
the proof of the lemma.

\begin{lemma} \label{lemma:formula63}
  Suppose $X$ is a finite set,
$f:X\to X$ is a cyclic permutation of $X$ with generator $a$, and
$c \not\in X$.  Suppose $\N$ is finite with $\s \k = \s \n$ and $k \in \Stem$.
Let $q$ be the order of $f$.  Suppose $q \neq \n$ and $q = \s t$ with $t \neq \n$,
and 
$$ g = f - \{\langle b,a \rangle\} \cup \{\langle b,c \rangle \} \cup \{ \langle c,a \rangle \}.$$
Then $tga = tfa = b$.
 \end{lemma}

\noindent{\em Proof}.
\begin{eqnarray*}
z:= tfa  && \mbox{\qquad definition of $z$}\\
tfa \in X && \mbox{\qquad by Lemma~\ref{lemma:xfmaps}}\\
z \in X  && \mbox{\qquad by the preceding two lines}\\
fz = f(t f a)  && \mbox{\qquad since $z = tfa$}\\
 fz = \s t f a   &&\mbox{\qquad  by Theorem~\ref{theorem:successorequation}}\\
 fz = q f a  &&\mbox{\qquad since $q = \s t$}\\
 fz = a    && \mbox{\qquad since $qfa = a$}  \\
  z = b  && \mbox{\qquad since $f$ is one-to-one and $fb = a$}\\
 t \prec \s t && \mbox{\qquad by Lemma~\ref{lemma:xprecsx}, since $t \neq \n$}\\
 t \prec q   && \mbox{\qquad since $\s t = q$}\\
 tfa \neq a  && \mbox{\qquad since $q$ is the order of $f$}\\
 z \neq a && \mbox{\qquad since $z = tfa$ }\\
  t \neq \ChurchZero && \mbox{\qquad since if $t = \ChurchZero$ then $z = tfa = a$} \\
 t = \s p  && \mbox{\qquad for some $p \in \N$ with $p \neq \n$, by Lemma~\ref{lemma:predecessornotn}} \\
  p \prec t  && \mbox{\qquad by Corollary~\ref{lemma:xprecsx}} \\
   u := p f a && \mbox{\qquad defining $u$}\\
 u = pga && \mbox{\qquad by (\ref{eq:E5752})} \\
gu = \s p g a &&\mbox{\qquad  by Theorem~\ref{theorem:successorequation}}\\
     = tga  && \mbox{\qquad since $t = \s p$} \\
     f u = f(p f a) && \mbox{\qquad since $u = pfa$}\\
 = \s p f a &&\mbox{\qquad  by Theorem~\ref{theorem:successorequation}} \\
 = tfa  && \mbox{\qquad since $\s p = t$}\\
 = z  = b  && \mbox{\qquad as shown above}\\
  f u \neq f b && \mbox{\qquad since $fu = b$ and $fb = a$ and $a \neq b$ }\\
 pf : X \to X && \mbox{\qquad by Lemma~\ref{lemma:xfmaps}} \\
 u \in X  && \mbox{\qquad since $u = pfa$ and $a \in X$ and $pf:X \to X$}\\
 u \neq b  && \mbox{\qquad since $f$ is one-to-one and $f u \neq f b$} \\
 gu = fu   && \mbox{\qquad by definition of $g$, since $u \neq b$}\\
 gu = b    && \mbox{\qquad since $fu = b$}\\
 z = b = t ga && \mbox{\qquad since $gu = tga$}\\
z = rga && \mbox{\qquad since $r=t$} 
\end{eqnarray*} 
That completes the proof of the lemma.
\smallskip

Now we have supplied the supporting lemmas required to prove
(\ref{eq:E5752}) and (\ref{eq:E5846}).  That completes the proof
of Lemma~\ref{lemma:orderstep}.

\begin{lemma}\label{lemma:simplestperm} Let $X = \{ a, b\}$, where $a \neq b$.
Let $f = \{ \langle a, b \rangle, \langle b, a \rangle \}$.  Then $f$ is 
a cyclic permutation with generator $a$ of $X$ with order $\s (\s \ChurchZero)$.  
\end{lemma} 

\noindent{\em Remark}.
There is no cyclic permutation of order $\one = \s \ChurchZero$,
as it would have to fix the generator $a$, so $X$ would have to be a singleton,
and the identity permutation on  a singleton has order $\ChurchZero$.
\medskip

\noindent{\em Proof}. 
  We omit the 610 simple steps of this proof. (The definitions involved
create many cases.)

\begin{lemma}\label{lemma:allorders}  
Suppose $\N$ is finite and Frege successor maps $\F$ to $\F$ (so $\F$ 
is infinite). 
Let $q \in \N$ with $q \neq \ChurchZero$ and $q \neq \s \ChurchZero$. Then not-not there exists a 
finite set $X$ and a cyclic permutation $f$ of $X$ with generator $a$ 
and order $q$.
\end{lemma}

\noindent{\em Remark}.  ``Every Church number is not-not the order of some
permutation.''  The double negation comes from Lemma~\ref{lemma:enlarge}.
\smallskip

\noindent{\em Proof}.  By Lemma~\ref{lemma:kinstem}, there exist $\k$ and $\n$
in $\N$ 
 with $\s \k = \s \n$ and $\k \in \Stem$ and $\n \neq \k$.
 By   induction on $q$ we prove that $q \neq \ChurchZero$ implies  not-not $q$ is the 
order of some cyclic permutation on a finite set.   We use ``finite induction'',
which means we get to assume $q \neq \n$ in the induction step.
\smallskip

Base case, there is nothing to prove because of the hypothesis $q \neq \ChurchZero$.
\smallskip

Induction step.  The assumptions for the induction step are 
\begin{eqnarray*}
 q \in \N   &&\\
 \s q \neq \ChurchZero &&\\
 \s q \neq \s \ChurchZero &&\\
 q \neq \n && \mbox{\qquad for ``finite induction''}
 \end{eqnarray*}
 Then also
 \begin{eqnarray*}
 q \neq \ChurchZero   && \mbox{\qquad since $\s q \neq \s \ChurchZero$}
 \end{eqnarray*} 

 We have to produce a set $X$ and a permutation 
$f$ of $X$ whose order is $\s q$.
 If $\s q = \s (\s \ChurchZero)$, or for short $\s q = 2$,  then by 
Lemma~\ref{lemma:simplestperm}, there is a permutation of order 2, so 
we are finished.   In other words, we may assume
\begin{eqnarray*}
q \neq \s \ChurchZero && \mbox{\qquad by Lemma~\ref{lemma:simplestperm}}
\end{eqnarray*}
Now that we have $q \neq \ChurchZero$ and $q \neq \s \ChurchZero$,
we may apply the induction hypothesis.
 By the induction 
hypothesis, not-not there is a finite set $X$ and a cyclic permutation of $X$ 
with generator $a$ such that $xfa = a$, and $f$ has order $q$.
   Suppose $X$ is such a set; 
by Lemma~\ref{lemma:enlarge},  not-not there exists $c \not\in X$. 
By 
Lemma~\ref{lemma:orderstep},  not-not there is an $X$, 
and a cyclic permutation of $X \cup \{c\}$, such that $\s x f a = a$,
and $f$ has order $\s q$.
  That completes the induction step.  That completes
the proof of the lemma.

\begin{theorem}\label{theorem:FinfiniteimpliesNinfinite}
Suppose Frege successor maps $\F$ to $\F$ (so $\F$ is infinite). 
Then $\N$ is infinite and Church successor is one-to-one. 
\end{theorem}

\noindent{\em Proof}.  Suppose Frege successor maps $\F$ to $\F$.
We have to prove Church successor is one-to-one.  By 
Theorem~\ref{theorem:Churchsuccessoroneone}, it suffices to 
prove that $\N$ is not finite.  Suppose $\N$ is finite; we must
derive a contradiction.   Since we have assumed $\N$ is finite, 
by Lemma~\ref{lemma:kinstem}  there exist $\k$ and $\n$
in $\N$ 
 with $\s \k = \s \n$ and $\k \in \Stem$ and $\n \neq \k$.  By Lemma~\ref{lemma:mexists}, there exists
$\m$ such that $\n = \k + \m$.  Since $\n \neq \k$, we have 
$\m \neq \ChurchZero$.  By Lemma~\ref{lemma:predecessor}, there exists
$r\in \N$ with $\s r = \m$.  We have $r \neq \ChurchZero$ and $r \neq \s \ChurchZero$,
by Lemma~\ref{lemma:mbig}.   By
 Lemma~\ref{lemma:allorders},  not-not there is a finite set $X$ and a permutation
 $f:X \to X$ with generator $a$ such that $rfa=a$ and $fa \neq a$.   Suppose $X$ and $f$
 are such a finite set and permutation.      Then 
\begin{eqnarray*}
\m f a = \s r f a  && \mbox{\qquad since $\m = \s r$}\\
= f (r f a)      && \mbox{\qquad by Theorem~\ref{theorem:successorequation}}\\
= f(a)           && \mbox{\qquad since $rfx = x$}
\end{eqnarray*}
On the other hand, by the Annihilation Theorem, $\m f a = a$. 
Hence $f(a) = a$, contradiction.
\smallskip

We reached that contradiction under the assumption that $X$ is a finite set 
with a permutation $f$ of order $r$;  but we have actually proved only the 
double negation of that.  But still,  if a proposition leads to a contradiction,
so does its double negation. 
That completes the proof of the theorem.

\begin{corollary} \label{corollary:Ninfinite}
Classical NF proves $\N$ is infinite and Church successor is 
one-to-one.
\end{corollary}

\noindent{\em Proof}.  Specker \cite{specker1953} proved, in 
classical NF,  that $\F$ is infinite.  By Theorem~\ref{theorem:FinfiniteimpliesNinfinite}, Specker's result
implies that $\N$ is infinite and Church successor is one-to-one.
That completes the proof of the corollary.
\medskip

{\em Remarks}.  We also proved that every inhabited finite set 
has a cyclic permutation; this proof uses the lemmas above in its 
induction step,  but it also requires proving that the order of a 
permutation cannot be $\n$, which needs the Annihilation theorem
and a bit more.  That focuses attention on the question whether an 
unenlargeable finite set $U$ can have a cyclical permutation $f$, 
since if that were impossible, then $\F$ could not be finite,  so 
both $\F$ and $\N$ would be infinite.   However, we could not 
derive a contradiction from the assumption that there is a cyclic 
permutation $f$ on an unenlargeable set $U$.

\section{Equivalence of Church and Rosser counting axioms}

Thanks are due to Thomas Forster, who overcame 
my initial skepticism about proving the equivalence of these two axioms.
\smallskip

We make use of $\T x$, defined so that for $x \in \F$, and $a \in x$,
we have $\T(x) = Nc(USC(a))$.

\begin{lemma} \label{lemma:T2} Suppose that $\T^2 x  = x$ for all $x \in \F$.
  Then $\T x = x$ holds for all $x \in \F$.
\end{lemma}

\noindent{\em Proof}.   By Lemma~\sref{lemma:Tfinite},
$\T x \in \F$.    
By Theorem~\sref{theorem:finitetrichotomy}, we have 
$$x < \T x \ \lor \ x = \T x \ \lor \ \T x < x.$$ 
We argue by cases.
\smallskip

Case~1, $x < \T x$.  Then 
\begin{eqnarray*}
\T x < \T^2 x && \mbox{\qquad by Lemma~\sref{lemma:Tlessthan}}\\
\T x < x      && \mbox{\qquad since $\T^2 x = x$}\\
x  < x  && \mbox{\qquad by Lemma~\sref{lemma:lessthan_transitive}}\\
x \not < x && \mbox{\qquad by Lemma ~\sref{lemma:xnotlessthanx}}
\end{eqnarray*}
That contradiction completes Case~1.
\smallskip

Case~2, $x = \T x$.  Then $\T x = x$ and we are done.
\smallskip

Case~3, $\T x < x$.  Then
\begin{eqnarray*}
\T^2 x < \T x && \mbox{\qquad by Lemma~\sref{lemma:Tlessthan}}\\
x < \T x && \mbox{\qquad since $\T^2 x = x$}\\
x  < x  && \mbox{\qquad by Lemma~\sref{lemma:lessthan_transitive}}\\
x \not < x && \mbox{\qquad by Lemma ~\sref{lemma:xnotlessthanx}}
\end{eqnarray*}
That contradiction completes Case~3.
That completes the proof of the lemma.

\begin{lemma} \label{lemma:T6} Suppose that $\T^6 x  = x$ for all $x \in \F$.
  Then $\T x = x$ holds for all $x \in \F$.
\end{lemma}

\noindent{\em Proof}.  (Similar to the proof of Lemma~\ref{lemma:T2}.)
 By Lemma~\sref{lemma:Tfinite},
$\T x \in \F$.    
By Theorem~\sref{theorem:finitetrichotomy}, we have 
$$x < \T x \ \lor \ x = \T x \ \lor \ \T x < x.$$ 
We argue by cases.
\smallskip

Case~1, $x < \T x$.  Then 
\begin{eqnarray*}
x < \T x < \T^2 x < T^3 x \ldots < \T^6 x && \mbox{\qquad by Lemma~\sref{lemma:Tlessthan}}\\
x  < \T^6 x && \mbox{\qquad by Lemma~\sref{lemma:lessthan_transitive}}\\
x < x  && \mbox{\qquad since $\T^6 x = x$}\\
x \not < x && \mbox{\qquad by Lemma ~\sref{lemma:xnotlessthanx}}
\end{eqnarray*}
That contradiction completes Case~1.
\smallskip

Case~2, $x = \T x$.  Then $\T x = x$ and we are done.
\smallskip

Case~3, $\T x < x$.  Then
\begin{eqnarray*}
\T^6 x < \T^5 x \ldots < \T^2 x < \T x && \mbox{\qquad by Lemma~\sref{lemma:Tlessthan}}\\
\T^6 x < x && \mbox{\qquad by Lemma~\sref{lemma:lessthan_transitive}}\\
x  < x  && \mbox{\qquad since $\T^6 x = x$}\\
x \not < x && \mbox{\qquad by Lemma ~\sref{lemma:xnotlessthanx}}
\end{eqnarray*}
That contradiction completes Case~3.
That completes the proof of the lemma.

Define $J(x):= \{z \in \F: z < x\}.$   Rosser stated his counting 
axiom in the form $J(x) \in x$ for $x \in \F$.  

\begin{lemma} \label{lemma:RosserT}
 Rosser's counting axiom is equivalent to $\T(x) = x$ for all $x \in \F$.
\end{lemma}

\noindent{\em Proof}.
Left to right:  Assume Rosser's counting axiom.  Let $x \in \F$.  Then
\begin{eqnarray*}
Nc(J(x)) = \T^2 x  && \mbox{\qquad by Lemma~\sref{lemma:Jcardinality}}\\
J(x) \in x   && \mbox{\qquad by   Rosser's counting axiom}\\
Nc(J(x)) = x  && \mbox{\qquad by Lemma~\sref{lemma:Ncdef}}\\
\T^2 x = x  && \mbox{\qquad by the preceding lines}
\end{eqnarray*}
Since $x$ was arbitrary, we have proved 
$\forall x \in \F\,(  \T^2 x = x)$.  Then  by Lemma~\ref{lemma:T2},
we have $\forall x \in \F\, (\T x = x)$.  That completes the left-to-right
direction.

Right to left: Assume for all $x \in \F$,  $\T(x) = x$.  Then 
$Nc(J(x)) = \T^2 x = \T x = x$, so $\J(x) \in x$ by Lemma~\sref{lemma:finitecardinals3}.
That completes the proof of the lemma.

\begin{definition} \label{definition:ChurchFrege}
 Let $\i$   be the intersection of all sets $w$
such that 
\begin{eqnarray*}
&&u \in w \imp \exists p,q\, (u = \langle p,q\rangle \ \land \ p \in \N \ \land \ q \in \F  \\
&&\langle \ChurchZero, \zero \rangle \in w \\
&& \langle p,q \rangle \in w \imp q^+ \in \F \imp
 \langle \s p, q^+ \rangle \in w.
\end{eqnarray*}
\end{definition}

\noindent{\em Remark}.  If $\N$ is finite,  then $\i$ may not be a function,
because if $\s \n = \s \k$ then $\i(\s \n)$ would have to be both 
$\i(\n)^+$ and $\i(\k)^+$.   And, if $\F$ is not finite, then the domain of
$\i$ might be a proper subset of $\F$.  But, if both $\N$ and $\F$ are infinite,
then $\i$ should turn out to be a similarity between them.

\begin{lemma} \label{lemma:ChurchFrege0} $\langle \ChurchZero, \zero \rangle \in \i$.
\end{lemma}

\noindent{\em Proof}.  Let $w$ be a set satisfying the conditions in 
the definition of $\i$.  Then $\langle \ChurchZero, \zero \rangle \in w$.
Since $w$ is arbitrary,  $\langle \ChurchZero, \zero \rangle \in \i$.
That completes the proof of the lemma.

\begin{lemma} \label{lemma:ChurchFrege1}
Suppose $\langle p,q \rangle \in \i$ and $q^+ \in \F$. Then $\langle \s p, q^+\rangle \in \i$. 
\end{lemma}

\noindent{\em Proof}. Suppose 
$\langle p, q\rangle \in \i$ and $q^+ \in \F$.  We must show $\langle \s p, q^+\rangle
\in \i$.  Let $w$ be any set satisfying the conditions in the definition.
Then $\langle p,q \rangle \in w$.  Since $q^+ \in \F$ we have 
$\langle \s p, q^+ \rangle \in w$.  Since $w$ was any set satisfying the 
conditions, we have $\langle \s p, q^+\rangle \in \i$ as desired.
That completes the proof of the lemma.

\begin{lemma} \label{lemma:ChurchFrege2}  
 Let $\i$ be as in Definition~\ref{definition:ChurchFrege}.  Suppose $\langle p, q\rangle\in \i$.
 Then either $p = \ChurchZero \ \land \ q = \zero$,  or for some $t \in\N$ and $r \in \F$,  we have $$\langle t,r \rangle \in \i  \ \land \ 
  p = \s t \ \land \ q = r^+.$$
\end{lemma}

\noindent{\em Proof}.   Define
$$ w := \{\, \langle p,q \rangle \in \i :  (p = \ChurchZero \ \land \ q = \zero) 
\ \lor \ \exists t,r\,(\langle t, r \rangle \in \i \ \land \ p = \s t \ \land \ q = r^+)\, \}.
$$
The formula is stratified, giving $p,q,r$ all index 0; $\i$ is a parameter.
Therefore the definition is legal.
\smallskip

Then $w$ satisfies the conditions in Definition~\ref{definition:ChurchFrege}.  Therefore
$\i \subseteq w$.   That completes the proof of the lemma. 

\begin{lemma} \label{lemma:ChurchFrege3}  
 Let $\i$ be as in Definition~\ref{definition:ChurchFrege}.  Suppose $\langle p, y^+\rangle\in \i$.
 Then   for some $t \in\N$ ,  we have $$\langle t,y \rangle \in \i  \ \land \ 
  p = \s t.$$
\end{lemma}

\noindent{\em Proof}.   Take $q$ in Lemma~\ref{lemma:ChurchFrege2} to be 
$y^+$.  By Lemma~\sref{lemma:Fregesuccessoromits0}, we do not have $y^+ = \zero$.
Hence there exists $t,r$ with $\langle t,r \rangle \in \i$ and 
$p = \s t$ and $y^+ = r^+$.  By Lemma~\sref{lemma:successoroneone}, we have 
$y = r$.  Then $\langle t, y \rangle \in i$.  That completes the proof of the lemma.

\begin{lemma} \label{lemma:ChurchRosserhelper} Suppose Church successor 
is one-to-one on $\N$.  Then for all $x \in \N$ and $z \in \F$, 
 $$
\langle x, y \rangle \in \i \imp \langle x \s \ChurchZero,z \rangle \in \i 
\imp \T^6 z = y.$$ 
\end{lemma}

\noindent{\em Proof}.
The formula is stratified, giving $x$ and $y$ index 6 and $z$ index 0,
with $\i$ as a parameter.  Therefore we may prove it by $\N$-induction on $x$.
\smallskip

Base case, $x = \ChurchZero$.  Assume 
$\langle \ChurchZero, y \rangle \in \i $ and $\langle \ChurchZero \s \ChurchZero,z \rangle \in \i$.
Then $y = z = \zero$, since $\ChurchZero \s \ChurchZero = \ChurchZero$. 
By Lemma~\sref{lemma:Tzero}, $\T^6 \ChurchZero = \ChurchZero$.
That completes the base case.
\smallskip

Induction step.  Suppose $\langle \s x, y \rangle \in \i$ and 
$\langle \s x \s \ChurchZero, z \rangle \in \i$.  
Then
\begin{eqnarray*}
\langle \s (x \s \ChurchZero), z \rangle \in \i  && \mbox{\qquad by Theorem~\ref{theorem:successorequation}}
\end{eqnarray*}
By Lemma~\ref{lemma:ChurchFrege3}, 
 $z = r^+$ for some $r \in \F$ and some $p$ with  
  $\langle p, r \rangle \in \i$ where $\s p = \s (x \s \ChurchZero)$.
Since Church successor is one-to-one (by hypothesis), we have 
$p = x \s \ChurchZero$.  Therefore $\langle x \s \ChurchZero,r \rangle\in \i$.
 Similarly, since  $\langle \s x, y \rangle \in \i$,
  $y = t^+$ for some $t \in \F$ with $\langle x,t\rangle \in \i$.
Since  $\langle x,t\rangle \in \i$ and $\langle x \s \ChurchZero,r \rangle\in \i$,
we have 
\begin{eqnarray*}
\T^6 r = t  && \mbox{\qquad  by the induction hypothesis}\\
(\T^6 r)^+ = t^+ && \mbox{\qquad by the previous line}\\
(\T^6 r)^+  = y && \mbox{\qquad since $y=t^+$}
\end{eqnarray*}
I say that
\begin{eqnarray}
(\T^6 r)^+ = \T^6 (r^+)  \label{eq:E5798}
\end{eqnarray}
We wish to justify (\ref{eq:E5798}) by six applications of
 Lemma~\sref{lemma:Tsuccessor}; let us consider the first step,
$(\T r)^+ = \T (r^+)$.  To use Lemma~\sref{lemma:Tsuccessor}, we 
need to show that $r^+$ is inhabited.   We have
\begin{eqnarray*}
z \in \F && \mbox{\qquad since $\langle \s x \s \ChurchZero, z \rangle \in \i$}\\
\exists u\,(u \in z) && \mbox{\qquad by Lemma~\sref{lemma:cardinalsinhabited}}\\
\exists u\, (u \in r^+)&& \mbox{\qquad since $z = r^+$} \\
(\T r)^+ = \T (r^+) && \mbox{\qquad by Lemma~\sref{lemma:Tsuccessor}}
\end{eqnarray*}
 Now to take the next steps, we 
need Lemma~\sref{lemma:successorT}, which says that $\T$ of anything in $\F$
has its successor in $\F$.   We have
\begin{eqnarray*}
\T r \in \F   && \mbox{\qquad by Lemma~\ref{lemma:Tfinite}}\\
(\T r)^+ \in \F && \mbox{\qquad by Lemma~\ref{lemma:successorT}}\\
\exists u\,(u \in (\T r)^+) && \mbox{\qquad by Lemma~\sref{lemma:cardinalsinhabited}}\\
 \T^2 r \in \F  && \mbox{\qquad by Lemma~\ref{lemma:Tfinite}}\\
(\T^2 r)^+ \in \F  &&\mbox{\qquad by Lemma~\ref{lemma:successorT}}\\
\exists u\, (u \in (\T^2 r)^+) && \mbox{\qquad by Lemma~\sref{lemma:cardinalsinhabited}}
\end{eqnarray*}
Continuing in this way we eventually prove $\T^m r$ is in $\F$ and is inhabited,
for $m$ up to and including 6.  Now we have all the side conditions 
necessary to apply Lemma~\sref{lemma:Tsuccessor} six times.  
That completes the proof of (\ref{eq:E5798}).  Continuing, we have
\begin{eqnarray*}
(\T^6 r)^+  = y   && \mbox{\qquad as shown above}\\
\T^6 (r^+) = y   && \mbox{\qquad by (\ref{eq:E5798}) }\\
 \T^6 z = y   && \mbox{\qquad since $z = r^+$}
 \end{eqnarray*}
   That completes the induction step.  That completes
the proof of the lemma.

\begin{lemma}\label{lemma:ChurchFregeonto}
The relation $\i$ is onto $\F$.   Explicitly,
\begin{eqnarray}
\forall y \in \F\, \exists x \in \N\, (\langle x, y \rangle \in \i) \label{eq:E5803}
\end{eqnarray}
\end{lemma}

\noindent{\em Proof}.  The formula is stratified, giving $x$ and $y$ index 0, 
so it can be proved
by $\F$-induction on $y$. 
\smallskip

Base case, $y = \zero$;  then $x = \ChurchZero$ will do.
\smallskip

Induction step.  Suppose $\langle x, y \rangle \in \i$ and $y^+$ is inhabited.
Then $y^+ \in \F$ by Lemma~\sref{lemma:successorF}, so $\langle \s x, y^+\rangle \in \i$.  That completes the proof of the lemma.

\begin{theorem}\label{theorem:ChurchimpliesRosser}
The Church counting axiom implies the Rosser counting axiom.
\end{theorem}

\noindent{\em Proof}.  Assume the Church counting axiom.  By 
Theorem~\ref{theorem:main}, $\N$ is not finite.  
Then by Theorem~\ref{theorem:infinity1}, Church successor is 
one-to-one on $\N$.  
Let $z \in \F$.  By Lemma~\ref{lemma:ChurchFregeonto},
 there exists $x \in \N$ 
such that $\langle x, z\rangle \in \i$. We have
\begin{eqnarray*}
x \s \ChurchZero = x  && \mbox{\qquad by the Church counting axiom} \\
\langle x \s \ChurchZero,z \rangle \in \i && \mbox{\qquad since $\langle x, z\rangle \in \i$} \\
\T^6 z = z && \mbox{\qquad by Lemma~\ref{lemma:ChurchRosserhelper}, with $y=z$} \\
\forall z\in \F\,(\T^6 z = z)  && \mbox{\qquad since $z$ was arbitrary}\\
\forall z \in \F\,(\T x = x) && \mbox{\qquad by Lemma~\ref{lemma:T6}}
\end{eqnarray*}
By Lemma~\ref{lemma:RosserT}, we have the Rosser counting axiom.
That completes the proof of the theorem.

\begin{lemma} \label{lemma:ChurchFregeoneone}
$\langle x,y \rangle \in \i \imp \langle z,y \rangle \in \i \imp x = z$.
\end{lemma}

\noindent{\em Proof}.  The formula is stratified, giving $x$, $y$, and $z$ 
all index 0.  We prove it by $\F$-induction on $y$.
\smallskip

Base case, $y = \zero$.  By Lemma~\ref{lemma:ChurchFrege2} and 
Lemma~\sref{lemma:Fregesuccessoromits0}, we have $x = \ChurchZero$
and $z = \ChurchZero$.  Therefore $x=z$.  That completes the base case.
\smallskip

Induction step.  Suppose $\langle x, y^+\rangle \in \i$ and 
$\langle z, y^+\rangle \in \i$.  We must show $x=z$.
By Lemma~\ref{lemma:ChurchFrege2} and Lemma~\sref{lemma:Fregesuccessoromits0},
there exist $p,q \in \N$ and $r,t \in \F$ such that 
$x = \s p$, $z = \s q$, and $r^+ = y^+ = t^+$, and $\langle p,r \rangle \in \i$
and $\langle q,t\rangle \in \i$.
  By Lemma~\sref{lemma:successoroneone},
$r = y =t$.   Then 
\begin{eqnarray*}
 \langle q,y\rangle \in \i &&\mbox{\qquad since $t = y$ and $\langle q,t\rangle \in \i$}\\
\langle p,y\rangle \in \i &&\mbox{\qquad since $\langle p,r\rangle \in \i$ and $y=r$}\\
p = q && \mbox{\qquad by the induction hypothesis}\\
x = z && \mbox{\qquad since $x = \s p$, $z = \s q$, and $p=q$}
\end{eqnarray*}
That completes the induction step.  That completes the proof of the lemma.

\begin{lemma} \label{lemma:RosserimpliesFinfinite}
The Rosser counting axiom implies that $x^+$ is inhabited for every $x \in \F$.
\end{lemma}

\noindent{\em Proof}.  Assume the Rosser counting axiom.
Define $J(x):= \{z \in \F: z < x\}$.   Then
  $x \not\in J(x)$, by Lemma~\sref{lemma:xnotlessthanx}. 
We will prove that
 $x^+$ is inhabited for every $x \in \F$.  The formula is 
 stratified, so we may prove it by induction on $x$.
\smallskip

Base case, $\zero^+$ is inhabited, since $\zero \in \zero^+$.
\smallskip

Induction step.  Suppose $x^+$ is inhabited; we must show $x^{++}$ is 
inhabited.  By the Rosser counting axiom, $J(x^+)\in x^+$. 
By Lemma~\sref{lemma:successorF}, $x^+ \in \F$.   
By Lemma~\sref{lemma:xnotlessthanx}, $x^+ \not\in J(x^+)$.
Then $J(x^+) \cup \{ x^+\} \in x^{++}$.  That completes the 
induction step.   That completes the proof of the lemma.

\begin{lemma} \label{lemma:Rosserimpliesitotal}
The Rosser counting axiom implies $\i$ is total.  Explicitly:
\begin{eqnarray*}
\forall x \in \N\, \exists y \in \F\, (\langle x, y \rangle \in \i)
\end{eqnarray*}
\end{lemma}

\noindent{\em Remark}.  $\i$ may still not be a function.
\medskip

\noindent{\em Proof}.  Assume the Rosser counting axiom.   
The formula is stratified, giving $x$ and $y$ index 0, with $\i$ as parameter.
We prove it by $\N$-induction on $x$.  
\smallskip

Base case, $x = \ChurchZero$.  Take $y = \zero$.
\smallskip

Induction step.   
By the induction hypothesis, there exists $y$ such that 
$\langle x, y\rangle \in \i$ and $y\in\F$.  Then
\begin{eqnarray*}
x \in \N \ \land \ y \in \F  && \mbox{\qquad since $\langle x, y\rangle \in \i$ }\\
\exists u\, (u \in y^+) && \mbox{\qquad by Lemma~\ref{lemma:RosserimpliesFinfinite}}\\ 
 y^+ \in \F && \mbox{\qquad by Lemma~\sref{lemma:successorF}}\\
 \langle \s x, y^+\rangle \in \i && \mbox{\qquad by  Lemma~\ref{lemma:ChurchFrege1}}
 \end{eqnarray*}
  That
completes the induction step.  That completes the proof of the lemma.

 \begin{lemma} \label{lemma:jFUNC}
Assume the Rosser counting axiom.  Let $\j$ be the converse relation of $\i$,
that is, $\langle u, v \rangle \in j \iff \langle v, u\rangle \in \i$.  Then
$j: \F \to \N$ is a function.
\end{lemma} 
 
\noindent{\em Proof}.   We have to prove
\begin{eqnarray*}
\forall y, z\,(\langle y,x \rangle \in \i \imp \langle z,x\rangle \in \i \imp y = z)
\end{eqnarray*}
That formula is stratified, giving $x$, $y$, and $z$ all index 0, with $\i$ as a 
parameter.  We prove it by $\F$-induction on $x$.
\smallskip

Base case, $x = \zero$.  Suppose $\langle y, \zero \rangle \in \i$ and 
$\langle z, \zero \rangle \in \i$.  Then $y = z = \ChurchZero$.  That completes
the base case.
\smallskip

Induction step.  Suppose $x^+$ is inhabited and $\langle y, x^+ \rangle \in \i$
and $\langle z, x^+\rangle \in i$.  Then 
\begin{eqnarray*}
y = \s r  \ \land \ \langle r, x\rangle \in i  && \mbox{\qquad for some $r$, by Lemma~\ref{lemma:ChurchFrege3}}\\
z = \s t  \ \land \ \langle t, x\rangle \in i  && \mbox{\qquad for some $r$, by Lemma~\ref{lemma:ChurchFrege3}}\\
r = t  && \mbox{\qquad by the induction hypothesis}\\
y = z  && \mbox{\qquad since $y = \s r$ and $z = \s t$}
\end{eqnarray*}
That completes the induction step.  That completes the proof of the lemma. 

\begin{lemma} \label{lemma:jhelper}
Assume the Rosser counting axiom.  Then for all $x \in \F$, 
$$(\j x)\s \ChurchZero  = \j (\T^6 x).$$
\end{lemma}

\noindent{\em Proof}.  The formula of the lemma is  stratified,
so we may prove it  by $\F$-induction on $x$.
\smallskip

Base case, $x = \zero$.  We have 
\begin{eqnarray*}
\j(\zero) = \ChurchZero && \mbox{\qquad by Lemma~\ref{lemma:ChurchFrege0}}\\
(\j \zero) \s \ChurchZero = \ChurchZero \s \ChurchZero && \mbox{\qquad by the preceding line}\\
= \ChurchZero  && \mbox{\qquad by Lemma~\ref{lemma:ApZero}}
\end{eqnarray*}
On the other hand 
\begin{eqnarray*}
\T^6 \zero = \zero && \mbox{\qquad by Lemma~\sref{lemma:Tzero}}\\
\j (\T^6 \zero) = \j(\zero)  && \mbox{\qquad by the preceding line}\\
 = \ChurchZero  && \mbox{\qquad by Lemma~\ref{lemma:ChurchFrege0}}
\end{eqnarray*} 
Therefore $(\j \zero) \s \ChurchZero = \j(\T^6 \zero)$, since
both sides are equal to $\ChurchZero$.  That 
completes the base case.
\smallskip

Induction step.  We have 
\begin{eqnarray*}
(\j (x^+)) \s \ChurchZero = (\s (\j x)) \s \ChurchZero &&\mbox{\qquad by Lemma~\ref{lemma:ChurchFrege3}} \\
 = \s ((\j x) \s \ChurchZero) && \mbox{\qquad by Theorem~\ref{theorem:successorequation}}\\
 = \s (\j (\T^6 x)) && \mbox{\qquad by the induction hypothesis}\\
\T x = x && \mbox{\qquad by Lemma~\ref{lemma:RosserT} and the Rosser counting axiom}\\
\T^6 x =x && \mbox{\qquad by six applications of the preceding line}\\
(\j (x^+)) \s \ChurchZero  = \s (\j x)   && \mbox{\qquad by the preceding three lines}\\
 = \j (\T^6 (x^+)) && \mbox{\qquad since $\T(x^+) = x^+$ by Lemma~\ref{lemma:RosserT}}
\end{eqnarray*}
That completes the induction step.   
 That completes the proof of the lemma.

\begin{theorem} \label{theorem:RosserimpliesChurch} 
The Rosser counting axiom implies the Church counting axiom.
\end{theorem}

\noindent{\em Proof}. Assume the Rosser counting axiom.
 Let $z \in \N$.  Then 
 \begin{eqnarray*}
z \s \ChurchZero \in \N  && \mbox{\qquad by Lemma~\ref{lemma:xsmapsN}}\\
 \langle z, x \rangle \in \i && \mbox{\qquad for some $x$, by 
  Lemma~\ref{lemma:Rosserimpliesitotal}}\\
 \j x = z && \mbox{\qquad since $\langle x, z \rangle \in \j \iff \langle z,x \rangle \in \i$}\\  
(\j x)\s \ChurchZero  = \j \T^6 x && \mbox{\qquad by Lemma~\ref{lemma:jhelper}}\\
\T^6 x =x  && \mbox {\qquad by  Lemma~\ref{lemma:RosserT} and the Rosser counting axiom}\\
  (\j x)\s \ChurchZero  = \j   x && \mbox{\qquad by the preceding two lines}
\end{eqnarray*}
Since $z = \j x$, we have $z \s \ChurchZero = z$.
Since $z$ was arbitrary, we have proved the Church counting axiom.
That completes the proof of the lemma.

\section{Conclusion}
We have thoroughly studied the Church numbers in NF set theory,
thus relating the logical work of the two great twentieth-century
logicians Quine and Church.  This analysis led to the result that
the set of Church numbers $\N$, if not infinite, must have a 
certain structure under successor, namely a stem $\Stem$ and a loop $\L$, 
connected at the unique double successor.  That situation led
to the Annihilation Theorem:  There is a Church number $\m$ such
that the $\m$-th iterate of any map from any finite set to itself
is the identity.  

This leads to a contradiction if we can show that there is a finite 
set with a permutation of order not divisible by $\m$.  
This we could do only under the additional assumption 
of the Church counting axiom, that $j \s \ChurchZero = j$
for every Church number $j$.

Assuming the Church counting axiom, we can use the Annihilation
Theorem to prove that the set $\N$ of Church numbers
is not finite.
That result is new even in NF with classical reasoning.  That is,
Specker's result shows that $\F$ is infinite, but that was not
known to imply anything about $\N$, even with the aid of the 
Rosser counting axiom. 

We also proved that if $\N$ is not finite, then $\N$ is infinite
and Church successor is injective.
While that result is classically trivial, it is far from obvious
intuitionistically.  That implicaation does not require the 
Church counting axiom.  Together, the two results show that Heyting's
arithmetic HA is interpretable in INF plus the Church counting axiom.
It remains open whether HA is interpretable in INF alone.  

However,  if we are willing to use classical 
logic, then we do not need the Church counting axiom, as we will
now explain.  Specker proved, using classical logic, that $\F$ is 
infinite.  We proved in Theorem~\ref{theorem:FinfiniteimpliesNinfinite} 
that
if $\F$ is infinite, but $\N$ is finite, 
then every Church number
is the order of some permutation, contradicting the Annihilation Theorem.
Therefore, if $\F$ is infinite, $\N$ is not finite, and therefore $\N$
is infinite.   Then, appealing to Specker's result, classical logic implies
$\N$ is infinite.  Summarizing, we showed that $\N$ is infinite if either
the Rosser counting axiom, or the Church counting axiom, or classical logic
holds.   It remains open whether INF proves $\N$ is infinite (without any
additional assumption). 

We also  proved that the Church counting axiom 
is equivalent to the Rosser counting axiom in INF.  This proof uses
our results that Church counting implies $\N$ is infinite.  And again,
although we proved the equivalent intuitionistically,  it is a new result
even classically.


\end{document}